\documentclass[11pt]{amsart}
\usepackage{amsmath,amsthm,amscd,amssymb, rlepsf, verbatim}
\newtheorem{thm}{Theorem}[section]
\newtheorem{lem}[thm]{Lemma}
\newtheorem{cor}[thm]{Corollary}

\newtheorem{lma}[thm]{Lemma}
\newtheorem{rem}[thm]{Remark}
\newenvironment{rmk}{\begin{rem}\rm}{\end{rem}}
\newtheorem*{clm}{Claim}

\newenvironment{pf}{\begin{proof}}{\end{proof}}

\theoremstyle{definition}
\newtheorem{dfn}[thm]{Definition}
\theoremstyle{remark}
\newtheorem{remark}[thm]{Remark}

\numberwithin{equation}{section}
\newcommand{\bbr}{\begin{remark}}        
\newcommand{\eer}{\end{remark}}

\font\bbb=msbm10 scaled 1100
\newcommand{\bea}{\begin{eqnarray}}
\newcommand{\eea}{\end{eqnarray}}
\newcommand{\bmini}{\begin{center}\begin{minipage}{5in}}
\newcommand{\emini}{\end{minipage}\end{center}}
\newcommand{\R}{{\mbox{\bbb R}}}
\newcommand{\C}{{\mbox{\bbb C}}}
\newcommand{\Z}{{\mbox{\bbb Z}}}

\newcommand{\Q}{{\mathbb{Q}}}

\newcommand{\A}{{\mathcal{A}}}

\newcommand{\conf}{{\mathcal{C}}}

\newcommand{\Ordo}{{\mathcal{O}}}
\newcommand{\cand}{{\mathcal{W}}}
\newcommand{\sblv}{{\mathcal{H}}}
\newcommand{\ix}{\operatorname{Index}}
\newcommand{\cokrn}{\operatorname{Coker}}
\newcommand{\aut}{{\mathcal{A}}}
\newcommand{\M}{{\mathcal{M}}}

\newcommand{\dbar}{{\bar{\partial}}}

\newcommand{\ZZ}{{\mathcal{Z}}}
\newcommand{\End}{\operatorname{End}}

\newcommand{\la}{\langle}
\newcommand{\ra}{\rangle}
\newcommand{\pa}{\partial}
\newcommand{\id}{\operatorname{id}}
\newcommand{\Lag}{\operatorname{Lag}}

\newcommand{\spa}{\operatorname{Span}}
\newcommand{\krn}{\operatorname{Ker}}
\newcommand{\img}{\operatorname{Im}}

\newcommand{\diag}{\operatorname{Diag}}

\begin{document}
\title[Orientations in contact homology]{Orientations in Legendrian Contact
Homology and Exact Lagrangian Immersions}
\author{Tobias Ekholm}
\address{Department of Mathematics,
University of Southern California,
Los Angeles, CA 90089-1113}
\author{John Etnyre}
\address{Department of Mathematics, University of Pennsylvania, 209 South 33rd Street,
        Philadelphia PA 19105-6395}
\author{Michael Sullivan}
\address{Department of Mathematics, University of Massachusetts,
Lederle Graduate Research Tower, Amherst, MA 01003-9305}

\begin{abstract}
We show how to orient moduli spaces of holomorphic disks with
boundary on an exact Lagrangian immersion of a spin manifold into
complex $n$-space in a coherent manner. This allows us to lift the
coefficients of the contact homology of Legendrian spin submanifolds
of standard contact $(2n+1)$-space from $\Z_2$ to $\Z.$ We
demonstrate how the $\Z$-lift provides a more refined invariant of
Legendrian isotopy. We also apply contact homology to produce lower
bounds on double points of certain exact Lagrangian immersions into
$\C^n$ and again including orientations strengthens the results.
More precisely, we prove that the number of double points of an
exact Lagrangian immersion of a closed manifold $M$ whose associated
Legendrian embedding has good DGA is at least half of the dimension
of the homology of $M$ with coefficients in an arbitrary field if
$M$ is spin and in $\Z_2$ otherwise.
\end{abstract}


\maketitle

\tableofcontents

\section{Introduction}
Legendrian contact homology has been an effective tool in studying
Legendrian submanifolds in $\R^{2n+1}.$ In $\R^3$, Chekanov
\cite{Chekanov} and Eliashberg and Hofer (unpublished but see
\cite{Eliashberg00}) used contact homology to show that Legendrian
knots are not determined up to Legendrian isotopy by the so-called
classical invariants (topological isotopy class,
Thurston-Bennequin number, and Maslov class). Subsequently,
contact homology has been used to greatly illuminate the nature of
Legendrian knots in $\R^3.$ The contact homology of Legendrian
submanifolds in $\R^{2n+1}\approx\C^n\times\R$ (for $n>1$) was
given a rigorous foundation in \cite{ees2} and its efficacy was
demonstrated in \cite{ees1}.
Very roughly speaking contact homology is the homology of a
differential graded algebra (DGA) associated to a Legendrian
submanifold $L\subset\C^n\times\R$. The algebra is generated by
double points in the (Lagrangian) projection of $L$ into $\C^n$ and
the differential counts rigid holomorphic disk with corners at these
double points and boundary on the projected Legendrian submanifold.
In the initial definition of contact homology the disks were counted
modulo 2 since in that version of the theory orientations and
orientability of spaces of holomorphic disks need not be considered.
A $\Z$-lift of contact homology of Legendrian knots in $\R^3$ have
been introduced in a purely combinatorial fashion in \cite{ENS}. It
is however still not known if the oriented version of the theory in
this case is any stronger that the unoriented version of the theory.
Orientations for the moduli space of certain Riemann surfaces
without boundary has been discussed in \cite{ FH, EGH, BM}.
In this paper we show how to lift the DGA of Legendrian
submanifolds, of $\R^{2n+1}$, which are spin to $\Z$. We demonstrate
that this lift gives a more refined invariant of Legendrian isotopy
than does the theory over $\Z_2$ in dimensions $2n+1\ge 9$. For
Legendrian knots in $\R^3$, our analytical approach to orientations
recovers the combinatorial sign rule of \cite{ENS} and furthermore
gives rise to another combinatorial sign rule not mentioned there.
We also use Legendrian contact homology to produce lower bounds on
the double points of exact Lagrangian immersions into $\C^n.$ (A
Lagrangian immersion $f\colon M\to \C^n$ is {\em exact} if the
closed form $f^\ast(\sum_{j=1}^n y_j\,dx_j)$, where
$(x_1+iy_1,\dots,x_n+iy_n)$ are standard coordinates on $\C^n,$ is
exact.) Generically an exact Lagrangian immersion can be lifted to a
Legendrian embedding. A DGA is called good if it is (tame)
isomorphic to a DGA without constant terms in its differential \cite{Chekanov}. We
show that
if $f\colon M\to \C^n$ is an exact self-transverse Lagrangian
immersion of a closed manifold such that the DGA associated to a
Legendrian lift of $f$ is good
then the number $R(f)$ of double points of $f$ satisfies
\begin{equation}\label{eqhalfhom}
R(f)\ge \frac12\dim(H_\ast(M;\Lambda)),
\end{equation}
where $\Lambda=\Q$ or $\Lambda=\Z_p$ for any prime $p$ if $M$ is spin
and where $\Lambda=\Z_2$ otherwise.
It is easy to construct exact Lagrangian immersions of spheres and
tori of arbitrary dimensions which shows that the estimate
\eqref{eqhalfhom} is the best possible. While the hypothesis on the
exact Lagrangian immersion seems somewhat unnatural it is frequently
satisfied and from anecdotal evidence one would expect exact
Lagrangian immersions with non-good DGA's to have more double points
than ones with good DGA's. Despite this evidence it does not seem
straightforward to use contact homology for estimates when the
algebra is not good. However, we prove that if one can establish an
estimate like \eqref{eqhalfhom} with any fixed constant subtracted
from the right hand side then \eqref{eqhalfhom} is true too.
The paper is organized as follows. In Section \ref{sec:basnot} we
introduce basic notions which will be used throughout the paper. In
Section \ref{sec:orimdli} we show how to orient moduli spaces of
holomorphic disks relevant to contact homology. To accomplish this
we discuss orientations of determinant bundles over spaces of
(stabilized) $\bar\pa$-operators associated to Legendrian
submanifolds and their interplay with orientions of spaces of
conformal structures on punctured disks. Similar constructions are
carried out in \cite{FOOO} but some of the details differ. In
Section \ref{sec:LegCH} we define the DGA associated to a Legendrian
spin submanifold $L$ as an algebra over $\Z[H_1(L)]$ with
differential $\pa$ and prove that $\pa^2=0$. Furthermore we prove
the invariance of contact homology under Legendrian isotopy by a
mixture of a homotopy method and the more direct bifurcation
analysis, making use of the stabilization mentioned above. (Over
$\Z_2$ this invariance proof gives an alternative to the invariance
proof given in \cite{ees2}.) We also describe how the contact
homology depends on the choice of spin structure of the Legendrian
submanifold and we derive diagrammatic sign rules for Legendrian
knots in $\R^3$. In Section \ref{sec:lmt}, we adapt a theorem of
Floer \cite{Floer89a} to our situation so that, in special cases,
the differential in contact homology can be computed. We also apply
these results to construct examples which demonstrates that contact
homology over $\Z$ is more refined than contact homology over
$\Z_2$. In Section \ref{sec:dbpt} we prove the results related to
the double point estimate for exact Lagrangian immersion mentioned
above.

Acknowledgments: The authors are grateful to Lenny Ng for many useful discussions concerning the
sign rules in dimension three. We also thank AIM who provided some support during a workshop where
part of this work was completed. Part of this work was done while TE was a research fellow of the
Swedish Royal Academy of Sciences sponsored by the Knut and Alice Wallenberg foundation. MS was
partially supported by an NSF VIGRE grant as well as NSF Grant DMS-0305825. He also thanks the
University of Michigan and MSRI for hosting him while working on this paper. JE was partially
supported by NSF Grant DMS-0203941, an NSF CAREER Award (DMS--0239600) and FRG-0244663.

\section{Basic notions}\label{sec:basnot}
In this section we introduce notation and briefly describe the
construction of the contact homology of a Legendrian submanifold in
$\R^{2n+1}$. For more details on this construction, see \cite{ees1,
ees2}. Let $(x_1,y_1,\dots, x_n,y_n,z)$ be standard coordinates on
$\R^{2n+1}.$ Throughout this paper we consider the standard contact
structure $\xi=\ker(dz-\sum y_j\, dx_j)$ on $\R^{2n+1}$ and we make
use of the following two projections: the {\em front projection}
$\Pi_F:\R^{2n+1}\to\R^{n+1}$ which forgets the $y_j$-coordinates and
the {\em Lagrangian projection} $\Pi_{\C}:\R^{2n+1}\to \R^{2n}$
which forgets the $z$-coordinate. Note that if $L\subset\R^{2n+1}$
is a Legendrian submanifold then $\Pi_\C\colon L\to\R^{2n}=\C^n$ is
a Lagrangian immersion. The {\em Maslov class} of a Lagrangian
immersion $\phi\colon M\to\C^n$ is the homomorphism $\mu\colon
H_1(M)\to\Z$ such that $\mu(A)$ is the Maslov index of the loop of
Lagrangian planes which consists of tangent planes to $\phi(M)$
along some loop $\gamma\subset M$ representing $A$. The {\em Maslov
number} $m_\phi>0$ of $\phi$ is the positive integer which is a
generator of the image of $\mu_\phi\colon H_1(L;\Z)\to\Z$ if
$\mu_\phi\ne 0$ and $m_\phi=0$ if $\mu_\phi=0$. If $L\subset
\C^n\times\R$ is a Legendrian immersion we sometimes write $m(L)$
for the Maslov number of the Lagrangian immersion $\Pi_\C\colon
L\to\C^n.$
We define the contact homology of a spin Legendrian submanifold $L$
in $(\R^{2n+1},\xi)$ equipped with a spin structure as an algebra
over the group ring $\Z[H_1(L)]$ if $L$ is connected and as an
algebra over $\Z$ otherwise in the following way. Consider first the
connected case. Assume that $L$ is generic with respect to the
Lagrangian projection and let $\A$ be the free associative
(non-commutative) algebra over $\Z[H_1(L)]$ generated by the
(transverse) double points of $\Pi_{\C}(L)\subset\R^{2n}=\C^n.$
There is a $\Z$ grading on this algebra defined as follows. For a
double point $c_i$ denote the two points in $L$ mapping to $c_i$ by
$c_i^+$ and $c_i^-$ where $c_i^+$ has the larger $z$-coordinate.
Choose, for each double point $c_i,$ a path $\gamma_i$ in $L$ that
runs from $c_i^+$ to $c_i^-.$ The Conley-Zehnder index $\nu(c_i)$ of
$c_i$ is the Maslov index of the loop of Lagrangian subspaces in
$\C^n$ which is the path of planes tangent to $\Pi_\C(L)$ along
$\gamma$ closed up in a specific way, see \cite{ees2}. The grading
of $c_i$ is defined as $|c_i|=\nu(c_i)-1$ and the grading of a
homology class $A\in H_1(L)$ is defined as the negative of the
Maslov index of the loop of Lagrangian subspaces in $\C^n$ which are
tangent to $\Pi_\C(L)$ along some loop $\gamma\subset L$
representing $A$. If one defines the algebra over $\Z$ instead of
$\Z[H_1(L)]$ the grading is only in $\Z/m(L)\Z$. In the disconnected
case $\A$ is graded over $\Z_2$, see Subsection \ref{defofal}.
(There is also a relative $\Z/m(L)\Z$-grading which we will not
discuss in this paper.)

The differential of the algebra $\A$, $\pa\colon \A\to\A$ lowers
grading by $1$ and is defined by counting holomorphic disks in
$\C^n$ with boundary on $\Pi_\C(L)$. More precisely, for $a$, a
double point of $\Pi_\C(L)$, ${\mathbf b}$, a word in the double
points of $\Pi_\C(L)$, and $A\in H_1(L)$, we consider the moduli
space $\M_A(a;{\mathbf b})$. This space consists of holomorphic maps
from the punctured unit disk to $\C^n$ which maps the boundary to
$\Pi_\C(L)$, which are asymptotic, see \cite{ees2}, to the double
points specified at the punctures, and which when restricted to the
boundary satisfies a certain homology condition specified by $A$. If
$L$ is connected then the differential on $\A$ is
\begin{equation}\label{eqdifferential}
\pa a= \sum_{\mathbf{b}} (\# \M_A(a;\mathbf{b})) A\mathbf{b},
\end{equation}
where the sum runs over all ${\mathbf b}$ such that
$\M_A(a;\mathbf{b})$ is $0$-dimensional. The moduli spaces which are
$0$-dimensional are compact manifolds. In Section~\ref{OMS} we
orient $\M_A(a;\mathbf{b})$. This means in particular that each
component of a $0$-dimensional moduli space comes equipped with a
sign and $\# \M_A(a;\mathbf{b})$ in \eqref{eqdifferential} is the
algebraic number of points in the moduli space. If $L$ is not
connected then the differential is defined as in
\eqref{eqdifferential} except that the homology class $A$ there
should be deleted. The {\em contact homology} of $L$ is
$CH(L)=\krn(\pa)/\img(\pa)$. Finally, we note that if
$L\subset\C^n\times\R$ is a Legendrian submanifold then double
points of $\Pi_\C\colon L\to\C^n$ correspond to segments in the
$\R$-direction of $\C^n\times\R$ with its endpoints on $L$. The
vector field $\pa_z$ is the Reeb field of the contact form
$dz-\sum_j y_j\,dx_j$ and thus such a segment is called a Reeb
chord. We will therefore use the words {\em Reeb chord} and {\em
double point} interchangeably below.
\section{Orienting the moduli spaces}\label{sec:orimdli}
In this section we orient moduli spaces of holomorphic disks in
$\C^n$ with boundary on $\Pi_\C(L)$, where $L\subset\C^n\times\R$ is
a Legendrian submanifold which is spin. To orient the moduli spaces
we find an orientation of the determinant bundle of a stabilized
version of the linearization of the defining $\bar\pa$-equation. The
source space of the linearization splits into an infinite
dimensional space and a finite dimensional space arising from
automorphisms or variations of the conformal structure of the source
disk. The restriction of the linearization to the infinite
dimensional space is a Fredholm operator. We orient the determinant
bundles over spaces of such Fredholm operators and orient spaces of
automorphisms and conformal structures separately.
The orientations we define depend on several choices. We make these
choices so that it is possible to endow the graded algebra discussed
in Section \ref{sec:basnot} with a differential. In particular, this
differential must respect the multiplication of the algebra in the
sense that it satisfies the graded Leibniz rule, see Equation
\eqref{eqleibniz}.

\subsection{Lagrangian boundary conditions for the $\bar\pa$-operator
on punctured disks}\label{ssec:boundcond}
Let $D_m$ denote the unit disk $D$ in $\C$ with $m$ distinct
punctures $\{p_1,\dots,p_m\}$ on the boundary, $m\ge 0$. The
orientation on $D$ induces an orientation on its boundary. If one
puncture, $p_1$ say, is distinguished then the orientation on $\pa
D_m$ induces an order of the punctures. The punctures subdivide $\pa
D_m$ into $m$ disjoint oriented open arcs. We denote their closures
by $I_1,\dots,I_m$, with notation such that $\pa I_j=p_{j+1}-p_j$
(where $p_{0}=p_m$ and $p_{m+1}=p_1$). For convenience we write
$p_j^+$ and $p_j^-$ for $p_j$ thought of as a point in $I_{j-1}$ and
$I_j$, respectively. (We let $I_{0}=I_m$, $I_{m+1}=I_1$, and in the
special case when $m=0$, $I_0=I_1=\pa D$.)
A {\em Lagrangian boundary condition} on $D_m$ is a collection of maps
$$
\lambda=(\lambda_1,\dots,\lambda_m)
$$
where $\lambda_j\colon I_j\to\Lag(n)$, where $\Lag(n)$ denotes the
space of Lagrangian subspaces of $\C^n$. For non-punctures $z\in\pa
D_m$ we write $\lambda(z)$ to denote $\lambda_j(z)$ where $j$ is the
unique subscript such that $z\in I_j$. If $p$ is a puncture then
$p\in I_{j-1}$ and $p\in I_j$. We write
$\lambda(p^+)=\lambda_{j-1}(p^+)$ and
$\lambda(p^-)=\lambda_j(p_j^-)$.
An {\em oriented Lagrangian boundary condition} on $D_m$ is a collection of maps
$$
\tilde\lambda=(\tilde\lambda_1,\dots,\tilde\lambda_m)
$$
where $\tilde\lambda_j\colon I_j\to\Lag_0(n)$, where $\Lag_0(n)$ denotes the space of oriented
Lagrangian subspaces of $\C^n$. Note that any Lagrangian boundary condition on a disk with at least
one puncture lifts (non-uniquely) to an oriented boundary condition, and that there are boundary
conditions on the $0$-punctured disk which do not lift.
A {\em trivialized Lagrangian boundary condition} for the $\bar\pa$-operator on $D_m$ is a
collection of maps
$$
A=(A_1,\dots,A_m),
$$
where $A_j\colon I_j\to U(n)$. A trivialized Lagrangian boundary condition induces a Lagrangian
boundary condition $\lambda$ via $\lambda_j(z)=A_j(z)\R^n$ and a trivialization
$(v_1(z),\dots,v_n(z))$ of a Lagrangian boundary condition (i.e. an ON-basis in $\lambda$) gives a
trivialized Lagrangian boundary condition by defining $A(z)$ as the matrix with column vectors
$(v_1(z),\dots,v_n(z))$. Fixing an orientation on $\R^n$, a trivialized boundary condition can be
considered as an oriented boundary condition.
\subsubsection{Angles, weights, and Fredholm properties}
Neighborhoods of the punctures of $D_m$ will be thought of as infinite half strips
$[0,\infty)\times[0,1]$ or $(-\infty,0]\times[0,1]$ with coordinates $\tau+it$, and we consider the
$\bar\pa$-operator
\begin{equation}\label{RH}
\bar\pa\colon \sblv_{2,\epsilon}[\lambda](D_m,\C^n)\to
\sblv_{1,\epsilon}[0](D_m, {T^\ast}^{0,1}D_m\otimes\C^n),
\end{equation}
where $\epsilon=(\epsilon_1,\dots,\epsilon_m)\in\R^m$, and
$\sblv_{2,\epsilon}[\lambda]$ is the closed subspace of the Sobolev
space $\sblv_{2,\epsilon}$, with weight $e^{\epsilon_j|\tau|}$ in
the neighborhood of the $j^{\rm th}$ puncture, which consists of
elements $u$ with $u(z)\in\lambda_j(z)$ for $z\in I_j$, and $\bar\pa
u=0$ on $\pa D_m$ and where $\sblv_{1,\epsilon}[0]$ is the closed
subspace of the Sobolev space $\sblv_{1,\epsilon}$ which consist of
elements with vanishing trace (restriction to the boundary). See
\cite{ees2}, Section 5. We denote the operator in \eqref{RH} by
$\bar\pa_{\lambda,\epsilon}$. Often in our applications, the weight
will be clear from the context and we drop it from the notation
writing simply $\bar\pa_\lambda$.

We recall the following notion from \cite{ees2}. Let $V_1$ and $V_2$ be Lagrangian subspaces of
$\C^n$. Define the {\em complex angle} $\theta(V_1,V_2)\in[0,\pi)^n$ inductively as follows. If
$\dim(V_1\cap V_2)=r\ge 0$ let $\theta_1=\dots=\theta_r=0$ and let $\C^{n-r}$ denote the Hermitian
complement of $\C\otimes(V_1\cap V_2)$ and let $V_j'=V_j\cap\C^{n-r}$ for $j=1,2$. If $\dim(V_1\cap
V_2)=0$ then let $V_j'=V_j$, $j=1,2$ and let $r=0$. Then $V'_1$ and $V'_2$ are Lagrangian subspaces.
Let $\alpha$ be smallest angle such that $\dim(e^{i\alpha}V'_1\cap V'_2)=r'>0$. Let
$\theta_{r+1}=\dots=\theta_{r+r'}=\alpha$. Repeat the construction until $\theta_n$ has been
defined.
Let $\lambda\colon \pa D_m\to\Lag(n)$ be a Lagrangian boundary condition. Let
$\theta^j=(\theta^j_1,\dots,\theta^j_n)\in [0,\pi)^n$ denote the complex angle of $\lambda(p_j^+)$
and $\lambda(p_j^-)$. As shown in \cite{ees1}, the operator $\bar\pa_{\lambda,\epsilon}$ is Fredholm
when
\begin{equation}\label{ellipt}
\epsilon_j\ne -\theta^j_k+l\pi,  \text{ for all $1\le j\le m$,
$1\le k\le n$, and all $l\in\Z$.}
\end{equation}
The {\em determinant line} $\det(F)$ of a Fredholm operator $F\colon B_1\to B_2$ between Banach
spaces is defined as the tensor product of the highest exterior power of its kernel and the highest
exterior power of the dual of its cokernel, respectively. That is,
$$
\det(F)=\Lambda^{\rm max}\krn(F)\otimes\Lambda^{\rm max}\cokrn(F)^\ast.
$$
Let
$$
X_m\subset C^{\infty}\left(\sqcup_{j=1}^m I_j;\,\Lag(n)\right)\times\R^m
$$
denote the space of all Lagrangian boundary conditions and weights
for which \eqref{ellipt} is satisfied. It is a standard result that
there exists a determinant line bundle $E$ over $X_m$ with fiber
over $(\lambda,\epsilon)$ equal to
$\det(\bar\pa_{\lambda,\epsilon})$. The bundle $E$ is not
orientable. In fact, the corresponding determinant bundle over
oriented Lagrangian boundary conditions on the $0$-punctured disk is
non-orientable. What is needed for orientability in that case are
trivialized boundary conditions, see Lemma \ref{lmacanor}.
\subsubsection{Parity and signs of punctures}\label{typesofpath}
We say a boundary condition $\lambda\colon\pa D_m\to\Lag(n)$, $m>0$ is {\em transverse} if
$\lambda(p^+)$ and $\lambda(p^-)$ are transverse for each puncture $p$. We will subdivide the
punctures on a disk with such boundary conditions into four classes. The first subdivision seems
rather ad hoc at this point but in our applications to Legendrian submanifolds it has a clear
geometric meaning: each puncture has a sign. That is, there are {\em positive} and {\em negative}
punctures. (In the applications, disks have one positive puncture and all other negative.)
The second subdivision comes from the {\em parity} of a puncture defined in the following way. Let
$V_1$ and $V_2$ be two {\em oriented} transverse Lagrangian subspaces of $\C^n$. Let
$\theta=\theta(V_1,V_2)$ be the complex angle of $(V_1,V_2)$. Then
there exists complex coordinates
$(z_1,\dots,z_n)=(x_1+iy_1,\dots,x_n+iy_n)$ on $\C^n$ so that in these
coordinates
$$
V_1=\spa(\pa_{x_1},\dots,\pa_{x_n})
$$
and
$$
V_2=\spa(e^{i\theta_1}\pa_{x_1},\dots, e^{i\theta_n}\pa_{x_n}).
$$
We will call such coordinates
{\em canonical coordinates of $(V_1,V_2)$}. Note that the canonical
coordinates are not unique. However, in the constructions below the
choice of specific canonical coordinates will be irrelevant. For
example, the unitary linear map of $\C^n$ with matrix
\begin{equation}\label{negrot+to+}
\diag\left(e^{-i(\pi-\theta_1)},e^{-i(\pi-\theta_2)},\dots,e^{-i(\pi-\theta_n)}\right),
\end{equation}
in canonical coordinates is well-defined and maps $V_1$ isomorphically to $V_2$.
The ordered pair $(V_1,V_2)$ is {\em even (odd)} if
\eqref{negrot+to+} is an orientation reversing (preserving) map from
$V_1$ to $V_2.$  Consider a transverse oriented boundary
condition $\lambda\colon\pa D_m\to\Lag_0(n)$. If $p$ is a negative puncture then $p$ is {\em even}
if the pair $(\lambda(p^+),\lambda(p^-))$ is even and $p$ is {\em odd} if
$(\lambda(p^+),\lambda(p^-))$ is odd. If $q$ is a positive puncture then $q$ is {\em even} ({\em
odd}) if the pair $(\lambda(q^-),\lambda(q^+))$ is even (odd).

\subsection{Gluing operations}\label{gluingx3}
We define gluing operations for the boundary conditions discussed in
Subsection \ref{ssec:boundcond} and explain how
orientations of the determinant lines over the pieces relate to an
orientation of the determinant line of the resulting glued boundary
condition. More precisely, the gluing operations give rise to exact
sequences involving kernels and cokernels of operators and
the orientations of the determinant lines are related via these.
In our applications, the exact form of
the relations is important. We therefore start out with explaining
our conventions for exact sequences of oriented vector spaces.
\subsubsection{Orientation conventions}\label{ssseclinalg}
All our gluing operations give rise to exact
sequences with at most four non-zero terms so we
discuss only this case. Let
\begin{equation}\label{arcseq}
\begin{CD}
0 @>>> V_1 @>{\alpha}>> W_1 @>{\beta}>> W_2 @>{\gamma}>> V_2 @>>> 0,
\end{CD}
\end{equation}
be an exact sequence of finite dimensional vector spaces. This
sequence induces an isomorphism
\begin{equation}\label{trivia}
\Lambda^{\rm max} V_1\otimes \Lambda^{\rm max} V_2^\ast\approx
\Lambda^{\rm max} W_1\otimes \Lambda^{\rm max} W_2^\ast.
\end{equation}
Our interest in this isomorphism is its effect on
orientations. Note that there is a natural correspondence between orientations on a
finite dimensional vector space $V$ and on its dual $V^\ast$. To see
this, think of an orientation of $V$ as a non-zero
element in $\Lambda^{\rm max} V$ up to multiplication by a positive
number. The correspondence can then be obtained as follows. Pick any
basis $v_1,\dots,v_n$ in $V$ and let $v_1^\ast,\dots,v_n^\ast$ be the
dual basis in $V^\ast$. Now identify the orientation in $V$ given by
$$
v_1\wedge\dots\wedge v_n
$$
with the orientation in $V^\ast$ given by
$$
v_1^\ast\wedge\dots\wedge v_n^\ast.
$$
It is easy to see that this identification is independent of the
choice of basis.
If $V$ and $W$ are finite dimensional vector spaces then the
correspondence just described gives rise to a correspondence between
orientations of the $1$-dimensional vector space
$$
\Lambda^{\rm max}V\otimes\Lambda^{\rm max}W^\ast
$$
and pairs $(o_V,o_W)$ of orientations of $V$ and $W$ respectively,
modulo the following equivalence relation: a pair represented by
$(\phi_V,\phi_W)\in \Lambda^{\rm max}V\times\Lambda^{\rm max}W$ is
identified with the pair represented by $(a\phi_V,a\phi_W)$ for any
non-zero $a\in\R$. We call an equivalence class an {\em orientation
pair}. Using this terminology, we interpret the isomorphism
\eqref{trivia} on the orientation level as follows: the sequence
\eqref{arcseq} gives a correspondence between orientation pairs of
$(V_1,V_2)$ with orientation pairs of $(W_1,W_2)$. We next give a
concrete description of this correspondence. Let $\xi$ be an
orientation pair on $(V_1,V_2)$. Pick bases $v_1^1,\dots,v_1^m$ in
$V_1$ and $v_2^1,\dots,v_2^s$ in $V_2$ such that
$$
\left(v_1^1\wedge\dots\wedge v_1^m,
v_2^1\wedge\dots\wedge v_2^s\right),
$$
represents $\xi$. Pick vectors $w^{m+1}_1,\dots,w_1^{n}$ such that
$\alpha(v_1),\dots,\alpha(v_m),w_1^{m+1},\dots,w_1^n$ is a basis of
$W_1$. By exactness, the vectors $\beta(w_1^{m+1}),\dots,\beta(w_1^n)$
are linearly independent in $W_2$. Let $w_2^1,\dots,w_2^s$ be any
vectors in $W_2$ such that $\gamma(w_2^j)=v_2^j$ for all $j$. Define
$\phi(\xi)$ to be the orientation pair represented by
$$
\left(\alpha(v_1^1)\wedge\dots\wedge \alpha(v_1^m)\wedge
w_1^{m+1}\wedge\dots\wedge w_1^n,
\beta(w_1^{n})\wedge\dots\wedge\beta(w_1^{m+1})\wedge
w_2^1\wedge\dots\wedge w_2^s\right).
$$
(Note the {\em reverse order} of the vectors $\beta(w_1^j)$.)
It is easy to see that $\phi(\xi)$ is well defined. The inverse $\psi$
of $\phi$ can be described as follows. Let $\eta$ be any orientation
pair of $(W_1,W_2)$. Pick bases $w_1^1,\dots,w^n_1$ in $W_1$ and
$w_2^1,\dots w_2^t$ in $W_2$ such that
$$
\left(w_1^1\wedge\dots\wedge w^n_1,
w_2^1\wedge\dots\wedge w_2^t\right)
$$
represents $\eta$. Pick any basis $v_1^1,\dots,v_1^m$ in $V_1$ then
there are vectors ${\hat w}_1^{m+1},\dots,{\hat w}_1^n$ in $W_1$ so that
$$
w_1^1\wedge\dots\wedge w^n_1\quad\text{ and }\quad
\alpha(v_1^1)\wedge\dots\wedge\alpha(v_1^m)\wedge
{\hat w_1}^{m+1}\wedge\dots\wedge{\hat w}_1^n
$$
represent the same orientation on $W_1$. Furthermore, there are vectors
${\hat w}_2^1,\dots,{\hat w}_2^s$ such that
$$
w_2^1\wedge\dots\wedge w_2^t\quad\text{ and }\quad
\beta({\hat w_1}^{n})\wedge\dots\wedge\beta({\hat w}_1^{m+1})\wedge
{\hat w}_2^1\wedge\dots\wedge{\hat w}_2^s
$$
represent the same orientation on $W_2$. Define $\psi(\eta)$ to be the
orientation pair on $(V_1,V_2)$ represented by
$$
\left(v_1^1\wedge\dots\wedge v_1^m,
\gamma({\hat w}_2^1)\wedge\dots\wedge\gamma({\hat w}_2^s)\right).
$$
It is easy to verify that $\psi$ is well defined and that $\psi$ and
$\phi$ are inverses of each other.

\subsubsection{Two gluing operations}\label{ssec2gluings}
Let $\lambda\colon\pa D_m\to\Lag(n)$ be a transverse boundary condition
on the punctured disk. Let $\bar\pa_\lambda$ be the operator with this
boundary condition and with trivial weights. Then, by Proposition~5.14 in
\cite{ees2},
$$
\ix(\bar\pa_\lambda)=n+\mu(\hat\lambda),
$$
where $\mu$ is the Maslov index and
where $\hat\lambda$ is the loop in $\Lag(n)$ obtained from the map
$\lambda\colon \pa D_m\to\Lag(n)$ as follows. Let
$\theta\in(0,\pi)^n$ denote the complex angle between $\lambda(p_j^+)$
and $\lambda(p_j^-)$ then there exists coordinates $(x_1+iy_1,\dots,x_n+iy_n)$
in $\C^n$ such that
$$
\lambda(p_j^+)=\spa(\pa_{x_1},\dots,\pa_{x_n})
$$
and
$$
\lambda(p_j^-)=\spa(e^{i\theta_1}\pa_{x_1},\dots,e^{i\theta_n}\pa_{x_n}).
$$
Let $\gamma_j(s)$, $0\le s\le 1$, be the path of Lagrangian subspaces
given in these coordinates as
$$
\spa(e^{-i(\pi-\theta_1)s}\pa_{x_1},\dots,e^{-i(\pi-\theta_n)s}\pa_{x_n}).
$$
Then $\hat\lambda$ is the loop obtained by concatenation of the paths
in $\lambda$ with the paths $\gamma_j$. That is,
$$
\hat\lambda=\lambda_1\ast\gamma_2\ast\lambda_2
\ast\dots\ast\lambda_m\ast\gamma_1.
$$
Let $A\colon \pa D_m\to U(n)$ and $B\colon \pa D_s\to U(n)$ be
trivialized Lagrangian boundary conditions. Let $p\in D_m$ and $q\in
D_s$ be points which are not punctures and assume that $A(p)=B(q)$.
We assume also that $A$ and $B$ are constant in neighborhoods of $p$
and $q$, respectively. (If this is not the case, we may homotope the
boundary conditions so that they become constant around $p$ and $q$.
Note that such homotopies may change the dimensions of the kernels
and cokernels of the operators involved. However, if the homotopy
can be chosen sufficiently $C^1$-small then the dimensions can be
kept constant.)
We glue $D_m$ and $D_s$ to one disk $D_{m+s}$ with a
trivialized Lagrangian boundary condition which is the concatenation
of $A$ and $B$.
To define this operation more accurately we proceed as follows.
Puncture the disks $D_m$ at $p$ and $D_s$ at
$q$. Use conformal coordinates $[0,\infty)\times[0,1]$ in a
neighborhood of $p$ in $D_m$ and conformal coordinates
$(-\infty,0]\times [0,1]$ in a neighborhood of $q$ in $D_s$.
Since the boundary conditions $A$ and $B$ are constant in
neighborhoods of $p$ and $q$, respectively, it follows that for all
sufficiently large $\rho>0$ the boundary condition $A$ is constant in
$[\rho,\infty)\times(\pa[0,1])$ and the boundary condition $B$ is
constant in $(-\infty,-\rho]\times(\pa[0,1])$. For such $\rho$ we
define a new disk $D(\rho)$ with $m+s$ punctures as follows. Let
$D'(\rho)=D_m\setminus([\rho,\infty)\times[0,1])$, let
$D''(\rho)=D_s\setminus((-\infty, -\rho]\times[0,1])$, and let
$D(\rho)$ be the disk which is obtained from gluing $D'(\rho)$ to
$D''(\rho)$ by identifying $\{\rho\}\times[0,1]$ with
$\{-\rho\}\times[0,1]$. Note that the boundary conditions glue in a
natural way to a boundary condition $A\sharp B$ on $D(\rho)$.
Since the Maslov index is additive under this gluing we find that
\begin{equation}\label{sum1index}
\ix(\bar\pa_{A\sharp_\rho B})=
n+\mu(\hat A)+\mu(\hat B)=\ix(\bar\pa_A)+\ix(\bar\pa_B)-n.
\end{equation}
We call this type of gluing {\em gluing at a boundary point}.

\begin{lma}\label{lmaglue1}
For all sufficiently large $\rho$ there exists an exact sequence
\begin{equation}\label{eqgluclosed}
\begin{split}
&\begin{CD}
0 @>>> \krn(\bar\pa_{A\sharp_\rho B}) @>{\alpha_\rho}>>
\krn(\bar\pa_A)\oplus\krn(\bar\pa_B)
@>{\beta_\rho}>>
\end{CD}\\
&\quad\quad
\begin{CD}
\cokrn(\bar\pa_A)\oplus\R^n\oplus\cokrn(\bar\pa_B)
@>{\gamma_\rho}>>
\cokrn(\bar\pa_{A\sharp_\rho B}) @>>> 0
\end{CD}
\end{split}
\end{equation}
where $\R^n$ is naturally identified with
$A(p)\R^n=B(q)\R^n$. (We describe the
maps  $\alpha_\rho$, $\beta_\rho$, and $\gamma_\rho$ and the identification in
Remark~\ref{rmkmapsinseq} below.)
In particular, together with an orientation on $\R^n$ the above
sequence induces an isomorphism between
$$
\Lambda^{\rm max}(\krn(\bar\pa_A)\oplus\krn(\bar\pa_B))\otimes
\Lambda^{\rm max}(\cokrn(\bar\pa_A)\oplus\cokrn(\bar\pa_B))^\ast
$$
and
$$
\Lambda^{\rm max}\krn(\bar\pa_{A\sharp B})\otimes
\Lambda^{\rm max}\cokrn(\bar\pa_{A\sharp B})^\ast.
$$
\end{lma}
\begin{rmk}\label{FO^3}
The special case $m=s=0$ of this lemma (together with Lemmas
\ref{lmacanor} and \ref{lmagluecanor} below) fixes a misstatement in
Lemma 23.5 of \cite{FOOO} (on page 206). Consider two bundles of
index $0$ with no kernel and no cokernel then Lemma 23.5 implies
that the kernel and cokernel of the glued problem are trivial as
well. However this cannot be correct since the problem has index
$-n$ by \eqref{sum1index}.
\end{rmk}

\begin{pf}
Using Sobolev spaces with a small negative exponential weight
$-\delta$, $\delta>0$ at the punctures $p\in D_m$ and $q\in D_s$
(i.e. in coordinates $\tau+it\in \R\times[0,1]$ the weight functions
$w'$ and $w''$
equal $e^{-\delta|\tau|}$ in neighborhoods of the punctures and equal
$1$ elsewhere), we get a canonical identification of the kernels
and cokernels on the punctured and the non-punctured disks, see Lemma~5.2 in \cite{ees2}.
Consider the operator
$$
\bar\pa\colon \sblv_{2,\rho}[A\sharp B](D(\rho),\C^n)\to,
\sblv_{1,\rho}[0](D(\rho),{T^\ast}^{0,1}D(\rho)\otimes\C^n),
$$
where the subscript $\rho$ indicates that we use Sobolev norms with
weight functions which are the gluing of the weight functions on
the two disks. (Since the support of the glued weight function in
$D(\rho)$ is compact this weight function is not
very important but it will be convenient to use that norm in the
estimates below.) More precisely, this weight function is a smoothing
of the function which equals $w'$ on
$D_m\setminus([\rho,\infty)\times[0,1])$ and $w''$ on
$D_s\setminus((-\infty,-\rho]\times[0,1])$.
Pick bases $a_1,\dots,a_{r_1}$ and $b_1,\dots,b_{r_2}$ in
$\krn(\bar\pa_A)$ and $\krn(\bar\pa_B)$, respectively. Also, pick bases
$\alpha_1,\dots,\alpha_{s_1}$ and $\beta_1,\dots,\beta_{s_2}$ in
$\cokrn(\bar\pa_A)$ and $\cokrn(\bar\pa_B)$, respectively. (A priori,
elements in the cokernels are elements in the dual of
$\sblv_1[0]$. However, they can be represented by $L^2$-pairing with smooth
functions,  see Lemma~5.1 in \cite{ees2}, and we will think of them as such.)
Let $\phi'_\rho$ be a cut-off function which equals $1$ on
$D_m\setminus([\frac12\rho,\infty)\times[0,1])$, which equals $0$ on
$D(\rho)\setminus (D_m\setminus[\rho,\infty)\times[0,1])$ and with
$|d^k\phi|=\Ordo(\rho^{-1})$ for $k=1,2$. Let $\phi''_\rho$ be a similar
cut-off function on with support in the part of $D(\rho)$
corresponding to $D_s$.
Define $V_\rho\subset\sblv_{2,\rho}(D(\rho),\C^n)$ to be the
$L^2$-complement of the $r_1+r_2$ functions
$$
\phi'_\rho a_1,\dots,\phi'_\rho a_{r_1},\phi''_\rho
b_1,\dots,\phi''_\rho b_{r_2}.
$$
We claim that there exists a constant $C$ such that for all
sufficiently large $\rho$
$$
\|v\|_{2,\rho}\le C\|\bar\pa v\|_{1,\rho}, \text{ for }v\in V_\rho.
$$
Assume not, then there exists a sequence of functions
$v_j\in V_{\rho(j)}$, $\rho(j)\to\infty$ as $j\to\infty$ such that
\begin{align}\label{glu1ker=1}
&\|v_\rho(j)\|_{2,\rho(j)}=1\\\label{glu1kerto0}
&\|\bar\pa v_{\rho(j)}\|_{1,\rho(j)}\to 0, \text{ as }j\to\infty.
\end{align}
Let $\Theta_{a,\rho}\subset D(\rho)$ denote the subset
$$
[a,\rho]\times[0,1]\cup [-\rho,-a]\times[0,1],
$$
where the first part is a subset of $D_m$ and the second of $D_s$.
Let $\psi_\rho$ be a cut-off which equals $1$ in $\Theta_{\frac14\rho,\rho}$,
which equals $0$ in the complement of $\Theta_{0,\rho}$
and which has $|d^k\psi|=\Ordo(\rho^{-1})$, $k=1,2$. Then using the
elliptic estimate for the $\bar\pa$-operator on $\R\times[0,1]$ with
small positive exponential weight $e^{\epsilon|\tau|}$ at both ends we
find using \eqref{glu1kerto0} that
\begin{equation}\label{glu1mid}
\|\psi_{\rho(j)}v_{\rho(j)}\|_{2,\rho(j)}\to 0, \text{ as }j\to\infty.
\end{equation}
Let $\psi'_\rho$ be a cut off function which is
$1$ on $D_m\setminus([\frac14\rho,\infty)\times[0,1])$, which is $0$
outside $D_m\setminus([\rho,\infty)\times[0,1])$, and with
$|d^k\psi'_\rho|=\Ordo(\rho^{-1})$. Using the elliptic estimate on the
$L^2$-complement of the kernel of $\bar\pa$ on $D_m$ we find
\begin{equation}\label{glu1'}
\|\psi'_{\rho(j)}v_{\rho(j)}\|_{2,\rho(j)}\to 0, \text{ as }j\to\infty.
\end{equation}
(Note since the support of the cut-off functions $\phi'_\rho$ are
contained in the support of the cut-off functions $\psi'_\rho$
we know $\psi'_{\rho(j)}v_{\rho(j)}$ is orthogonal to the kernel of $\bar\pa$ on $D_m.$)
With $\psi''_\rho$ similarly defined but with support on $D_s$ instead we
find
\begin{equation}\label{glu1''}
\|\psi''_{\rho(j)}v_{\rho(j)}\|_{2,\rho(j)}\to 0, \text{ as }j\to\infty.
\end{equation}
However, \eqref{glu1mid}, \eqref{glu1'}, and \eqref{glu1''}
contradicts \eqref{glu1ker=1} and we find the estimate on $V_\rho$ holds as
claimed.
The restriction of $\bar\pa$ to $V_\rho$ is a Fredholm operator of
index
$$
-\dim(\cokrn(\bar\pa_A))-\dim(\cokrn(\bar\pa_B))-n.
$$
Consider the elements $\phi'\alpha_1,\dots,\phi'\alpha_{s_1}$ and
$\phi''\beta_1,\dots,\phi''\beta_{s_2}$ in the dual of
$\sblv_{1,\rho}[0]$. Furthermore, we pick in this dual constant
functions $c_1,\dots,c_n$ on the gluing region with values in the
orthogonal complement of $A(1)\R^n$ and consider $\psi_\rho
c_1,\dots,\psi_\rho c_n$ as elements in the dual. Moreover, choose
them so that their $L^2$-norms (in the dual weight, which is large
in the gluing region) are $1$.
Now pick elements (in $\sblv_{1,\rho}[0]$) dual to these elements and
all of norm one. Then they stay a uniform positive distance away from
the intersection
$$
W=\bigcap_j\krn(\psi'_\rho\alpha_j)\cap
\bigcap_k\krn(\psi''_\rho\beta_k)\cap
\bigcap_l\krn(\psi_\rho c_l),
$$
and thus give a direct sum decomposition
$$
\sblv_{1,\rho}[0]=W\oplus\R^{s_1+s_2+n}.
$$
Let $\pi$ be the projection induced by the above direct sum
decomposition.
We claim that
$$
\pi\circ\bar\pa\colon V_\rho\to W,
$$
is an isomorphism. Assume not, then there is a
sequence $u_\rho$ with
$$
\|u_\rho\|_{V_\rho} =1 \text{ and } \|\bar\pa u_\rho\|_W\to 0.
$$
We need only note that also $\xi(\bar\pa u_\rho)\to 0$ for $\xi$ one
of the chosen elements in the dual to see that
this contradicts the previous estimate. But this is clear since
$\pa \hat\xi=0$ where $\hat\xi$ is the function which was cut-off to
get $\xi$.

We thus find the exact sequence
\begin{align*}
&\begin{CD} 0 @>>> \krn(\bar\pa_{A\sharp_\rho B}) @>{\alpha_\rho}>>
\krn(\bar\pa_A)\oplus\krn(\bar\pa_B)
@>{\beta_\rho}>>
\end{CD}\\
&\quad\quad
\begin{CD}
\cokrn(\bar\pa_A)\oplus\R^n\oplus\cokrn(\bar\pa_B)
@>{\gamma_\rho}>> \cokrn(\bar\pa_{A\sharp_\rho B}) @>>> 0
\end{CD}
\end{align*}
where $\alpha_\rho$ is the inclusion followed by orthogonal projection to
the subspace spanned by
$\{\psi'_\rho a_j\}_{j=1}^{r_1}\cup\{\psi''_\rho b_k\}_{k=1}^{r_2}$,
where $\beta_\rho$ is the $\bar\pa$-operator followed by second
projection in the direct sum decomposition of $\sblv_{1,\rho}[0]$, and
where $\gamma_\rho$ is the inclusion followed by the projection to the
quotient $\sblv_{1,\rho}[0]/\img(\bar\pa)$.
\end{pf}

\begin{rmk}\label{rmkmapsinseq}
We note that the natural identification of $\R^n$ in the above
sequence is obtained by identifying the cokernel of the
$\bar\pa$-problem on the gluing strip with the constant functions
taking values in $i\cdot A(1)\R^n$.
The map $\alpha_\rho$ can be described as follows. First we map the
space $(\krn(\bar\pa_A)\oplus\krn(\bar\pa_B))$ into
$\sblv_{2}(D_\rho,\C^n)$ by cutting off kernel elements with cut-off
functions as in the proof. The obtained map is clearly injective.
The map $\alpha_\rho$ is then simply the
inclusion of $\krn(\bar\pa_{A\sharp B})$ into $\sblv_2(D_\rho,\C^n)$
followed by the $L^2$-orthogonal to the image of
$(\krn(\bar\pa_A)\oplus\krn(\bar\pa_B))$.
The map $\beta_\rho$ can be described as follows. First we map the
space $\cokrn(\bar\pa_A)\oplus\R^n\oplus\cokrn(\bar\pa_B)$ into the
space $\sblv_1(D_\rho, {T^\ast}^{0,1} D_\rho\otimes\C^n)$ using
cut-off functions as in the proof above. The map is
injective. We then consider $(\krn(\bar\pa_A)\oplus\krn(\bar\pa_B))$
as a subspace of $\sblv_{2}(D_\rho,\C^n)$ and map it to the subspace
of $\sblv_1(D_\rho, {T^\ast}^{0,1} D_\rho\otimes\C^n)$ with the
$\bar\pa$-operator followed by $L^2$-orthogonal projection.
The map $\gamma_\rho$ is simply mapping
$\cokrn(\bar\pa_A)\oplus\R^n\oplus\cokrn(\bar\pa_B)$, included into
$\sblv_1(D_\rho, {T^\ast}^{0,1} D_\rho\otimes\C^n)$ as above, to
$\cokrn(\bar\pa_{A\sharp B})$ by projection to the quotient.
\end{rmk}
We next consider another type of gluing. Let $A\colon \pa D_m\to U(n)$
and $B\colon \pa D_s\to U(n)$ be as above.
Let $p$ be a puncture of $D_m$ and let $q$ be a puncture on $D_r$.
Assume that $A(p^+)=B(q^-)$ and that
$A(p^-)=B(q^+)$ and that $A$ and $B$ are
constant in neighborhoods of $p^{\pm}$ and $q^{\pm}$, respectively.
(In general, our assumption that boundary conditions are $C^2$-small in
standard coordinates around the punctures implies that we can actually
homotope them to constant without changing the kernel/cokernel dimension.)
Choose conformal coordinates
$[0,\infty)\times[0,1]$ in a neighborhood of $p$ and conformal
coordinates $(-\infty,0]\times[0,1]$ in a neighborhood of $q$. For
$\rho>0$ large enough we can glue
$D_m\setminus([\rho,\infty)\times[0,1])$ to
$D_s\setminus((-\infty,-\rho]\times[0,1])$ to obtain a disk
$D(\rho)$ with $m+r-2$ punctures. Also the boundary conditions
glue in a natural way to a boundary condition $A\sharp B$ on
$D(\rho)$. The index of the corresponding operator satisfies
$$
\ix(\bar\pa_{A\sharp B})=n+\mu(\widehat{A\sharp B})=
\ix(\bar\pa_A)+\ix(\bar\pa_{B}).
$$
(Here the equality follows from the
fact that to obtain $\widehat{A\sharp B}$ from $\hat A$ and
$\hat B$, $n$ negative $\pi$-rotations are removed at their
common puncture.)
We call this type of gluing, {\em gluing at a puncture}.
\begin{lma}\label{lmaglue2}
For all sufficiently large $\rho>0$ there exist an exact sequence
\begin{equation}\label{eqglutv}
\begin{split}
&\begin{CD}
0 @>>> \krn(\bar\pa_{A\sharp B}) @>{\alpha_\rho}>>
\krn(\bar\pa_A)\oplus\krn(\bar\pa_{B})  @>{\beta_\rho}>>
\end{CD}\\
&\quad\quad
\begin{CD}
\cokrn(\bar\pa_A)\oplus\cokrn(\bar\pa_{B}) @>{\gamma_\rho}>>
\cokrn(\bar\pa_{A\sharp B}) @>>> 0
\end{CD}
\end{split}
\end{equation}
where $\alpha_\rho$, $\beta_\rho$, and $\gamma_\rho$ are defined similarly to
the maps in Remark~\ref{rmkmapsinseq}.
In particular, orientations of the determinants
spaces
$$
\Lambda^{\rm max}\Bigl(\krn(\bar\pa_A)
\oplus\krn(\bar\pa_{B})\Bigr)
\otimes\Lambda^{\rm max}
\Bigl(\cokrn(\bar\pa_A)
\oplus\cokrn(\bar\pa_{B})\Bigl)^\ast
$$
induces, via \eqref{eqglutv}, an orientation on
$$
\Lambda^{\rm max}\krn(\bar\pa_{A\sharp B})\otimes
\Lambda^{\rm max}\cokrn(\bar\pa_{A\sharp B})^\ast.
$$
\end{lma}
\begin{pf}
The proof of this result is similar to the proof of Lemma
\ref{lmaglue1}. However, in the present situation the proof is simpler since
the operator in the gluing region is Fredholm of index $0$ with
trivial kernel and cokernel. This is also the reason for the absence of
the $\R^n$ summand in the third term of the gluing sequence.
\end{pf}
\subsubsection{Associativity of orientations under gluing}
We show that orientations behave associatively with respect to both
gluings at punctures and at boundary points. We will often deal with
many vector spaces. For convenience we employ the following
notation: if $(V_1, V_2,\dots, V_n)$ are ordered vector spaces then
we denote their direct sum by a column vector containing the vector
spaces. That is,
$$
V_1\oplus\dots\oplus V_n=
\left[
\begin{matrix}
V_1\\
\vdots\\
V_n
\end{matrix}
\right].
$$
(Note that, when dealing with {\em oriented} vector spaces the
ordering of the terms in a direct sum is important.)
We first consider the case of two gluings at punctures.
Let $A\colon \pa D_{m}\to U(n)$, $B\colon \pa D_r\to U(n)$, and
$C\colon \pa D_s\to U(n)$ be trivialized transverse boundary conditions
for the $\bar\pa$-operator on three punctured disks. Assume that $A$
can be glued to $B$ at a puncture and that $B$ can be glued to $C$ at
another puncture. Performing these two gluings we obtain a boundary
condition $E\colon \pa D_{m+r+s-4}\to U(n)$.
Fix orientations $o_A$, $o_B$, and $o_C$ on
$\det(\bar\pa_A)$, $\det(\bar\pa_B)$, and $\det(\bar\pa_C)$,
respectively.
Via Lemma \ref{lmaglue2} these orientations induce orientations $o_{A\sharp B}$ and
$o_{B\sharp C}$ on
$\det(\bar\pa_{A\sharp B})$ and $\det(\bar\pa_{B\sharp C})$.
\begin{lma}\label{lmaassocpp}
Let $o'_E$ be the orientation on $\det(\bar\pa_E)$
induced via gluing (as in Lemma \ref{lmaglue2}) from the orientations
$o_{A\sharp B}$ and $o_C$ on
$\det(\bar\pa_{A\sharp B})$ and $\det(\bar\pa_C)$, respectively.
Let $o''_E$ be the orientation on $\det(\bar\pa_E)$
induced via gluing (as in Lemma \ref{lmaglue2}) from the orientations of
$o_{A}$ and $o_{B\sharp C}$ on
$\det(\bar\pa_{A})$ and $\det(\bar\pa_{B\sharp C})$, respectively. Then
$$
o'_E=o''_E=o_E,
$$
where $o_E$ is the orientation induced from the exact sequence
\begin{equation*}
\begin{split}
0\longrightarrow & \, \,  \begin{CD}
\krn(\bar\pa_{E}) @>{\alpha_\rho}>>
\left[\begin{matrix}
\krn(\bar\pa_A)\\
\krn(\bar\pa_B)\\
\krn(\bar\pa_C)
\end{matrix}\right]
@>{\beta_\rho}>>\end{CD}\\
&\begin{CD}
\left[\begin{matrix}
\cokrn(\bar\pa_A)\\
\cokrn(\bar\pa_B)\\
\cokrn(\bar\pa_C)
\end{matrix}\right]
@>{\gamma_\rho}>>
\cokrn(\bar\pa_{E})
\end{CD} \longrightarrow 0,
\end{split}
\end{equation*}
$\alpha_\rho$, $\beta_\rho$, and $\gamma_\rho$ are defined similarly to
the maps in Remark~\ref{rmkmapsinseq}.
\end{lma}
\begin{pf}
To see this consider the diagram
$$
\begin{CD}
\left[\begin{matrix}
\krn(\bar\pa_{A\sharp B})\\
\krn(\bar\pa_C)
\end{matrix}\right]
@>{\id}>>
\left[\begin{matrix}
\krn(\bar\pa_{A\sharp B})\\
\krn(\bar\pa_C)
\end{matrix}\right]
@>>>
\left[\begin{matrix}
\cokrn(\bar\pa_{A\sharp B})\\
\cokrn(\bar\pa_C)
\end{matrix}\right]
@>{\id}>>
\left[\begin{matrix}
\krn(\bar\pa_{A\sharp B})\\
\krn(\bar\pa_C)
\end{matrix}\right]
\\
@AAA @VVV @AAA @VVV\\
\krn(\bar\pa_{E}) @>{\alpha_\rho}>>
\left[\begin{matrix}
\krn(\bar\pa_A)\\
\krn(\bar\pa_B)\\
\krn(\bar\pa_C)
\end{matrix}\right]
@>{\beta_\rho}>>
\left[\begin{matrix}
\cokrn(\bar\pa_A)\\
\cokrn(\bar\pa_B)\\
\cokrn(\bar\pa_C)
\end{matrix}\right]
@>{\gamma_\rho}>>
\cokrn(\bar\pa_{E})
\end{CD},
$$
where the $0$'s at the ends of the horizontal rows are dropped as are
the $0$'s above the two middle terms in the upper horizontal
row. Completing the diagram with an arrow from the left middle entry
in the lower row to the right middle entry in the upper row we find the direct sum of
the gluing sequence of $\bar\pa_{A\sharp B}$ and the trivial sequence
for $\bar\pa_{C}$.
We check orientations. Let
$$
(\sigma,\bar\sigma)\in
\Lambda^{\rm max}(\krn(\bar\pa_E))
\times\Lambda^{\rm max}(\cokrn(\bar\pa_E))
$$
represent an orientation pair for $\det(\bar\pa_E)$. Pick a wedge of vectors
$\omega\wedge\gamma$ on the complement of its image in
$\krn(\bar\pa_{A\sharp B})\oplus\krn(\bar\pa_C)$ where $\omega$ is
a wedge of vectors in $\krn(\bar\pa_{A\sharp B})$ and $\gamma$ a wedge
of vectors in $\krn(\bar\pa_{C})$. Then the induced orientation on the
pair
$$
\left(
\left[\begin{matrix}
\krn(\bar\pa_{A\sharp B})\\
\krn(\bar\pa_{C})
\end{matrix}\right],
\left[\begin{matrix}
\cokrn(\bar\pa_{A\sharp B})\\
\cokrn(\bar\pa_{C})
\end{matrix}\right]
\right)
$$
is represented by
$$
(\sigma\wedge\omega\wedge\gamma,
\bar\gamma\wedge\bar\omega\wedge\bar\sigma),
$$
where $\bar\alpha$ is the image of the form $\alpha$ under the appropriate map in
Remark~\ref{rmkmapsinseq} with the vectors in the
opposite order, see Subsection~\ref{ssseclinalg}.
This latter orientation pair induces the orientation pair
$$
(\sigma\wedge\omega\wedge\gamma\wedge\eta,
\bar\eta\wedge\bar\gamma\wedge\bar\omega\wedge\bar\sigma).
$$
on the pair
$$
\left(
\left[\begin{matrix}
\krn(\bar\pa_{A})\\
\krn(\bar\pa_{B})\\
\krn(\bar\pa_{C})
\end{matrix}\right],
\left[\begin{matrix}
\cokrn(\bar\pa_{A})\\
\cokrn(\bar\pa_{B})\\
\cokrn(\bar\pa_{C})
\end{matrix}\right]
\right),
$$
where $\eta$ is any wedge of vectors in
$\krn(\bar\pa_A)\oplus\krn(\bar\pa_B)$ on the complement of the image
of $\omega$.
On the other hand the orientation pair $(\sigma,\bar\sigma)$ induces
the orientation pair
$$
(\sigma\wedge\phi,\bar\phi\wedge\bar\sigma)
$$
on the same pair of spaces. Since $\phi$ and $\omega\wedge\gamma\wedge\eta$
represents the same orientation on the complement of $\sigma$ if and
only if $\bar\phi$ represents the same orientation as
$\bar\eta\wedge\bar\gamma\wedge\bar\omega$ on the complement of
$\bar\sigma$ we see
that the orientation pairs agree and the two gluing sequences give the
same orientation on $\det(\bar\pa_E)$.
An identical argument shows that also the gluing sequence for
$\bar\pa_A$ and $\bar\pa_{B\sharp C}$ induces this orientation on
$\det(\bar\pa_E)$.
\end{pf}

Second we consider the mixed case.
Let $A\colon \pa D_{m}\to U(n)$, $B\colon \pa D_r\to U(n)$, and
$C\colon \pa D_s\to U(n)$ be trivialized transverse boundary conditions
for the $\bar\pa$-operator on three punctured disks. Assume that $A$
can be glued to $B$ at a boundary point and that $B$ can be glued to $C$ at
puncture. Performing these two gluings we obtain a boundary
condition $E\colon \pa D_{m+r+s-2}\to U(n)$. Fix orientations $o_A$, $o_B$, and $o_C$ on
$\det(\bar\pa_A)$, $\det(\bar\pa_B)$, and $\det(\bar\pa_C)$,
respectively.
Via Lemmas \ref{lmaglue1} and \ref{lmaglue2} these orientations induce
orientations $o_{A\sharp B}$ and
$o_{B\sharp C}$ on
$\det(\bar\pa_{A\sharp B})$ and $\det(\bar\pa_{B\sharp C})$, respectively.
\begin{lma}\label{lmaassocbp}
Let $o'_E$ be the orientation on $\det(\bar\pa_E)$
induced via gluing (as in Lemma \ref{lmaglue2}) from the orientations of
$o_{A\sharp B}$ and $o_C$ on
$\det(\bar\pa_{A\sharp B})$ and $\det(\bar\pa_C)$, respectively.
Let $o''_E$ be the orientation on $\det(\bar\pa_E)$
induced via gluing (as in Lemma \ref{lmaglue1}) from the orientations of
$o_{A}$ and $o_{B\sharp C}$ on
$\det(\bar\pa_{A})$ and $\det(\bar\pa_{B\sharp C})$, respectively. Then
$$
o'_E=o''_E=o_E,
$$
where $o_E$ is the orientation induced from the exact sequence
\begin{equation*}
\begin{split}
0\longrightarrow & \, \,  \begin{CD}
\krn(\bar\pa_{E}) @>{\alpha_\rho}>>
\left[\begin{matrix}
\krn(\bar\pa_A)\\
\krn(\bar\pa_B)\\
\krn(\bar\pa_C)
\end{matrix}\right]
@>{\beta_\rho}>>
\end{CD}\\
&\begin{CD}
\left[\begin{matrix}
\cokrn(\bar\pa_A)\\
\R^n\\
\cokrn(\bar\pa_B)\\
\cokrn(\bar\pa_C)
\end{matrix}\right]
@>{\gamma_\rho}>>
\cokrn(\bar\pa_{E})
\end{CD} \longrightarrow 0,
\end{split}
\end{equation*}
where $\alpha_\rho$, $\beta_\rho$, and $\gamma_\rho$ are the naturally
defined maps in a simultaneous gluing process (see Remark~\ref{rmkmapsinseq}).
\end{lma}
\begin{pf}
The proof is similar to the proof of Lemma \ref{lmaassocpp}.
\end{pf}

Finally we consider the case of two gluings at boundary points.
Let $A\colon \pa D_{m}\to U(n)$, $B\colon \pa D_r\to U(n)$, and
$C\colon \pa D_s\to U(n)$ be trivialized transverse boundary conditions
for the $\bar\pa$-operator on three punctured disks. Assume that $A$
can be glued to $B$ at a boundary point and that $B$ can be glued to
$C$ at a boundary point. Performing these two gluings we obtain a boundary
condition $E\colon \pa D_{m+r+s}\to U(n)$. Fix orientations $o_A$,
$o_B$, and $o_C$ on
$\det(\bar\pa_A)$, $\det(\bar\pa_B)$, and $\det(\bar\pa_C)$,
respectively.
Via Lemma \ref{lmaglue1} these orientations induce
orientations $o_{A\sharp B}$ and
$o_{B\sharp C}$ on
$\det(\bar\pa_{A\sharp B})$ and $\det(\bar\pa_{B\sharp C})$, respectively.
\begin{lma}\label{lmaassocbb}
Let $o'_E$ be the orientation on $\det(\bar\pa_E)$
induced via gluing (as in Lemma \ref{lmaglue1}) from the orientations of
$o_{A\sharp B}$ and $o_C$ on
$\det(\bar\pa_{A\sharp B})$ and $\det(\bar\pa_C)$, respectively.
Let $o''_E$ be the orientation on $\det(\bar\pa_E)$
induced via gluing (as in Lemma \ref{lmaglue1}) from the orientations of
$o_{A}$ and $o_{B\sharp C}$ on
$\det(\bar\pa_{A})$ and $\det(\bar\pa_{B\sharp C})$, respectively. Then
$$
o'_E=o''_E=o_E,
$$
where $o_E$ is the orientation induced from the exact sequence
\begin{equation*}
\begin{split}
0\longrightarrow & \, \,  \begin{CD}
\krn(\bar\pa_{E}) @>{\alpha_\rho}>>
\left[\begin{matrix}
\krn(\bar\pa_A)\\
\krn(\bar\pa_B)\\
\krn(\bar\pa_C)
\end{matrix}\right]
@>{\beta_\rho}>>\end{CD}\\
&\begin{CD}\left[\begin{matrix}
\cokrn(\bar\pa_A)\\
\R^n\\
\cokrn(\bar\pa_B)\\
\R^n\\
\cokrn(\bar\pa_C)
\end{matrix}\right]
@>{\gamma_\rho}>>
\cokrn(\bar\pa_{E})
\end{CD} \longrightarrow 0,
\end{split}
\end{equation*}
where $\alpha_\rho$, $\beta_\rho$, and $\gamma_\rho$ are naturally
defined maps in a simultaneous gluing process (see Remark~\ref{rmkmapsinseq}).
\end{lma}
\begin{pf}
The proof is similar to the proof of Lemma \ref{lmaassocpp}.
\end{pf}
\subsection{Canonical orientations and capping disks}\label{COCD}
We describe, following \cite{FOOO}, the canonical
orientation of the determinant bundle over the space of trivialized
Lagrangian boundary conditions on the closed ($0$-punctured) disk.
We also relate trivialized Lagrangian boundary conditions on a punctured
disk with trivialized boundary conditions on the closed disk using
capping disks.

\subsubsection{The canonical orientation of the determinant bundle over
trivialized boundary conditions on the disk}
The space of trivialized Lagrangian boundary conditions over the
$0$-punctured disk is the space $\Omega(U(n))$ of (free) loops in $U(n)$.
\begin{lma}\label{lmacanor}
The determinant bundle over $\Omega(U(n))$ is orientable. In fact, an
orientation of $\R^n$ induces an orientation of this bundle.
\end{lma}
\begin{pf}
Consider the fibration $\Omega(U(n))\to U(n)$ which is evaluation at
$1\in\pa D$. We find that if $Y$ is any component of $\Omega(U(n))$
then $\pi_1(Y)=\Z$ since $\pi_1(U(n))=\Z$ and $\pi_2(U(n))=0$. In
fact, a generator of $\pi_1(Y)$ can be described as follows. Fix any
element $A\colon\pa D\to U(n)$ in $Y$ and a loop $B(s)$ in $U(n)$
which starts and ends at $\id$ and which generates $\pi_1(U(n))$.
Then the loop $s\mapsto A\cdot B(s)$ generates
$\pi_1(\Omega(U(n)))$. Let $E\to Y$ be the determinant bundle. We
need to check that $w_1(E)=0$, where $w_1$ is the first
Stiefel-Whitney class. To see this note that the monodromy of the
orientation bundle of $E$ along the loop described is the following.
If $v_1,\dots,v_r$ is a basis in the kernel and $w_1,\dots,w_s$ one
in the cokernel of $\bar\pa_A$ then the monodromy is given by the
orientation pairs
$$
\Bigl(B(s)v_1\wedge\dots\wedge B(s)v_r,
B^\ast(s) w_1\wedge\dots\wedge B^\ast(s)w_s\Bigr).
$$
Hence it preserves orientation and $w_1(E)=0$.
To prove the second statement, we follow the argument in
\cite{FOOO}, \S 21. The trivialization $A$ on the boundary allow us to choose
a trivialization of the trivial bundle $\C^n$ near the boundary of
$\pa D$ such that the real plane field is constantly equal to
$\R^n\subset\C^n$.
Extending this trivialization inwards towards the center of the disk
and stretching the neck, we eventually split off a $\C P^1$ with a
complex vector bundle at $0\in D$. Using gluing arguments similar to
those above we find that orientations of the determinant line over the
disk with constant boundary condition together with the complex
orientation of the determinant line of the complex bundle over $\C P^1$
induces an orientation of $\det(\bar\pa_A)$. Since the determinant line
of the $\bar\pa$-operator is canonically (by evaluation at any point
in the boundary) identified with $\Lambda^n\R^n$ we get an induced
orientation as claimed.
\end{pf}
\begin{rmk}
As shown in \cite{FOOO} the orientations on the index space induced by
different trivializations (different ${\rm mod}\,\,2$ if $n=2$) are
different.
\end{rmk}
Fix an orientation on $\R^n$.
\begin{dfn}
The orientation of the determinant bundle over $\Omega(U(n))$ induced
from the fixed orientation on $\R^n$
will be called the {\em canonical orientation}.
\end{dfn}
Consider two trivialized boundary conditions $A$ and $B$ on the
$0$-punctured disk as in Lemma \ref{lmaglue1} and construct the glued boundary
condition $A\sharp B$.
\begin{lma}\label{lmagluecanor}
The gluing sequence in Lemma \ref{lmaglue1} induces from the canonical
orientations on $\det(\bar\pa_A)$ and on $\det(\bar\pa_B)$, and the
orientation $(-1)^{\frac12(n-1)(n-2)}$ times the standard orientation on
$\R^n$,  the canonical orientation on $\det(\bar\pa_{A\sharp B})$.
\end{lma}
\begin{pf}
As in the proof of Lemma \ref{lmacanor} we pick special trivializations of the
bundles near the boundary. Note that these trivializations glue in a
natural way. When pushing inwards we obtain two complex bundles over
two copies of $\C P^1$ sitting over the origins of the two disks
glued. Using a gluing argument similar to the one in Lemma
\ref{lmaglue1}, we see that the lemma will follow as soon as it is
proved for the disk with constant boundary conditions
$\R^n\subset\C^n$.
In this case
$\krn(\bar\pa_A)$, $\krn(\bar\pa_B)$, and $\krn(\bar\pa_{A\sharp B})$
are all isomorphic to $\R^n$ and the cokernels of all three operators
are trivial. Moreover, the fixed orientation on $\R^n$ gives the
canonical orientation on the determinants. The gluing sequence is
$$
\begin{CD}
0 @>>> \R^n @>{\alpha_\rho}>> \R^n\oplus\R^n
@>{\beta_\rho}>>\R^n @>>> 0,
\end{CD}
$$
where, as is easily seen, $\alpha_\rho(v)=(v,v)$. To determine
$\beta_\rho$, let $w\in\R^n$ be thought of as an element in the kernel
on the first disk. We then cut this constant function off by a cut-off
function $\phi$ which equals $1$ on the first disk and equals $0$ on
the second.
The cokernel of the $\bar\pa$-operator on the strip with positive
weights and constant boundary conditions is a complement of the
intersections of the kernel of the $L^2$-pairing with constant
functions. We can thus represent it as the subspace spanned by
constant functions in the strip which are cut-off {\em outside} the
supports of the cut-off functions $\phi$. Now,
$\bar\pa(\phi w)=(\frac{\pa\phi}{\pa\tau}+i\frac{\pa\phi}{\pa t}) w$,
and taking the $L^2$-inner product with the cut-off constant functions
we find that the sign of $\frac{\pa\phi}{\pa\tau}$ essentially
determines the map. Since the orientation of the $\tau$-axis is from
the first disk to the second we find that
$\beta_\rho(w_1,w_2)=-w_1+w_2$. It follows from our conventions in
Section \ref{ssseclinalg} that the standard
orientation on $\R^n\oplus\R^n$ and the orientation
$(-1)^{\frac12(n-1)(n-2)}$ times the standard orientation the second
$\R^n$ induces the standard orientation on the first $\R^n$. The lemma
follows.
\end{pf}
\begin{rmk}
Note that it follows from the above Lemma \ref{lmagluecanor} that
picking two
disks with constant boundary conditions and gluing these we get a disk
with constant boundary conditions and the orientation induced on the
determinant of the later from orientations on the determinants of the
former two is independent of the ordering of the former
two. The reason for this is that the change of order of the summands
is accompanied by a change of the map $\beta_\rho$ to $-\beta_\rho$. If
the dimension $n$ is even then both these changes preserve
orientation and, if $n$ is odd both reverse orientation. Thus in either
case the induced orientation is not affected by the change of order.
\end{rmk}
\subsubsection{Transverse Lagrangian subspaces and unitary operators}
Let $V_1$ and $V_2$ be two transverse Lagrangian subspaces of
$\C^n$. Let $\theta=(\theta_1,\dots,\theta_n)\in (0,\pi)^n$ be the
complex angle of $(V_1,V_2)$ and recall that there exists canonical
complex coordinates $(x_1+iy_1,\dots,x_n+iy_n)$ in $\C^n$ such that
$V_1=\spa(\pa_{x_1},\dots,\pa_{x_n})$ and
$V_2=(e^{i\theta_1}\pa_{x_1},\dots,e^{i\theta_n}\pa_{x_n})$.
Let $Q[V_1,V_2](s)$ be the $1$-parameter family of unitary
transformations of $\C^n$ given by the matrix
$$
\diag(e^{i\theta_1 s},\dots,e^{i\theta_n s})
$$
in canonical coordinates. Then $(V_1,V_2)\mapsto Q[V_1,V_2](s)$
defines a map of the space of pairs of transverse Lagrangian subspaces
in $\C^n$ into the path space of $U(n)$. (Note that $Q[V_1,V_2]$ is
independent of the choice of canonical coordinates.) The next lemma
shows that this map is continuous.
\begin{lma}
The $1$-parameter family of unitary transformations $Q[V_1,V_2](s)$
depends continuously on $(V_1,V_2)$.
\end{lma}
\begin{pf}
The space of Lagrangian subspaces of $\C^n$ transverse to $\R^n$
is identified with the space of symmetric linear matrices
$L$ via
$$
L\mapsto V(L)= \spa(L\pa_1+i\pa_1,\dots, L\pa_n+i\pa_n).
$$
If $\cot^{-1}\colon (-\infty,\infty)\to (0,\pi)$ is the inverse of
$\cot=\frac{\cos}{\sin}$ then it is easily checked that
$Q[\R^n,V(L)](s)$ is given by the
matrix
\begin{equation}\label{eqarccot}
e^{is\cot^{-1}(L)}
\end{equation}
in the standard basis of $\C^n$. The lemma is a straightforward
consequence of \eqref{eqarccot}.
\end{pf}
Recall, see Subsection~\ref{typesofpath}, that we subdivided the set of pairs of
transverse oriented Lagrangian subspaces into two subsets: even pairs
and odd pairs. We associate to such pairs, with complex angle
meeting certain conditions, $1$-parameter families of unitary
transformations. We define these unitary transformations by writing
their matrices in canonical coordinates. Note that the extra
conditions on the complex angle ensures that the transformations are
independent of the choice of canonical coordinates. Let
$\theta(V_1,V_2)=(\theta_1,\dots,\theta_n)$ denote the complex angle
of the ordered pair $(V_1,V_2)$.
\begin{itemize}
\item[({\bf ne})] If $(V_1,V_2)$ is even and if
$\theta_1<\theta_j$ for all $j=2,\dots,n$ then define
$$
R_{ne}[V_1,V_2](s)=
\diag\left(e^{i\theta_1s},e^{-i(\pi-\theta_2)s},
\dots,e^{-i(\pi-\theta_n)s}\right).
$$
\item[({\bf no})] If $(V_1,V_2)$ is odd and if
$\theta_1<\theta_2<\theta_j$, $j=3,\dots,n$ then define
$$
R_{no}[V_1,V_2](s)=
\diag\left(e^{i\theta_1 s}, e^{-i(2\pi-\theta_2)s},
e^{-i(\pi-\theta_3)s},\dots,e^{-i(\pi-\theta_n)s}\right).
$$
\item[({\bf pe})] If $(V_2,V_1)$ is even and if $\theta_n>\theta_j$,
$j=1,\dots,n-1$ then define
$$
R_{pe}[V_1,V_2](s)=
\diag\left(e^{-i(2\pi-\theta_1)s},\dots,e^{-i(2\pi-\theta_{n-1})s},
e^{-i(\pi-\theta_n)s}\right).
$$
\item[({\bf po})] If $(V_2,V_1)$ is odd and if
$\theta_n>\theta_{n-1}>\theta_j$, $j=1,\dots,n-2$ then define
$$
R_{po}[V_1,V_2](s)=
\diag\left(e^{-i(2\pi-\theta_1)s},\dots,e^{-i(2\pi-\theta_{n-2})s},
e^{-i(\pi-\theta_{n-1})s},e^{-i(\pi-\theta_n)s}\right).
$$
\end{itemize}
Let $\hat\pi=(\pi,\dots,\pi)$ and note that if the complex angle of
$(V_1,V_2)$ equals $\theta$ then the complex angle of
$(V_2,V_1)$ equals $\hat\pi-\theta$. It is then easily seen that if
$(V_1,V_2)$ satisfies {\bf ne} and
$(V_2,V_1)$ satisfies {\bf pe} then
$R_{pe}[V_2,V_1](1)\circ R_{ne}[V_1,V_2](1)=\id$. Similarly, if
$(V_1,V_2)$ satisfies {\bf no} and
$(V_2,V_1)$ satisfies {\bf po} then
$R_{po}[V_2,V_1](1)\circ R_{no}[V_1,V_2](1)=\id$.

\subsubsection{Stabilization of Lagrangian subspaces}
Let $V$ be an oriented Lagrangian subspace of $\C^n$ and let
$0\le\beta<\pi$. Let $L^u(\beta)$ and $L^l(\beta)$ be the oriented
Lagrangian subspaces of $\C^2$, with standard coordinates
$(x_1+iy_1,x_2+i y_2)$, given by the orienting basis
\begin{align*}
L^u(\beta) &= \left(e^{i\tfrac{\beta}{2}}\pa_{x_1},
e^{i\beta}\pa_{x_2}\right),\\
L^l(\beta) &= \left(e^{-i\tfrac{\beta}{2}}\pa_{x_1},
e^{-i\beta}\pa_{x_2}\right),\\
\end{align*}
Define the {\em upper $\beta$-stabilization} of $V$ to be the oriented
Lagrangian subspace $V^u$ of $\C^{n+2}=\C^n\times\C^2$ given by
$$
V^u(\beta)=V\times L^u(\beta).
$$
Define the {\em lower $\beta$-stabilization} of $V$ to be the oriented
Lagrangian subspace $V^l$ of $\C^{n+2}=\C^n\times\C^2$ given by
$$
V^l(\beta)=V\times L^l(\beta).
$$

Let $\bigl(V_1(\lambda),V_2(\lambda)\bigr)$,
$\lambda\in\Lambda$ be a continuous family of
transverse Lagrangian subspaces parameterized by a compact space
$\Lambda$. Let
$\theta(\lambda)=\bigl(\theta_1(\lambda),\dots,\theta_n(\lambda)\bigr)$
be the complex angle of $\bigl(V_1(\lambda),V_2(\lambda)\bigr)$. By
compactness of $\Lambda$ there exists $\beta>0$ such that
$2\beta<\theta_j(\lambda)$ and $\pi-2\beta>\theta_j(\lambda)$ for all
$\lambda\in\Lambda$ and
$j=1,\dots,n$. Fix such a $\beta>0$ and let
$\tilde V_1(\lambda)=[V_1(\lambda)]^u(\beta)$ and
$\tilde V_2(\lambda)=[V_2(\lambda)]^l(\beta)$.
Note that $\bigl(\tilde V_1(\lambda),\tilde V_2(\lambda)\bigr)$ is even (odd) if
and only if $\bigl(V_1(\lambda),V_2(\lambda)\bigr)$ is even (odd).
Moreover, by the choice of $\beta$,
$\bigl(\tilde V_1(\lambda),\tilde V_2(\lambda)\bigr)$ satisfies
the condition {\bf ne} ({\bf no}) if
$\bigl(V_2(\lambda),V_1(\lambda)\bigr)$ is even (odd) and
$\bigl(\tilde V_2(\lambda),\tilde V_1(\lambda)\bigr)$ satisfies
the condition {\bf pe} ({\bf po}) if
$\bigl(V_1(\lambda),V_2(\lambda)\bigr)$ is even (odd), for
all $\lambda\in\Lambda$. Thus we can construct the corresponding
$\Lambda$-families of unitary operators.
\begin{lma}
Let $\bigl(V_1(\lambda),V_2(\lambda)\bigr)$ and $\beta>0$ be as
above. If $\bigl(V_1(\lambda),V_2(\lambda)\bigr)$ is even then the
families of unitary operators
$R_{ne}[\tilde V_1(\lambda),\tilde V_2(\lambda)](s)$ and
$R_{pe}[\tilde V_2(\lambda),\tilde V_1(\lambda)](s)$ depend continuously
on $\lambda\in\Lambda$. If
$\bigl(V_1(\lambda),V_2(\lambda)\bigr)$ is odd then the families of
unitary operators
$R_{no}[\tilde V_1(\lambda),\tilde V_2(\lambda)](s)$ and
$R_{po}[\tilde V_2(\lambda),\tilde V_1(\lambda)](s)$ depend continuously
on $\lambda\in\Lambda$.
\end{lma}
\begin{pf}
This is a straightforward consequence of \eqref{eqarccot}.
\end{pf}
\subsubsection{Determinant bundles over stabilized Lagrangian
subspaces}
Let $\bigl(V_1(\lambda),V_2(\lambda)\bigr)$, $\lambda\in\Lambda$ be a
continuous family of
transverse Lagrangian subspaces of $\C^n$ parameterized by a
compact simply connected space $\Lambda$. Fix $\beta>0$ small enough
and consider the stabilized family
$\bigl(\tilde V_1(\lambda),\tilde V_2(\lambda)\bigr)$ of transverse
Lagrangian subspaces in $\C^{n+2}$.
Assume that $\tilde V_1(\lambda)$ and $\tilde V_2(\lambda)$ are
equipped with positively oriented frames $X_1(\lambda)$ and
$X_2(\lambda)$ which vary continuously with $\lambda\in\Lambda$.
We associate to this family two families of trivialized
Lagrangian boundary conditions on the $1$-punctured disk.
To simplify notation, if $\bigl(V_1(\lambda),V_2(\lambda)\bigr)$ is
even then let $\ast=e$ and if $\bigl(V_1(\lambda),V_2(\lambda)\bigr)$
is odd then let $\ast=o$. Note that
$$
R_{n\ast}[\tilde V_1(\lambda),\tilde V_2(\lambda)](1)X_1(\lambda)
$$
is a framing of $\tilde V_2(\lambda)$. Hence there exists
$\alpha(\lambda)\in SO(n+2)$ such that
$$
R_{n\ast}[\tilde V_1(\lambda),\tilde V_2(\lambda)](1)X_1(\lambda)
=X_2(\lambda)\cdot\alpha(\lambda).
$$
Since $\Lambda$ is simply connected the map
$\alpha\colon \Lambda\to SO(n+2)$ lifts to a map
$\tilde\alpha\colon \Lambda\to PSO(n+2)$, where $PSO(n+2)$ is
the space of paths in $SO(n+2)$ with initial endpoint at the identity
matrix and which projects to $SO(n+2)$ by evaluation at the final
endpoint. Pick such a lift. Identify $\pa D_1$ with $[0,1]$. Define
two families of trivialized boundary
conditions for the $\bar\pa$-operator on $D_1$ as follows
\begin{align}
A_{n\ast}[\lambda](s) &=
R_{n\ast}[\tilde V_1(\lambda),\tilde V_2(\lambda)]
X_1(\lambda)\cdot\tilde\alpha[\lambda](s),\\
A_{p\ast}[\lambda](s) &=
R_{p\ast}[\tilde V_2(\lambda),\tilde V_1(\lambda)]
X_2(\lambda)\cdot\tilde\alpha^{-1}[\lambda](s),
\end{align}
where $\tilde\alpha^{-1}[\lambda](s)$ is the inverse of the matrix
$\tilde\alpha[\lambda](s)$.
Let $p$ be the puncture on $D_1$ and note that
$A_{n\ast}[\lambda](p^\pm)=A_{p\ast}[\lambda](p^\mp)$ so that these
boundary conditions can be glued. Let $\bar\pa_{{n\ast},\lambda}$ and
$\bar\pa_{p\ast,\lambda}$ denote the $\bar\pa$-operators with boundary
conditions $A_{n\ast}[\lambda]$ and $A_{p_\ast}[\lambda]$,
respectively. Then the bundles
$E_{n\ast}=\det(\bar\pa_{n\ast,\lambda})\to\Lambda$ and
$E_{p\ast}=\det(\bar\pa_{p\ast,\lambda})\to\Lambda$ are orientable since
$\Lambda$ is simply connected. We orient them as follows. Let
$0\in\Lambda$. Pick an orientation of
$\det(\bar\pa_{n\ast,0})$. Together with the canonical orientation of
the $\bar\pa$-operator on the $0$-punctured disk and Lemma~\ref{lmaglue2} this
orientation determines an orientation on
$\det(\bar\pa_{p\ast,0})$. Since the bundles
$E_{n\ast}$ and $E_{p\ast}$ are orientable the orientations of
$\det(\bar\pa_{n\ast,0})$ and $\det(\bar\pa_{p\ast,0})$ induce
orientations on $\det(\bar\pa_{n\ast,\lambda})$ and
$\det(\bar\pa_{p\ast,\lambda})$, respectively, for any
$\lambda\in\Lambda$.
\begin{lma}
The orientations on $\det(\bar\pa_{n\ast,\lambda})$ and
$\det(\bar\pa_{p\ast,\lambda})$ as defined above glue to the
canonical orientation on the $0$-punctured disk.
\end{lma}
\begin{pf}
After adding a finite dimensional vector space we may assume that all
operators are surjective. The lemma then follows from properties of
the direct sum operation on vector bundles.
\end{pf}
We call the operators $\bar\pa_{ne,\lambda}$, $\bar\pa_{pe,\lambda}$,
$\bar\pa_{no,\lambda}$, and $\bar\pa_{po,\lambda}$ which arises as
above {\em capping operators} and we call an orientation pair on
$E_{n\ast}$, $E_{p\ast}$ with the properties above a pair of capping
orientations.
\subsubsection{Capping orientations of trivialized boundary
conditions}\label{caportrivbdcond}
Let $\Lambda$ be any compact simply connected space and let
$\bigl(V_1^j(\lambda),V_2^j(\lambda)\bigr)$, $j=0,\dots,m$ be
$\Lambda$-families of transverse Lagrangian
subspaces. Construct as in the previous section the stabilizations
$\bigl(\tilde V_1^j(\lambda),\tilde V_2^j(\lambda)\bigr)$ and assume
that these families are equipped with positively oriented
frames $X_1(\lambda)$ and $X_2(\lambda)$, respectively. Construct the families
$\bar\pa^j_{ne,\lambda}$, $\bar\pa^j_{pe,\lambda}$,
$\bar\pa^j_{no,\lambda}$, and $\bar\pa^j_{po,\lambda}$ of operators on
the $1$-punctured disks with oriented determinant bundles as there.
Let $X$ be any topological space and consider an
$(X\times\Lambda)$-family of oriented Lagrangian
boundary conditions $a[x,\lambda]\colon\pa D_{m+1}\to\Lag(n)$,
$(x,\lambda)\in X\times\Lambda$ on the
$(m+1)$-punctured disk with the following properties.
For all $(x,\lambda)$, if $j\ge 1$ then
$a[x,\lambda](p_j^+)=V^1_j(\lambda)$,
$a[x,\lambda](p_j^-)=V_2^j(\lambda)$,
$a[x,\lambda](p_0^+)=V_2^0(\lambda)$,
and $a[x,\lambda](p_0^-)=V_1^0(\lambda)$.
Let $E\to X\times\Lambda$ denote the determinant bundle of
$\bar\pa_{a[x,\lambda]}$. We construct, under some additional
trivialization conditions, the {\em capping orientation} of
$E$ in the following way.
First stabilize the family $a[x,\lambda]\colon\pa D_m\to\Lag(n)$ to
$\tilde a[x,\lambda]\colon\pa D_m\to\Lag(n+2)$ by letting
$$
\tilde a[x,\lambda](\zeta)=
a[x,\lambda](\zeta)\times
\spa\left(e^{i\tfrac{\gamma}{2}(\zeta)}\pa_{x_1},
e^{i\gamma(\zeta)}\pa_{x_2}\right),\quad\zeta\in\pa D_m,
$$
where $\C^{n+2}=\C^n\times\C^2$, $\C^2=\{(x_1+iy_1,x_2+iy_2)\}$, and
where $\gamma\colon\pa D_m\to [-\beta,\beta]$ satisfies
$\gamma(\zeta)=\beta$ in neighborhoods of $p_j^+$, $j\ge 1$ and
$p_0^-$, and $\gamma(\zeta)=-\beta$ in neighborhoods of $p_j^-$ and
$p_0^+$. Noting that the $\bar\pa$-operator on $D_m$ with boundary
conditions
$$
\zeta\mapsto \diag\left(e^{i\tfrac{\gamma}{2}(\zeta)}\pa_{x_1},
e^{i\gamma(\zeta)}\pa_{x_2}\right)
$$
has both trivial kernel and trivial cokernel we find that the kernels
and cokernels of $\bar\pa_{a[x,\lambda]}$ and
$\bar\pa_{\tilde a[x,\lambda]}$ are canonically isomorphic.
Thus to orient the determinant bundle of $\bar\pa_{a[x,\lambda]}$ it
suffices to orient the determinant bundle of
$\bar\pa_{\tilde a[x,\lambda]}$. We find such an orientation in the
case when the boundary conditions $\tilde a[x,\lambda]$ are equipped
with a certain type of trivialization. Thus, assume that
the boundary conditions $\tilde a[x,\lambda]\colon \pa D_m\to\Lag(n+2)$
are trivialized, i.e. they are represented by
$\tilde A[x,\lambda]\colon \pa D_m\to U(n+2)$. Assume moreover that
this trivialization satisfies the following conditions
for all $(x,\lambda)$, if $j\ge 1$ then
$\tilde A[x,\lambda](p_j^+)=X^j_1(\lambda)$,
$\tilde A[x,\lambda](p_j^-)=X_2^j(\lambda)$,
$\tilde A[x,\lambda](p_0^+)=X_2^0(\lambda)$,
and $\tilde A[x,\lambda](p_0^-)=X_1^0(\lambda)$.

Glue to $\bar\pa_{\tilde A[x,\lambda]}$ first the capping operator
$\bar\pa_{pe,\lambda}^0$ ($\bar\pa_{po,\lambda}^0$) at $p_0$ if the
puncture $p_0$ is even (odd). Then glue to
$\bar\pa_{\tilde A[x,\lambda]}$ the operators $\bar\pa_{ne,\lambda}^j$
($\bar\pa_{ne,\lambda}^j$) to $\tilde A[x,\lambda]$ at its
even (odd) negative puncture $p_j$, $j\ge 1$, {\em in the order opposite} to that
induced by the boundary orientation of $\pa D_{m+1}$. (We glue in the opposite order
so that the Leibniz rule works out appropriately, see Section~\ref{sec:LegCH}.)
We obtain in this way a
trivialized boundary condition $\hat A[x,\lambda]\colon \pa D\to
U(n+2)$ on the closed disk. Lemma \ref{lmaglue2} and repeated
application of Lemma \ref{lmaassocpp} gives the gluing sequence, where
we write $\bar\pa_{q^\pm,\lambda}$ for the capping operator at the
puncture $q$ at $\lambda\in\Lambda$ and where the sign indicates the
sign of the puncture,
\begin{equation}\label{eqcapseq}
\begin{split}
0\longrightarrow & \, \,  \begin{CD}
\krn(\bar\pa_{\hat A[x,\lambda]}) @>>>
\left[\begin{matrix}
\krn(\bar\pa_{p_1^-,\lambda})\\
\vdots\\
\krn(\bar\pa_{p_m^-,\lambda})\\
\krn(\bar\pa_{p_0^+,\lambda})\\
\krn(\bar\pa_{\tilde A[x,\lambda]})
\end{matrix}\right] @>>>
\end{CD}\\
&\begin{CD}
\left[\begin{matrix}
\cokrn(\bar\pa_{p_1^-,\lambda})\\
\vdots\\
\cokrn(\bar\pa_{p_m^-,\lambda})\\
\cokrn(\bar\pa_{p_0^+,\lambda})\\
\cokrn(\bar\pa_{\tilde A[x,\lambda]})
\end{matrix}\right] @>>>
\cokrn(\bar\pa_{\hat A[x,\lambda]})
\end{CD} \longrightarrow 0.
\end{split}
\end{equation}
We give the determinant of $\bar\pa_{\tilde A[x,\lambda]}$ the unique
orientation $o$ which together with the chosen orientations for
the capping disks, via the
the gluing sequence \eqref{eqcapseq}, give the canonical orientation on
$\det(\bar\pa_{\hat A[x,\lambda]})$. It is clear that the orientation
so defined gives an
orientation of the determinant bundle over
$\det(\bar\pa_{\tilde A[x,\lambda]})\to X\times\Lambda$ and thus by
the above mentioned isomorphism also the bundle
$E=\det(\bar\pa_{a[x,\lambda]})$ gets oriented.
\begin{dfn}\label{dfncap}
We call
the orientation of the determinant bundle $E\to X\times\Lambda$ the
{\em capping orientation}.
\end{dfn}
Note that the same construction can be applied when $X\times\Lambda$
is replaced by a locally trivial fibration $Y\to\Lambda$.
\subsubsection{Kernels and cokernels of capping operators}
Let $(V_1,V_2)$ be a pair of transverse Lagrangian subspaces in
$\C^{n-2}$ and let $(\tilde V_1,\tilde V_2)$ in $\C^n$ be its
stabilization. Note
that in canonical coordinates of $(\tilde V_1,\tilde V_2)$ the boundary
conditions of the capping operators constructed from
$(\tilde V_1,\tilde V_2)$ are split and we may determine the
dimensions of the kernel and
cokernel from properties of the classical Riemann-Hilbert problem.
To this end let $\theta=(\theta_1,\dots,\theta_n)$,
$0<\theta_1\le\theta_2\le\dots\le\theta_n<\pi$, be the complex
angle of $(\tilde V_1,\tilde V_2)$. We think of the
$1$-punctured disk $D_1$ as of the unit disk in $\C$ punctured at $1$
boundary parameterized by $e^{is}$, $0\le s\le 2\pi$.
\begin{itemize}
\item
The {\em negative even} boundary condition is, in
canonical coordinates, given by
$$
s\mapsto\diag\left(e^{i\theta_1 \tfrac{s}{2\pi}},
e^{-i(\pi-\theta_2)\tfrac{s}{2\pi}},
\dots,e^{-i(\pi-\theta_n)\tfrac{s}{2\pi}}\right).
$$
Thus,
\begin{align*}
&\ix(\bar\pa_{ne})=n-(n-1)=1,\\
&\dim(\krn(\bar\pa_{ne}))=1,\\
&\dim(\cokrn(\bar\pa_{ne}))=0.
\end{align*}
The kernel is spanned by the function $w(z)=(w_1(z),\dots,w_n(z))$,
where
\begin{align*}
w_1(z)&=(i(z-1))^{\tfrac{\theta_1}{\pi}},\\
w_j(z)&=0,\quad j>1.
\end{align*}
\item
The {\em negative odd} boundary condition is, in
canonical coordinates, given by
$$
s\mapsto\diag\left(e^{i\theta_1 \tfrac{s}{2\pi}},
e^{-i(2\pi-\theta_2)\tfrac{s}{2\pi}},
e^{-i(\pi-\theta_3)\tfrac{s}{2\pi}},
\dots,e^{-i(\pi-\theta_n)\tfrac{s}{2\pi}}\right).
$$
Thus,
\begin{align*}
&\ix(\bar\pa_{no})=n-n=0,\\
&\dim(\krn(\bar\pa_{no}))=1,\\
&\dim(\cokrn(\bar\pa_{no}))=1.
\end{align*}
The kernel spanned by the function $w(z)=(w_1(z),\dots,w_n(z))$,
where
\begin{align*}
w_1(z)&=(i(z-1))^{\tfrac{\theta_1}{\pi}},\\
w_j(z)&=0,\quad j>1.
\end{align*}
The cokernel is spanned by the function
$u(z)=(u_1(z),\dots,u_n(z))$, where
\begin{align*}
u_2(z)&=(i(\bar z-1))^{\tfrac{\theta_2}{\pi}},\\
u_j(z)&=0,\quad j\ne 2,
\end{align*}
where we view this function as a linear functional on the target space
of $\bar\pa$ via the $L^2$-pairing and where the cokernel is a one
dimensional subspace complementary to the kernel of this
functional.
\item
The {\em positive even} boundary condition is, in
canonical coordinates, given by
$$
s\mapsto\diag\left(e^{-i(2\pi-\theta_1)\tfrac{s}{2\pi}},
\dots,e^{-i(2\pi-\theta_{n-1})\tfrac{s}{2\pi}},
e^{-i(\pi-\theta_n)\tfrac{s}{2\pi}}\right).
$$
Thus,
\begin{align*}
&\ix(\bar\pa_{ne})=n-2(n-1)-1=1-n,\\
&\dim(\krn(\bar\pa_{ne}))=0,\\
&\dim(\cokrn(\bar\pa_{ne}))=n-1.
\end{align*}
The cokernel is spanned by the functions
$u(z)=(u_1(z),\dots,u_n(z))$, where for $k\ne n$
\begin{align*}
u_k(z)&=(i(\bar z-1))^{\tfrac{\theta_k}{\pi}},\\
u_j(z)&=0,\quad j\ne k.
\end{align*}
\item
The {\em positive odd} boundary condition is, in
canonical coordinates, given by
$$
s\mapsto\diag\left(e^{-i(2\pi-\theta_1)\tfrac{s}{2\pi}},
\dots,e^{-i(2\pi-\theta_{n-2})\tfrac{s}{2\pi}}
,e^{-i(\pi-\theta_{n-1})\tfrac{s}{2\pi}},
e^{-i(\pi-\theta_n)\tfrac{s}{2\pi}}\right).
$$
Thus,
\begin{align*}
&\ix(\bar\pa_{ne})=n-2(n-2)-2=2-n,\\
&\dim(\krn(\bar\pa_{ne}))=0,\\
&\dim(\cokrn(\bar\pa_{ne}))=n-2.
\end{align*}
The cokernel is spanned by the functions
$u(z)=(u_1(z),\dots,u_n(z))$, where for $k< n-1$
\begin{align*}
u_k(z)&=(i(\bar z-1))^{\tfrac{\theta_k}{\pi}},\\
u_j(z)&=0,\quad j\ne k.
\end{align*}
\end{itemize}
\subsection{Orientations of moduli spaces}\label{OMS}
We orient moduli spaces of holomorphic disks with boundary on a
generic Legendrian submanifold equipped with a spin structure. The
orientation is
obtained by comparing capping orientations with fixed
orientations of spaces of automorphisms and of conformal structures.

\subsubsection{Orientations of spaces of automorphisms and conformal
structures}
Let $D_{m+1}$ denote the unit disk $D\subset \C$ with $1$ positive, and $m$
negative punctures on the boundary, $m\ge 1$. Let $p_0$ be the
positive puncture and let $\{p_1,\dots,p_m\}$ be the negative punctures.
As mentioned above the positive puncture $p_0$ and the orientation of
$\pa D$ induces an ordering of the negative punctures
$(p_1,p_2,\dots,p_m)$.
Let $\conf_m$ denote the space of conformal structures on $D_m$. If
$m\le 3$ then $\conf_m$ is a one-point space. For $m<3$ let $\aut_m$
denote the group of conformal automorphisms.
We orient $\conf_{m+1}$, $m>2$ in the following way. Let $D_m$ with
conformal structure $\kappa$ be represented by a disk with positive
puncture $p_0$ and ordered negative punctures $(p_1,\dots,p_m)$.
Then any conformal structure $\kappa'$ in a neighborhood of $\kappa$
can be represented uniquely by a disk with positive puncture at
$p_0$ its first two negative punctures at $p_1$ and $p_2$ and with
the rest of its negative punctures at $(p_3',\dots,p_m')$. Thus, the
tangent space of $\conf_{m+1}$ at $\kappa$ is identified with the
direct sum of the tangent spaces of $\pa D$ at $p_3,\dots,p_m$. We
orient $\conf_m$ by declaring the oriented basis
$$
\Bigl\{\pa_{p_3}\oplus 0\dots\oplus 0,\,\, 0\oplus\pa_{p_2}\oplus
0\dots\oplus 0,\,\dots\,, 0\oplus\dots\oplus 0\oplus
\pa_{p_m}\Bigr\},
$$
where $\pa_{p_j}$ is the positive unit tangent to $\pa D$ at $p_j$,
to be positive basis in $T_\kappa\conf_m$.
We orient the one-point space $\conf_3$ by declaring it positively
oriented.
Next consider $\aut_1$ and $\aut_2$. To orient $\aut_1$, consider
$D_1$ punctured at $p$. Pick two points $q_1$ and $q_2$ in $\pa D_1$
a small distance from $p$. The tangent space of $\aut_1$ at the
identity is the $2$-dimensional space of holomorphic vector fields
on $D$ tangent to $\pa D$ along the boundary and vanishing at $p$.
Evaluation of such a vector fields at $q_1$ and $q_2$ gives a map
$$
T_{\id}\aut_1\to T_{q_1}D\oplus T_{q_2}D.
$$
We use this map to orient $T_{\id}\aut_1$ and the group structure of
$\aut_1$ to orient $\aut_1$.
To orient $\aut_2$, consider $D_2$ as $D$ with positive puncture at
$1$ and negative at $-1$. Pick a point $q$ in the {\em lower}
hemisphere of the two into which $\pa D$ is subdivided by $-1$ and
$1$ and orient $T_{\id}\aut_1$ by evaluation at $q$ as above.

\subsubsection{Stably trivialized boundary conditions and Legendrian
submanifolds with spin structures}\label{cappingrulegeq3}
Let $M$ be an orientable manifold of dimension $n$. Let ${\tilde
T}M=TM\oplus\R^2$ denote the stabilized tangent bundle of $M$, where
$\R^2$ denotes the trivial bundle over $M$. Fix some
triangulation of $M$. Then a spin-structure on $M$ can be viewed as a
trivialization of ${\tilde T}M$ restricted to the 1-skeleton that
extends to the 2-skeleton. Note that the extension to the 2-skeleton
is homotopically unique if it exists since $\pi_2(SO(n+2))=0$. For the
same reason, any trivialization over the 2-skeleton
automatically extends over the 3-skeleton but this extension is in
general not homotopically unique.
Let $L$ be an oriented manifold equipped with a spin structure and
let $\Lambda$ be a compact simply connected space. (In fact, for our
applications $\Lambda=\{{\rm point}\}$ and $\Lambda=[0,1]$ are
sufficient.) Let $\Phi_\lambda\colon L\to\C^n\times\R$,
$\lambda\in\Lambda$ be a family of {\em chord generic} Legendrian
embeddings. The chord genericity implies that we get a continuous
family of Reeb chords $c_1(\lambda),\dots,c_k(\lambda)$ of
$\Phi_\lambda$ and that their endpoints vary continuously with
$\lambda$ in $L$. Fix $0\in\Lambda$ and choose a family of
diffeomorphisms $\phi_\lambda\colon L\to L$ such that
$\phi_\lambda(c_j(0)^\pm)=c_j(\lambda)^\pm$ for all $j=1,\dots,k$.
For each Reeb chord $c(0)$ of $\Phi_0$, fix a capping path
$\gamma_{c(0)}\subset L$ connecting its upper end point
to the lower one, see \cite{ees1}, Section~2.3. Fix a
triangulation $\Delta=\Delta^{(0)}\cup\dots\cup\Delta^{(n)}$ of $L$,
where $\Delta^{(j)}$ denotes the $j$-skeleton of $\Delta$, such that
each Reeb chord endpoint $c_j^{\pm}(0)$ lies in $\Delta^{(0)}$ and each
capping path lies in $\Delta^{(1)}$.
As mentioned above the spin structure gives a trivialization of the
restriction of ${\tilde T}L$ to $\Delta^{(3)}$. Fix such a
trivialization. As $\lambda\in\Lambda$ varies we move the
triangulation and capping paths by $\phi_\lambda$ and the
trivialization by the bundle isomorphism ${\tilde T}M\to{\tilde T}M$
covering $\phi_\lambda$ given by
$$
\left(
\begin{matrix}
d\phi_\lambda & 0\\
0 & \id
\end{matrix}
\right).
$$
Using the notation from \cite{ees2}, Section 4, let
$\cand_{A,\Lambda}(a;b_1,\dots,b_m)$ denote the
space of candidate maps.
Such a trivialization enables us to construct the capping
orientation of the determinant bundle over the space of linearized
boundary conditions over $\cand_{A;\Lambda}(a;b_1,\dots,b_m)$
in the following way. Consider first the trivial $\C^n$-bundle over
$L$ which is the the pull-back $(\Pi_\C\circ\Phi_\lambda)^\ast
T\C^n$ and its stabilization $\C^n\times\C^2$. Associate to each
point $q\in L$ the Lagrangian subspace
$$
\Pi_\C\circ d\Phi_\lambda (T_p L)\oplus
\spa\left(e^{i\tfrac{\alpha_\lambda(p)}{2}}\pa_{x_1},
e^{i\alpha_\lambda(p)}\pa_{x_2}\right),
$$
where $\alpha_\lambda\colon L\to[-\beta,\beta]$ is a family of smooth
functions such that $\alpha_\lambda(c_j^\pm(\lambda))=\pm\beta$ for
all Reeb chords $c_j$ and where $\beta$ is such that $2\beta$ is
smaller than any component of any complex angle at any Reeb chord and
$\pi-2\beta$ is larger than such components.
Let
$(w,h,\kappa)\in\cand_{A,\Lambda}(a;b_1,\dots,b_m)$.
Then the
restriction of $(w,h)$ to the part of the boundary of $D_{m+1}$
lying between $p_j$ and $p_{j+1}$ is a path in $L$ connecting two
Reeb chord endpoints and we obtain from the field of Lagrangian
subspaces just defined a family of Lagrangian subspaces over
$\cand_{A;\Lambda}(a;b_1,\dots,b_m)$ which is the stabilization of
the tangent plane family and which satisfies the conditions in
Section~\ref{caportrivbdcond}. Pick a homotopy, fixing endpoints, of
this path to a path which lies in $\Delta^{(1)}$. Then the
trivialization $\eta$ of ${\tilde T}L$ along $\Delta^{(1)}$ induces
a trivialization of $(w,h)^\ast{\tilde T}L$ on the corresponding
parts. This in turn induces a trivialization of the stabilized
Lagrangian boundary condition and we define the capping orientation
of the determinant bundle of the linearized $\bar\pa$-operator over
$\cand_{A,\Lambda}(a;b_1,\dots,b_m)$ as in Definition \ref{dfncap}.
We must check that this orientation is well-defined. Choosing a
different homotopy to some path in $\Delta^{(1)}$, the two end paths
in $\Delta^{(1)}$ can be connected with a homotopy in
$\Delta^{(2)}$. Using the capping orientation over this homotopy proves
that the orientation is well-defined. We call also the orientation of
$\det(\bar\pa_\kappa)$ over $\cand_A(a;b_1,\dots,b_m)$ the
{\em capping orientation}.

\subsubsection{Moduli space orientations}
Using the following lemma we orient all moduli spaces of
holomorphic disks. As in \cite{ees2} we let $\Gamma$ denote the full
$\bar\pa$-operator.
\begin{lma}\label{lmadGamma}
Let $d\Gamma$ denote the full linearization of the $\bar\pa$-operator
at some holomorphic $(w,h,\kappa)\in\cand_A(a;b_1,\dots,b_m)$.
The determinant bundle over
$\cand_A(a;b_1,\dots,b_m)$ with fiber $\det(d\Gamma)$ is canonically
isomorphic to the tensor product of the determinant bundle of
$\bar\pa_\kappa$ with the highest exterior power of the tangent bundle
to $\conf_{m+1}$.
\end{lma}
\begin{pf}
Let $(w,h,\kappa)\in\cand_A$. We must show that (at $(w,h,\kappa)$)
\begin{align*}
&\Lambda^{\rm max}\krn(d\Gamma)\otimes\left(\Lambda^{\rm
max}\cokrn(d\Gamma)\right)^\ast=\\
&\Lambda^{\rm max}\krn(\bar\pa_\kappa)\otimes\left(\Lambda^{\rm
max}\cokrn(\bar\pa_\kappa)\right)^\ast\otimes \Lambda^{\max}
T\conf_{m+1}.
\end{align*}
Note that
$$
T\cand_A= T\cand_A(\kappa)\oplus T\conf_{m+1}= V\oplus\R^N,
$$
and that
$$
d\Gamma=\bar\pa_\kappa\oplus\psi,
$$
for some map
$$
\psi\colon
T_\kappa\conf_{m+1}=\R^N\to\sblv_1({T^\ast}^{0,1}D_{m+1}
\otimes\C^n)=W.
$$
Let $\bar\psi\colon \R^N\to\cokrn(\bar\pa_\kappa)$ be the map induced
by projection. Then the following sequence is exact
$$
0\to\krn(\bar\psi)\to\R^N\to\cokrn(\bar\pa_\kappa)\to\cokrn(\bar\psi)\to 0,
$$
and induces a canonical isomorphism
$$
\Lambda^{\rm max}\krn(\bar\psi)\otimes\Lambda^{\rm max}\cokrn(\bar\pa_\kappa)
=
\Lambda^{\rm max}\R^N\otimes\Lambda^{\rm max}\cokrn(\bar\psi).
$$
On the other hand $\cokrn(d\Gamma)=\cokrn(\bar\psi)$ and
$\krn(d\Gamma)=\krn(\bar\pa_\kappa)\oplus\krn(\bar\psi)$.
Hence
\begin{align*}
&\Lambda^{\rm max}\krn(d\Gamma)\otimes\left(\Lambda^{\rm max}\cokrn(d\Gamma)\right)^\ast=\\
&\Lambda^{\rm max}\krn(\bar\pa_\kappa)\otimes
\Lambda^{\rm max}\krn(\bar\psi)\otimes
\left(\Lambda^{\rm max}\cokrn(\bar\psi)\right)^\ast=\\
&
\Lambda^{\rm max}\krn(\bar\pa_\kappa)\otimes
\left(\Lambda^{\rm max}\cokrn(\bar\pa_\kappa)\right)^\ast\otimes
\Lambda^{\rm max}\R^N.
\end{align*}
\end{pf}
Let $(w,h,\kappa)\in\cand_A(a;b_1,\dots,b_m)$ be a transversely cut
out holomorphic disk. Using Lemma \ref{lmadGamma} we orient the moduli
space to which $(w,h,\kappa)$ belongs as follows.
\begin{itemize}
\item[{\rm (a)}] If $m\le 1$ then $d\Gamma$ simply agrees with the
ordinary $\bar\pa$-operator and the transversality condition implies
that this operator is surjective. Moreover,
$v\mapsto dw\cdot v$ for $v\in T\aut_{m+1}$ gives an injection
$T\aut_{m+1}\subset\krn(\bar\pa)$. The quotient $\krn(\bar\pa)/
T\aut_{m+1}$
can be identified with the tangent space to the moduli space
$T\M$. The capping orientation on $\det(\bar\pa)$ together with the
orientation on $\aut_{m+1}$ thus gives an orientation of the moduli space.
\item[{\rm (b)}] If $m\ge 2$ then an orientation of
$\det(\bar\pa_\kappa)$ and $T\conf_{m+1}$ together give an orientation of
$\det(d\Gamma)$. The surjectivity assumption implies that
$\cokrn(d\Gamma)$ is trivial hence we get an orientation on
$\krn(d\Gamma)$ which is the tangent space of the moduli space.
\end{itemize}
\begin{rmk}
An oriented connected $0$-dimensional manifold is a point with a
sign. The above definition says that, to get the sign of a rigid disk
we compare the capping orientation of the kernel/cokernel of the
$\bar\pa_\kappa$-operator with the orientation on $\aut_{m+1}$ or $\conf_{m+1}$
depending on the number of punctures $m.$
\end{rmk}

\section{Legendrian contact homology over $\Z$}\label{sec:LegCH}
In this section we associate to any Legendrian submanifold
$L\subset\C^n\times\R$ which is equipped with a spin structure a
graded algebra $\A(L)$ over the group ring $\Z[H_1(L)]$ if $L$ is
connected and over $\Z$ otherwise. We define a map
$\pa\colon\A(L)\to\A(L)$ and prove that it is a differential. With
this established we prove that the stable tame isomorphism class of
the differential graded algebra $(\A(L),\pa)$ remains invariant
under Legendrian isotopies. This implies in particular that the
contact homology $\krn(\pa)/\img(\pa)$ is a Legendrian isotopy
invariant. We then show how the differential $\pa$ depends on the
particular spin structure on $L$ and in the final subsection discuss
the relation of our approach to contact homology over $\Z[H_1(S^1)]=\Z[t,t^{-1}]$ for
Legendrian $1$-knots with the completely combinatorial approach
taken in \cite{ENS}.
\subsection{The algebra and its differential}\label{defofal}
Let $L\subset\C^n\times\R$ be (an admissible chord generic) oriented
connected Legendrian submanifold.  Define $\A(L)$ as the free
associative algebra over $\Z[H_1(L)]$ generated by the Reeb chords
of $L$. That is
$$
\A(L)=\Z[H_1(L)]\la c_1,\dots,c_m\ra.
$$
Recall from \cite{ees1} (also see Section~\ref{sec:basnot} above) that each generator $c$ comes equipped with
a grading $|c|$ and a capping path $\gamma_i.$ Elements $A\in H_1(L)$ have gradings
$|A|$. So, $\A(L)$ is a $\Z$-graded algebra. We note that if $L$ is an
orientable manifold then $|A|$ is even for any $A\in H_1(L)$.

If the Legendrian submanifold $L$ is not connected (see also
\cite{Mishachev, Ng03}) then we will use a simpler version of the
theory: we let
$$
\A(L)=\Z\la c_1,\dots, c_m\ra.
$$
In this case the algebra $\A(L)$ has a natural $\Z_2$-grading. There is also a
relative $\Z/m(L)\Z$-grading which we will not discuss here. (Note
that the orientability implies that the Maslov number is even, see
also Remark \ref{rmkdcgr}). In the connected case this simpler version
corresponds to setting all $A\in H_1(L)$ to $1$ and reducing the
grading modulo $2$.
Assume that $L$ is equipped with a spin structure. Define the {\em
differential} $$\pa\colon\A(L)\to\A(L)$$ by requiring that it is
linear over $\Z[H_1(L)]$ (over $\Z$ in the disconnected case), that
it satisfies the graded Leibniz rule on products of monomials
\begin{equation}\label{eqleibniz}
\pa(\alpha\beta)=(\pa\alpha)\beta+(-1)^{|\alpha|}\alpha(\pa\beta),
\end{equation}
and define it on generators as
$$
\pa a =\sum_{\dim\M_A(a;{\mathbf b})=0}
(-1)^{(n-1)(|a|+1)}\left|\M_A(a;{\mathbf b})\right|\,A{\mathbf b},
$$
where $\left|\M_A(a;{\mathbf b})\right|$ is the algebraic number of
rigid disks in the moduli space (where the sign of rigid disk is
defined as in the previous section). It follows from \cite{ees2}
Lemma~1.5, that $\pa$ decreases grading by $1$ in the connected
case, see Remark \ref{rmkdclower1} for the disconnected case.
The purpose of the next subsection is to complete the proof of the
following theorem.
\begin{thm}\label{thmd^2=0}
The map $\pa\colon\A(L)\to\A(L)$ is a differential. That is,
\begin{equation}\label{eqd^2=0}
\pa\circ \pa=0.
\end{equation}
\end{thm}
\begin{pf}
It is a straightforward consequence of the signed Leibniz rule that
the lemma follows once \eqref{eqd^2=0} has been established for generators.
The fact that it holds for generators will be established below.
\end{pf}
\begin{rmk}\label{rmkcoeff}
Note that the above definition also make sense for other coefficient
rings. For example one could replace $\Z$ above with $\Z/k\Z$, for
any $k\in\Z$ or by $\Q$.
\end{rmk}
\subsection{Orientations of $1$-dimensional moduli spaces and their boundaries}\label{subsec:orion1d}
In order to prove Theorem \ref{thmd^2=0}, we will determine the
relations between orientations on $1$-dimensional moduli spaces and
the signs of the pairs of rigid disks which are their boundaries.

\subsubsection{Even and odd punctures and grading}
We connect our abstract definitions of even and odd to the
geometrical situation under study. Let $c$ be a Reeb chord of an
oriented connected chord generic Legendrian submanifold
$L\subset\C^n\times\R$. Then the two tangent spaces of $L,$ the
projections of which intersect at $\Pi_\C(c)=c^\ast,$ are oriented.
We order these tangent spaces by taking the one with the {\em
largest} $z$-coordinate {\em first}. The holomorphic disks we study
have punctures mapping to Reeb chords. Translating the definition of
even and odd punctures to the present situation, we say that a Reeb
chord $c$ is {\em even} if the even negative boundary condition
$R_{ne}(1,\theta)$ takes the orientation of the upper tangent space
to that of the lower. Otherwise we say that $c$ is {\em odd} and, as
is easy to see, in this case the odd negative boundary condition
$R_{no}(1,\theta)$ takes the orientation of the upper tangent space
to the orientation of the lower.
\begin{lma}\label{eochords}
A Reeb chord $c$ is even (odd) if and only if its grading $|c|$ is
even (odd).
\end{lma}
\begin{pf}
Recall that
$$
|c|=\mu(\gamma\ast\lambda)-1,
$$
where $\gamma$ is a capping path of $c$ and where $\lambda$ is the
$J$-trick path connecting the lower to the upper tangent space at
$c^\ast$. Consider instead the inverse of the path
$\gamma\ast\lambda$. This path is
a negative rotation $\hat\pi-\theta$ where
$\theta$ is the complex angle of the upper and lower tangent
spaces followed by $\gamma$ backwards (i.e. from bottom to
top). The Maslov index of this path is even (odd) if and only if
this path preserves (reverses) the orientation of the upper tangent plane. Now,
closing up with the negative even boundary condition path instead of the inverse
$J$-trick path changes the Maslov index by one. Thus the negative even
boundary condition path preserves orientation if and only if
$|c|$ is even.
\end{pf}
\begin{rmk}\label{rmkdcgr}
In the case of an oriented disconnected Legendrian submanifold we
use this notion of even and odd punctures to define the
$\Z_2$-grading discussed above.
\end{rmk}
\begin{rmk}\label{rmkdclower1}
To see that the differential decreases grading by $1$ in the
disconnected case note that (see Proposition~5.14 in \cite{ees2}) the formal
dimension of a component of a moduli space $\M(a;b_1,\dots,b_k)$
with boundary mapping to the collection of paths
$\lambda=(\lambda_1,\dots,\lambda_{k+1})$ on the Legendrian
submanifold $L\subset\C^n\times\R$ equals
$$
n+\mu(\hat\lambda)+k-2,
$$
where $\mu(\hat\lambda)$ is the Maslov index of the closed path
$\hat\lambda$ obtained by closing up the paths $\lambda$ at the
corners by rotating the Lagrangian subspace of the incoming edge to
the Lagrangian subspace at the outgoing edge in the negative
direction. Note that at a negative corner the incoming Lagrangian
subspace is the upper one. Thus, the negative close up preserves
orientation if and only if the puncture is odd. At the positive
puncture the incoming Lagrangian subspace is the lower one and the
negative close up preserves the orientation if and only if $n$ is even
and the puncture is even or $n$ is odd and the puncture is odd. Since
the Maslov index of a loop of Lagrangian subspaces is even if and only
if it is orientation preserving we find that
$$
n+\mu(\hat\lambda)+k-2\equiv |a|_2+\sum_{j=1}^k|b_j|_2+1  \mod{2},
$$
where $|\cdot|_2$ denotes the modulo $2$ grading. We conclude that the
differential changes grading.
\end{rmk}
\subsubsection{Gluing conformal structures and orientations}
We shall determine the relation between orientations on spaces of
conformal structures on a disk which is close to splitting up into
two disks and the conformal structures on the two disks into which
it splits. To this end we first describe the standard orientation on
$\conf_{m+1}$ in various local coordinates. Consider a disk
$D_{m+1}$ with $m+1>3$ punctures $(p_0,\dots,p_m)$ on the boundary
representing the conformal structure $\kappa$ and where $p_0$ is the
positive puncture. Then the standard orientation of $\conf_{m+1}$ at
$\kappa$ is given by the $(m-2)$-tuple of vectors
$\pa_2,\dots,\pa_m$ where $\pa_j$ denotes the vector tangent to $\pa
D_{m+1}$ at $p_j$ and points in the positive direction along $\pa
D_{m+1}$. On the other hand we can coordinatize a neighborhood of
$\kappa$ in $\conf_{m+1}$ by fixing any three punctures on $\pa
D_{m+1}$ and letting the remaining punctures move. We will however
restrict attention to coordinates in which the positive puncture
remains fixed. For such coordinates the positive orientation is
given by the following lemma, the proof of which is straightforward.
\begin{lma}\label{lmashuffle}
The positive orientation of $\conf_{m+1}$ at $\kappa$ in coordinates
obtained by fixing the ordered points $(p_0,p_r,p_{r+s})$ is given
by the $(m-2)$-tuple of vectors
$$
\pa_1,\dots,\pa_{r-1},-\pa_{r+1},-\pa_{r+2},\dots,-\pa_{r+s-2},-\pa_{r+s-1},\pa_{r+s+1},\dots,\pa_m.
$$
\end{lma}
Secondly, we look at the outward normal to the space of conformal
structures at a disk which is obtained by gluing two disks. Let
$m_j\ge 2$, $j=1,2$. Consider $D_{m_1+1}$ and $D_{m_2+1}$ with punctures
$(q_0,\dots,q_{m_1})$ and $(p_0,\dots,p_{m_2})$, respectively. Let
$D_m$ be the disk obtained from gluing $D_{m_1+1}$ to $D_{m_2+1}$ by
identifying $q_0$ and $p_j$, $1\le j\le m_2$. To see the outward
normal of the conformal structure of the glued disk we use coordinates
on $D_{m_2+1}$ which fixes $p_0$, $p_j$, and one of the punctures next
to $p_j$, and coordinates on $D_{m_1+1}$ which fixes $q_0$, $q_1$, and
$q_{m_1}$. Note that $D_{m}$ has punctures corresponding to all $p_k$
except $p_j$ and all $q_k$ except $q_0$. We use coordinates on $T\conf_{m}$ fixing
$p_0$, $q_1$, and $q_{m_1}$. The outward normal is then represented by
the tangent vector to the circle at the puncture formerly next
to $p_j$ directed away from $p_j$.
We next note that there are natural inclusions
$T\conf_{m_j+1}\to T\conf_m$, $j=1,2$. These induce the decomposition
$T\conf_m=\R\oplus T\conf_{m_1+1}\oplus T\conf_{m_2+1}$, where the
$\R$-direction is spanned by the normal direction discussed above. Let
$o_{m_j+1}$ be the standard orientation on $\conf_{m_j+1}$, $j=1,2$ and
let $o$ be the orientation of the outward normal then
\begin{lma}\label{lmaglueconfstr}
The direct sum map above induces the orientation
$-(-1)^{(m_1-1)j}$. In other words,
$$
o\wedge o_{m_1+1}\wedge o_{m_2+1}=-(-1)^{(m_1-1)j}o_m.
$$
\end{lma}
\begin{pf}
We must consider two cases separately according to whether the fixed
puncture in $\pa D_{m_2+1}$ is $p_{j-1}$ or $p_{j+1}$. Let $\pa_{p_k}$
and $\pa_{q_k}$ denote the
positive tangent vector of the boundary of the glued disk at $p_k$ and
$q_k$, respectively.
Assume first that $p_{j-1}$ is fixed. Then by Lemma \ref{lmashuffle}
\begin{align*}
o_{m}=&
\Bigl(\pa_{p_1}\wedge\dots\wedge\pa_{p_{j-2}}\Bigr)
\wedge\pa_{p_{j-1}}
\wedge\Bigl(-\pa_{q_2}\wedge\dots\wedge-\pa_{q_{m_1-1}}\Bigr)
\wedge\Bigl(\pa_{p_{j+1}}\wedge\dots\wedge\pa_{p_{m_2}}\Bigr)\\
=&-(-1)^{(m_1-1)(j-2)}\,\,
\Bigl(-\pa_{p_{j-1}}\Bigr)
\wedge\Bigl(-\pa_{q_2}\wedge\dots\wedge-\pa_{q_{m_1-1}}\Bigr)\\
&\quad\quad\quad \quad\quad \quad\wedge\Bigl(\pa_{p_1}\wedge\dots\wedge\pa_{p_{j-2}}\Bigr)
\wedge\Bigl(\pa_{p_{j+1}}\wedge\dots\wedge\pa_{p_{m_2}}\Bigr)\\
=&-(-1)^{(m_1-1)j}\,\,o\wedge o_{m_1+1}\wedge o_{m_2+1}.
\end{align*}
The computation for $p_{j+1}$ fixed is similar.
\end{pf}
\subsubsection{Disks with few punctures and marked points}
When proving the gluing theorem, Section~7 in \cite{ees2}, holomorphic
disks with $\le 2$ punctures and holomorphic disks with $\ge 3$
punctures were treated simultaneously by adding marked points to disks
with few punctures. We will take the same approach here and thus we
need to discuss orientations on determinants in the presence of marked
points.
Recall that we introduced marked points on a holomorphic disk
$u\colon D_m\to\C^n$ with boundary on a Legendrian submanifold
$L\subset\C^n\times\R$ by picking a codimension one submanifold $H$
of $L$ such that $u(p)\in H$ and such that $du(\pa_p)$ is transverse
to $H$ for some $p\in\pa D_m$. We then consider the
$\bar\pa$-problem for maps $w\colon D_{m+1}\to\C^n$ in a
neighborhood of $u$, which takes the additional marked point $p$
into $H$ and we use a neighborhood of $u$ in that Sobolev space to
find local coordinates on the moduli space. In our study of
orientations below, we will take $H$ as codimension one spheres in
$L$ around one of the endpoints of some Reeb chord to which some
puncture of $u$ maps. We start by describing the corresponding
situation on the level of the linearized equation.
Let $A\colon\pa D_m\to U(n)$ be a trivialized boundary condition
for the $\bar\pa$-operator. Pick $d$ points $q_1,\dots,q_d$ in
$\pa D_m$ which are not punctures. Let
$W_A=\sblv_{2}[A](D_m,\C^n)$ and let
$V_A=\sblv_{1}[0]$ then $\bar\pa_A\colon W_A\to V_A$. For
each $q_j$, pick a linear form $l_j\colon A(q_j)\R^n\to\R$ and
consider the linear functionals
$$
\alpha_j\colon W_A\to\R,\quad \alpha_j(w)=l_j(w(q_j)).
$$
Let $W'_A=\bigcap_{j=1}^d\krn(\alpha_j)$ and define the operator
$$
\bar\pa_A'\colon W'_A\to V_A
$$
as the restriction of $\bar\pa_A$. Then the index of $\bar\pa'_A$
satisfies
$$
\ix(\bar\pa_A')=\ix\bar\pa_A-d.
$$
Pick $d$ elements $\alpha^1,\dots,\alpha^d\in W$ such that
$\alpha_j(\alpha^k)=\delta_j^k$. This gives a direct sum decomposition
$$
W_A=W'_A\oplus\R^d,
$$
which induces the exact sequence
\begin{equation}
\begin{split}
0\longrightarrow & \, \,  \begin{CD}
\krn(\bar\pa_A) @>{\alpha}>> \krn(\bar\pa_A')\oplus\R^d @>{\beta}>>\
\end{CD}\\
&\begin{CD}
\cokrn(\bar\pa_A') @>{\gamma}>> \cokrn(\bar\pa_A)
\end{CD} \longrightarrow 0,
\end{split}
\end{equation}
where $\alpha$ is the direct sum decomposition followed by the
$L^2$-projection to $\krn(\bar\pa'_A)$ in the first summand, where
$\beta$ is $\bar\pa_A$ followed by projection to $\cokrn(\bar\pa_A')$, and
where $\gamma$ is the natural projection induced form the inclusion
$\img(\bar\pa_A')\subset\img(\bar\pa_A)$.
This sequence induces an isomorphism
$$
\Lambda^{\rm max}\krn(\bar\pa_A)\otimes\Lambda^{\rm max}\cokrn(\bar\pa_A)^\ast=
\Lambda^{\rm max}(\krn(\bar\pa'_A)\oplus\R^d)\otimes\Lambda^{\rm
max}\cokrn(\bar\pa'_A)^\ast.
$$
To facilitate our sign discussions we will assume that $d$ is
{\em even}. In that case the position of $\R^d$ in the direct sum
decomposition is of no importance for orientations and we have a
canonical isomorphism
$$
\Lambda^{\rm max}(\krn(\bar\pa'_A)\oplus\R^d)=
\Lambda^{\rm max}\krn(\bar\pa'_A)\otimes\Lambda^{\rm max}(\R^d).
$$
Thus, for $d$ even, an orientation on $\R^d$ gives a canonical
isomorphism between orientations on
$\det(\bar\pa_A')$ and orientations on $\det(\bar\pa_A)$.
We next consider gluing isomorphisms in the presence of marked points.
Let $A\colon D_m\to U(n)$ and $B\colon D_s\to U(n)$ be trivialized
boundary conditions such that the positive puncture of $A$ can be
glued to some negative puncture of $B$.
Assume that there are $2a$ marked points near the positive puncture on
$D_m$ and $2b$ marked points near the negative puncture of $B$ to which
$A$ is glued. Consider the glued boundary condition
$A\sharp B\colon D_{m+s-2}\to U(n)$
and note that $D_{m+s-2}$ inherits the marked points from $D_m$ and
$D_s$, and thus comes equipped with $2(a+b)$ marked points.
Consider the gluing sequence in Lemma \ref{lmaglue2}. Note that the
cut-off of a function in $\krn(\bar\pa_A)$ ($\krn(\bar\pa_B)$)
vanishes on the part of $D_{m+s-2}$ which correspond to $D_s$ ($D_m$).
This observation implies that the gluing sequence obtained by
replacing the operators $\bar\pa_A$, $\bar\pa_B$, and
$\bar\pa_{A\sharp B}$ in Lemma \ref{lmaglue2} with $\bar\pa_A'$,
$\bar\pa_B'$, and $\bar\pa_{A\sharp B}'$, respectively, is also exact.
Let $o_A$ and $o_B$ be orientations on
$\det(\bar\pa_A)$ and $\det(\bar\pa_B)$. Orient $\R^{2a}$ and
$\R^{2b}$. Let $o_A'$ and $o_B'$ be the induced orientations on
$\det(\bar\pa_A')$ and $\det(\bar\pa_B')$ and endow
$\R^{2a}\oplus\R^{2b}$ with the orientation induced from the
orientations of its summands. Let $o_{A\sharp B}$ be the orientation
induced on $\det(\bar\pa_{A\sharp B})$ from the gluing sequence of $\bar\pa_A$
and $\bar\pa_B$. Then $o_{A\sharp B}$ and the orientation on
$\R^{2a}\oplus\R^{2b}$ induces an orientation $o_{A\sharp B}'$ on
$\det(\bar\pa_{A\sharp B}')$.
\begin{lma}
The orientation on $\det(\bar\pa_{A\sharp B}')$ induced from the
gluing sequence of $\bar\pa_A'$ and $\bar\pa_B'$ and the orientations
$o_A'$ and $o_B'$ equals $o'_{A\sharp B}$.
\end{lma}
\begin{pf}
To show this we consider the following diagram.
$$
\begin{CD}
\left[\begin{matrix}
\krn\bar\pa_A\\
\krn\bar\pa_B
\end{matrix}\right]
@>>>
\left[\begin{matrix}
\krn\bar\pa_A'\\
{\mathbb R}^{2a}\\
\krn\bar\pa_B'\\
{\mathbb R}^{2b}
\end{matrix}\right]
@>>>
\left[\begin{matrix}
\cokrn\bar\pa_A'\\
\cokrn\bar\pa_B'
\end{matrix}\right]
@>>>
\left[\begin{matrix}
\cokrn\bar\pa_A\\
\cokrn\bar\pa_B
\end{matrix}\right]\\
@AAA @AAA @VVV @VVV \\
\krn\bar\pa_{A\sharp B}
@>>>
\left[\begin{matrix}
\krn\bar\pa_{A\sharp B}'\\
{\mathbb R}^{2a}\\
{\mathbb R}^{2b}
\end{matrix}\right]
@>>> \cokrn\bar\pa_{A\sharp B}'
@>>> \cokrn\bar\pa_{A\sharp B}\\
@AAA @AAA @VVV @VVV
\end{CD}
$$
where the $0$'s on the left and right in the first and second
horizontal rows have been dropped and where the $0$'s below the lower
row of arrows have
been dropped as well. The upper horizontal row is
the direct sum of the sequences inducing $o'_A$ and $o'_B$ from $o_A$
and $o_B$, respectively. The lower horizontal row is the sequence
inducing $o'_{A\sharp B}$ from $o_{A\sharp B}$.
The gluing sequence for $\bar\pa_{A\sharp B}$ is obtained from the
above diagram by adding an arrow from the top left entry to the top
right entry. (Thinking of it as a half-circular arc the gluing sequence
is then the outer boundary of the diagram.)
The gluing sequence for $\bar\pa_{A\sharp B}'$ is obtained from the
above diagram by adding an arrow from the middle left entry in top row
to the middle right entry in the same row and forgetting the $\R^{2a}$
and $\R^{2b}$ summands. (Thinking also of this arrow as a
half-circular arc the gluing sequence is the inner boundary of the
diagram.)
Let $(\sigma,\bar\sigma)$ be wedges of vectors representing the
orientation pair on the pair
$$
\left(\krn(\bar\pa_{A\sharp B}),
\cokrn(\bar\pa_{A\sharp B})\right).
$$
The induced orientation on the pair
$$
\left(
\left[\begin{matrix}
\krn(\bar\pa_{A})\\
\krn(\bar\pa_{B})
\end{matrix}\right],
\left[\begin{matrix}
\cokrn(\bar\pa_{A})\\
\cokrn(\bar\pa_{B})
\end{matrix}\right]
\right)
$$
is then represented by
$$
(\sigma\wedge\omega,\bar\omega\wedge\bar\sigma),
$$
where $\omega$ is a wedge of vectors on the complement of the image of
$\sigma$. This pair in turn induces the orientation pair
$$
(\sigma\wedge\omega\wedge\phi,\bar\phi\wedge\bar\omega\wedge\bar\sigma),
$$
on the pair
\begin{equation}\label{eqmidpair}
\left(
\left[\begin{matrix}
\krn(\bar\pa_{A}')\\
\R^{2a}\\
\krn(\bar\pa_{B}')\\
\R^{2b}
\end{matrix}\right],
\left[\begin{matrix}
\cokrn(\bar\pa_{A}')\\
\cokrn(\bar\pa_{B}')
\end{matrix}\right]
\right)
\end{equation}
where $\phi$ is a wedge of vectors on the complement of the image of
$\sigma\wedge\omega$.
On the other hand the orientation pair $(\sigma,\bar\sigma)$ induces
the orientation pair
$$
(\sigma\wedge\theta,\bar\theta\wedge\bar\sigma)
$$
on the pair
$$
\left(
\left[\begin{matrix}
\krn(\bar\pa_{A\sharp B}')\\
\R^{2a}\\
\R^{2b}
\end{matrix}\right],
\cokrn(\bar\pa_{A\sharp B}')\right).
$$
This in turn induces an orientation
$$
(\sigma\wedge\theta\wedge\tau,\bar\tau\wedge\bar\theta\wedge\bar\sigma)
$$
on the pair in \eqref{eqmidpair}. Noting that the orientations of
$\omega\wedge\phi$ and $\theta\wedge\tau$ agree if and only if those
of $\bar\omega\wedge\bar\phi$ and $\bar\theta\wedge\bar\tau$ do and
that the vertical maps on the $\R^{2a}\oplus\R^{2b}$ summand is the
identity we find that the lemma holds.
\end{pf}

We now turn to induced orientations on moduli spaces. We show that
there is a natural way to pick marked points and an orientation on the
corresponding copies of $\R$ so that the orientation of a moduli space
of a disk with marked points is determined exactly as for disks with
many punctures. In our applications we use only moduli spaces of
dimensions $0$ and $1$ so we consider these two cases. We also
assume we are in a generic situation where all moduli spaces are
transversely cut out.
Let $u\colon D_m\to\C^n$ be a holomorphic disk in a moduli space of
dimension $0$ or $1$. Note then that the corresponding
$\bar\pa_\kappa$-operator which is a part of the linearization of the
$\bar\pa$-map at $u$ ($\kappa$ denotes the conformal structure) has
either only kernel or only cokernel. Consider
the case of only cokernel first. In that case we choose an even number
of marked points $q_1,\dots,q_d$ near a puncture of $u$. We obtain in
that way a new source disk $\tilde D_{m,d}$ with $m$ punctures and $d$
marked points. Pick local $1$-parameter families of diffeomorphisms
$\phi_j^\tau$ around
the marked points which move them
in the positive direction along the boundary and set
$\alpha_j(w)=\la du(\pa_\tau\phi_j(q_j)),w(q_j)\ra$.
In this case the orientation inducing sequence reduces to
$$
\begin{CD}
0 @>>> \R^d @>>> \cokrn(\bar\pa_\kappa') @>>> \cokrn(\bar\pa_\kappa)
@>>> 0.
\end{CD}
$$
Viewing  $\cokrn(\bar\pa_\kappa)$ as a subspace of
$\sblv_{1,\epsilon}[0]$, we obtain a splitting
$\cokrn(\bar\pa_\kappa)\to\cokrn(\bar\pa_\kappa')$ and a
corresponding direct sum decomposition
$\cokrn(\bar\pa_\kappa')=\cokrn(\bar\pa_\kappa)\oplus\R^d$.
We obtain the diagram
$$
\begin{CD}
T\conf_{m}\oplus\R^d @>>> T\conf_{m+d}\\
@VVV @VVV\\
\cokrn\bar\pa_\kappa\oplus\R^d @>>>\cokrn(\bar\pa_\kappa'),
\end{CD}
$$
where the top arrow is an isomorphism with inverse which maps the
tangent vector of a moving marked point to the corresponding
factor in $\R^d$. Note that the marked points in $\pa D_m$ lies in an
arc which contains no punctures. Thus, this map is orientation
preserving since $d$ is even. Likewise the leftmost vertical map is
the identity on the $\R^d$-component. It follows that the orientation
on the moduli space induced by $T\conf_m\to\cokrn(\bar\pa_\kappa)$ is
the same as the one induced by
$T\conf_{m+d}\to\cokrn(\bar\pa_\kappa')$.
Consider next the case when $\cokrn(\bar\pa)=0$. In this case, if $d$
is sufficiently large then $\krn(\bar\pa')=0$ and the orientation
inducing sequence reduces to
$$
\begin{CD}
0 @>>> \krn(\bar\pa) @>>> \R^d @>>> \cokrn(\bar\pa') @>>> 0.
\end{CD}
$$
Thus we get $\cokrn(\bar\pa')=\R^d/\krn(\bar\pa)$. We use the
convention that marked points are added near the positive puncture in
the negative direction of it. We then get the diagram
$$
\begin{CD}
\R^d/T\A_m @>>> T\conf_{3-m + d}\\
@VVV  @VVV\\
\R^d/\krn(\bar\pa) @>>> \cokrn(\bar\pa')
\end{CD},
$$
where the top horizontal map is the quotient of the (oriented)
automorphism group acting on the positively oriented tangent vectors
to the marked points moving on the boundary which is then orientation
preserving, and where the left vertical arrow is the map induced by
$\A_m\to\krn(\bar\pa)$. It follows again that the orientation on the
moduli space induced by $T\A_m\to\krn(\bar\pa)$ is the same as that
induced by $T\conf_{3-m + d}\to\cokrn(\bar\pa')$.

\subsubsection{Gluing linearized operators corresponding to rigid disks}
We now determine how capping orientations behave under gluing as in
Lemma \ref{lmaglue2}. In fact, we concentrate on the case important for our
applications: when two {\em rigid} disks $u$ and $v$ with boundary on
a Legendrian submanifold $L\subset\R\times\C^n$ are being glued.
As we have seen gluings of disks with few punctures (when
kernels of $\bar\pa$-operators are present) can be treated as
gluings of disks with many punctures. We concentrate on the case of
many punctures first and explain in Remark \ref{rmkfewpunct} how to modify the
arguments in the presence of marked points.
Let $A$ be a trivialized stabilized boundary condition on the $m$-punctured
disk (corresponding to $u$ above), with positive puncture at a Reeb
chord $b_k$ and negative punctures at Reeb chords $f_1,\dots,f_{m-1}$.
Let $B$ be a trivialized stabilized  boundary condition on the $(r+1)$-punctured
disk (corresponding to $v$ above), with
positive puncture at a Reeb chord $a$ and negative punctures at Reeb
chords $b_1,\dots,b_k,\dots,b_r$. Let
$\bar\pa_{A\sharp B}$ be the operator which is obtained by gluing
$\bar\pa_A$ and $\bar\pa_B$ at $b_k$. Since we model the case when the
disks $u$ and $v$ are rigid and have many punctures
we assume that
$\krn(\bar\pa_A)=0=\krn(\bar\pa_B)$ and that
$$
\ix(\bar\pa_A)=-(m-3),\qquad
\ix(\bar\pa_B)=-(r-2).
$$
In particular, it follows from the dimension formula that
\begin{align}\label{eqdimu}
&|b_k|-\sum_{j=1}^{m-1}|f_j|=1-|A_u|\\\label{eqdimv}
&|a|-\sum_{j=1}^r|b_j|=1-|B_v|,
\end{align}
where $A_u, B_v\in H_1(L)$ are the homology classes of the boundary
paths of $u$ and $v$ capped off with the capping paths at each
puncture. Since $L$ is orientable, $|A_u|$ and $|B_v|$ are even.
Let $o_A$, $o_B$ and $o_{A\sharp B}$ denote the capping orientations
of the determinant lines of $\bar\pa_A$, $\bar\pa_B$, and
$\bar\pa_{A\sharp B}$.
\begin{lma}\label{lmaglueridgop}
Let $\sigma_j=\pm 1$, $j=1,2$.
The gluing map in Lemma \ref{lmaglue2} induces from orientations
$\sigma_1o_A$ and $\sigma_2o_B$ the orientation
$$
(-1)^{n+1}(-1)^{mk}(-1)^{\sum_{j=1}^{k-1}|b_j|}(-1)^{(|b_k|+1)(n-1)}
(\sigma_1\sigma_2)o_{A\sharp B}.
$$
\end{lma}
\begin{pf}
Consider the boundary condition $\hat A\sharp \hat B$ obtained by
gluing the two boundary conditions $\hat A$ and $\hat B$.
This boundary condition can be obtained in two ways.
\begin{itemize}
\item[(i)]
Glue all capping disks to $A$ and $B$, respectively and
then join $\hat A$ and $\hat B$.
\item[(ii)]
Glue the operators $\bar\pa_{b_k+}$ and $\bar\pa_{b_k-}$ to obtain an
operator $\bar\pa_E$ on the closed disk. Glue all relevant capping
disks to $A\sharp B$ giving the boundary condition $\widehat{A\sharp
B}$ on the closed disk.  Then glue $\bar\pa_E$ to
$\bar\pa_{\widehat{A\sharp B}}$.
\end{itemize}
Repeated application of Lemmas \ref{lmaassocpp} and \ref{lmaassocbp}
gives the following gluing sequence corresponding to the gluing (i):
\begin{equation*}
\begin{split}
0\longrightarrow & \, \,  \begin{CD}
\krn(\bar\pa_{\hat A\sharp \hat B}) @>>>
\left[
\begin{matrix}
\bigoplus_{j=1}^{m-1}\krn(\bar\pa_{f_j-}) \\
\bigoplus_{j=1}^{r}\krn(\bar\pa_{b_j-})
\end{matrix}
\right]
@>>>
\end{CD}\\
&\begin{CD}
\left[
\begin{matrix}
\bigoplus_{j=1}^{m-1}\cokrn(\bar\pa_{f_j-})\\
\cokrn(\bar\pa_{b_k+})\\
\cokrn(\bar\pa_A)\\
\R^n\\
\bigoplus_{j=1}^{r}\cokrn(\bar\pa_{b_j-})\\
\cokrn(\bar\pa_{a+})\\
\cokrn(\bar\pa_B)
\end{matrix}
\right]
@>>>
\cokrn(\bar\pa_{\hat A\sharp \hat B})
\end{CD} \longrightarrow 0.
\end{split}
\end{equation*}
Note that the capping orientation $o_A$ together with the orientations
on the capping disks of $A$ induce the canonical orientation on
$\det(\bar\pa_{\hat A})$, by definition. Similarly,
the capping orientation $o_B$ together with the orientations
on the capping disks of $B$ induce the canonical orientation on
$\det(\bar\pa_{\hat B})$. Thus Lemma \ref{lmagluecanor} implies that
that these orientations in the above sequence induce the canonical
orientation on $\bar\pa_{\hat A\sharp \hat B}$.

To find the corresponding gluing sequence for (ii) we first look at
the gluing sequence for $\bar\pa_{A\sharp B}$. This sequence is the
following.
$$
\begin{CD}
0 @>>> \cokrn(\bar\pa_A)\oplus\cokrn(\bar\pa_B) @>>>
\cokrn(\bar\pa_{A\sharp B})@>>> 0,
\end{CD}
$$
It allows us to identify the two spaces involved.
Second the gluing sequence for $\bar\pa_E$ is
$$
\begin{CD}
0 @>>> \krn(\bar\pa_E) @>>>
\krn(\bar\pa_{b_k-}) @>>>\\
\left[\begin{matrix}
\cokrn(\bar\pa_{b_k+})\\
\cokrn(\bar\pa_{b_k-})
\end{matrix}\right]
@>>> \cokrn(\bar\pa_E) @>>> 0,
\end{CD}
$$
where both the leftmost and the rightmost non-trivial maps are
isomorphisms. Moreover, by definition, the chosen orientation pair on
$$
(\krn(\bar\pa_{b_k-}),\cokrn(\bar\pa_{b_k+})\oplus\cokrn(\bar\pa_{b_k-}))
$$
corresponds to the canonical orientation on $\bar\pa_E$.
With these two identifications made we apply Lemmas \ref{lmaassocpp}
and \ref{lmaassocbp} and find that the gluing sequence in case (ii) is
\begin{equation*}
\begin{split}
0\longrightarrow & \, \,  \begin{CD}
\krn(\bar\pa_{\hat A\sharp\hat B})
@>>>
\left[
\begin{matrix}
\krn(\bar\pa b_k-)\\
\bigoplus_{j=1}^{k-1}\krn(\bar\pa_{b_j-})\\
\bigoplus_{j=1}^{m-1}\krn(\bar\pa_{f_j-})\\
\bigoplus_{j=k+1}^{r}\krn(\bar\pa_{b_j-})
\end{matrix}
\right]
@>>>
\end{CD}\\
&\begin{CD}
\left[
\begin{matrix}
\cokrn(\bar\pa_{b_k+})\\
\cokrn(\bar\pa_{b_k-})\\
\R^n\\
\bigoplus_{j=1}^{k-1}\cokrn(\bar\pa_{b_j-})\\
\bigoplus_{j=1}^{m-1}\cokrn(\bar\pa_{f_j-})\\
\bigoplus_{j=k+1}^{r}\cokrn(\bar\pa_{b_j-})\\
\cokrn(\bar\pa_{a+})\\
\cokrn(\bar\pa_A)\\
\cokrn(\bar\pa_B)
\end{matrix}
\right]
@>>>
\krn(\bar\pa_{\hat A\sharp\hat B})
\end{CD} \longrightarrow 0.
\end{split}
\end{equation*}
Lemma \ref{lmagluecanor} implies that the chosen orientations on all
capping disks
together with the capping orientation $o_{A\sharp B}$ on
$\cokrn(\bar\pa_{A\sharp B})=\cokrn(\bar\pa_A)\oplus\cokrn(\bar\pa_B)$
induce the canonical orientation on $\bar\pa_{\hat A\sharp \hat B}$.
We are now in position to determine the gluing sign.
Note first that in both sequences the middle map takes the
complement of $\krn(\bar\pa_{b_k-})$ to $0$ and maps this vector
non-trivially into the $\R^n$ summand in the next term in the sequence. To compare the signs we first move
$\R^n$ so that it becomes the first summand in its respective sum in both
sequences.
Note that because of \eqref{eqdimu}, the dimension of the space
$$
\left[\begin{matrix}
\bigoplus_{j=1}^{m-1}\cokrn(\bar\pa_{f_j-})\\
\cokrn(\bar\pa_{b_k+})
\end{matrix}\right]
$$
equals $(n-2)+{\rm even}$ if $b_k$ is odd and $(n-1)+{\rm odd}$ if
$b_k$ is even. Thus, the dimension of this space is congruent to
$n$ modulo $2$ in either case. In sequence (i), moving $\R^n$ to the first position
thus gives an orientation change of sign
$$
(-1)^{n(m-3)}(-1)^{n^2}=(-1)^{n(m-2)}.
$$
To find the corresponding change in the sequence (ii) note that the
dimension of the space
$$
\left[\begin{matrix}
\cokrn(\bar\pa_{b_k+})\\
\cokrn(\bar\pa_{b_k-})
\end{matrix}\right]
$$
equals $n-1$ and we find that the orientation change is
$$
(-1)^{n(n-1)}=1.
$$
The orientation difference is thus
$$
(-1)^{n(m-2)}.
$$
With $\R^n$ moved to the first position we compute the sign change by
changing the order of summands in sequence (i) to agree with those in
sequence (ii):
The orientation difference in the two sequences arising from their
second terms are
$$
(-1)^{m-1+k-1}\cdot (-1)^{(m-1)(k-1)}=(-1)^{mk-1},
$$
where the first factor comes from moving $\krn(\bar\pa_{b_k-})$ to the
left and the second from moving $\oplus_j\krn(\bar\pa_{f_j-})$ to its
position.
The orientation difference arising in the third term in the sequences
(with $\R^n$ already moved to the first position) is calculated as follows.
First we move the $\cokrn(\bar\pa_{b_k+})$-term in sequence (i) to the
left. This subspace has dimension $n-1$ if $|b_k|$ is even and $n-2$ if $|b_k|$
is odd. Since $|b_k|-\sum_j|f_j|=1$ modulo 2, by \eqref{eqdimu}, we
see that if $|b_k|$ is odd then
there is no sign arising from this move since
$\dim(\bigoplus_j\cokrn(\bar\pa_{f_j}))$ is even. On the other hand if
$|b_k|$ is even we find the sign $(-1)^{(n-1)}$. In general the sign
is thus
$$
(-1)^{(|b_k|+1)(n-1)}.
$$
Second we move the subspace
$\bigoplus_j\cokrn(\bar\pa_{b_j-})\oplus\cokrn(\bar\pa_{a+})$ over
$\cokrn(\bar\pa_A)$. The latter space has dimension $m-3$. We must
calculate the dimension modulo $2$ of the former. Since
$|a|-\sum_j|b_j|=1$ modulo $2$, we find that the dimension of this
space is $(n-1+{\rm odd})$ if $a$ is even and $(n-2+{\rm even})$ if
$a$ is odd. Thus in any case we get a sign
$$
(-1)^{(m-3)(n-2)}
$$
from this motion.
Finally we must move out $\cokrn(\bar\pa_{b_k-})$ if this term is
non-zero and move $\bigoplus_j\cokrn(\bar\pa_{f_j-})$ to the right
position. We consider two cases. First if $b_k$ is even then
$\cokrn(\bar\pa_{b_k-})$ is zero dimensional and
$\bigoplus_j\cokrn(\bar\pa_{f_j-})$ is odd dimensional giving a sign
$(-1)^{\sum_{j=1}^{k-1}|b_j|}$. On the other hand if $b_k$ is odd
then $\bigoplus_j\cokrn(\bar\pa_{f_j-})$ is even dimensional and the
sign is still $(-1)^{\sum_{j=1}^{k-1}|b_j|}$ since we must move the
$1$-dimensional space $\cokrn(\bar\pa_{b_k-})$. Thus in any case we
get
$$
(-1)^{\sum_{j=1}^{k-1}|b_j|}.
$$
Multiplying out all the signs gives the over all sign
\begin{align*}
&-(-1)^{(|b_k|+1)(n-1)}(-1)^{mk}(-1)^{(\sum_{j=1}^{k-1}|b_j|)}(-1)^{(n-2)(m-3)}(-1)^{n(m-2)}\\
&= (-1)^{n+1}(-1)^{mk}(-1)^{\sum_{j=1}^{k-1}|b_j|}(-1)^{(|b_k|+1)(n-1)},
\end{align*}
as claimed.
\end{pf}

\begin{rmk}\label{rmkfewpunct}
In the presence of marked points we simply need to replace the
operators $\bar\pa_A$, $\bar\pa_B$, and $\bar\pa_{A\sharp B}$ above
with $\bar\pa_A'$, $\bar\pa_B'$, and $\bar\pa_{A\sharp B}'$,
respectively. In particular by adding sufficiently large even numbers
of  marked points we can assure that none of the latter operators has
non-trivial
kernel. Of course adding marked points effects the dimensions of the
cokernels of the operators but the final expression for the sign is
independent of these dimensions so the lemma holds in the presence of
marked points as well.
\end{rmk}

\subsubsection{Orientations of $1$-dimensional moduli spaces}
Let $L\subset\C^n\times\R$ be a generic Legendrian submanifold
oriented and with a fixed spin structure. Then as shown above, all
moduli spaces of holomorphic disks with boundary on $L$ come equipped
with orientations. We show the following result.
Let $\M$ be a component of a $1$-dimensional moduli space with
boundary $\pa \M= \M_0\cup\M_1$. For $\M_j$, $j=0,1,$
let $\sigma(\M_j)=1$ if the orientation of $\M$ is outwards at $\M_j$
and let $\sigma(\M_j)=-1$ otherwise.
A boundary component $\M_j$ is a broken holomorphic disk. That is,
two rigid holomorphic disks $u$ and $v$ such that the
positive puncture of $u$ is identified with some negative puncture of
$v$. Assume that the positive puncture of $u$ is $b_k$ and that its
negative punctures are $f_1,\dots,f_{m-1}$. Assume also that the
positive puncture of $v$ is $a$ and that its negative punctures are
$b_1,\dots,b_k,\dots,b_r$.
\begin{lma}\label{lmaMorient}
Let $\mu_1$ and $\mu_2$ be the signs of $u$ and $v$, respectively.
Then
$$
\sigma(\M_j)=\mu_1\mu_2(-1)^n(-1)^{\sum_{j=1}^{k-1}|b_j|}(-1)^{(n-1)(|b_k|+1)}.
$$
\end{lma}
\begin{pf}
It is enough to consider the case of many punctures, see Remark
\ref{rmkfewpunct}, so assume $m\ge 3$ and $r\ge 2$. Let $A$
and $B$ be the
boundary conditions corresponding to $u$ and  $v,$
respectively. Consider the commutative diagram
$$
\begin{CD}
\R\oplus T\conf_{m}\oplus T\conf_{r+1} @>>> T\conf_{m+r-1}\\
@VVV @ VVV\\
\cokrn(\bar\pa_{A})\oplus\cokrn(\bar\pa_{B}) @>>>
\cokrn(\bar\pa_{A\sharp B}).
\end{CD}
$$
Here $T\conf_{m}$, $T\conf_{r+1}$, and $T\conf_{m+r-1}$ are endowed
with their standard orientations and $\R$ is endowed with the
orientation corresponding to the {\em outward} normal $\nu$. Moreover, the
spaces appearing in the lower horizontal row are endowed with their
capping orientations. We compute the orientations.
The leftmost vertical map restricted to the complement of $\R$ is an
isomorphism between oriented spaces of sign $\mu_1\mu_2$. The lower
horizontal arrow is an isomorphism of oriented vector spaces with sign
$$
\sigma=(-1)^{n+1}(-1)^{mk}(-1)^{\sum_{j=1}^{k-1}|b_j|}(-1)^{(n-1)(|b_k|+1)}
$$
by Lemma \ref{lmaglueridgop}. The upper horizontal arrow is an
isomorphism of sign
$\tau=-(-1)^{mk}$ by Lemma \ref{lmaglueconfstr}.
Now, the oriented tangent space to $\M$ is identified with the
oriented kernel of the rightmost vertical map. Chasing the diagram we
find that this oriented kernel is identified with $\R$ oriented by
\begin{align*}
\tau\sigma\mu_1\mu_2\cdot \nu=
(-1)^n(-1)^{\sum_{j=1}^{k-1}|b_j|}(-1)^{(n-1)(|b_k|+1)}\mu_1\mu_2\cdot \nu,
\end{align*}
as claimed.
\end{pf}

Finally we are in position to complete the proof of Theorem \ref{thmd^2=0}.
\begin{pf}[Completion of the proof of Theorem \ref{thmd^2=0}]
Let $a$ be a Reeb chord of $L$ and let $C b_1,\dots,b_r$ be a word
appearing in $\pa\pa a$. Consider a component of a
$1$-dimensional moduli space
$$
\M=\M_C(a;b_1,\dots,b_r),
$$
with oriented boundary $\pa \M=\M_1-\M_0$. Assume that $\M_0$
consists of two broken disks in
\begin{equation}\label{eqdisks1}
\M_A(a;b_1,\dots,b_{k-1},f,b_{k+s},\dots,b_r)\text{ and }
\M_B(f;b_k,\dots,b_{k+s-1}),
\end{equation}
respectively, where $A+B=C$,
and that $\M_1$ consists of two broken disks in
\begin{equation}\label{eqdisks2}
\M_{A'}(a;b_1,\dots,b_{k'-1},f',b_{k'+s'},\dots,b_r)\text{ and }
\M_{B'}(f';b_{k'},\dots,b_{k'+s'-1}),
\end{equation}
respectively, where $A'+B'=C$.
By the definition of the differential the disks in \eqref{eqdisks1}
contribute
\begin{align*}
&(-1)^{(n-1)(|a|+1)}(-1)^{|A|}
(-1)^{\sum_{j=1}^{k-1}|b_j|}(-1)^{(n-1)(|f|+1)}\mu_1\mu_2
(A+B)b_1\dots b_r=\\
&(-1)^{(n-1)(|a|+1)}(-1)^{\sum_{j=1}^{k-1}|b_j|}(-1)^{(n-1)(|f|+1)}\mu_1\mu_2
Cb_1\dots b_r,
\end{align*}
where $\mu_1$ and $\mu_2$ are the signs of the two disks and where we
use the fact that $|A|$ is even.
On the other hand the disks in \eqref{eqdisks2}
contribute
\begin{align*}
&(-1)^{(n-1)(|a|+1)}(-1)^{|A'|}
(-1)^{\sum_{j=1}^{k'-1}|b'_j|}(-1)^{(n-1)(|f'|+1)}\mu_1'\mu_2'
(A'+B')b_1\dots b_r=\\
&(-1)^{(n-1)(|a|+1)}(-1)^{\sum_{j=1}^{k'-1}|b_j|}(-1)^{(n-1)(|f'|+1)}\mu_1'\mu_2'
Cb_1\dots b_r,
\end{align*}
where $\mu_1'$ and $\mu_2'$ are the signs of the two disks and where we
use the fact that $|A'|$ is even.
By Lemma \ref{lmaMorient}
\begin{align*}
-(-1)^n=(-1)^n\sigma(\M_0)=(-1)^{\sum_{j=1}^{k-1}|b_j|}(-1)^{(n-1)(|f|+1)}\mu_1\mu_2,\\
(-1)^n=(-1)^n\sigma(\M_1)=(-1)^{\sum_{j=1}^{k'-1}|b_j|}(-1)^{(n-1)(|f'|+1)}\mu'_1\mu'_2.
\end{align*}
Thus these two terms cancel in $\pa\pa a$. Since all terms
contributing to $\pa\pa a$ arise in this way we conclude $\pa\pa
a=0$.
\end{pf}

\subsection{Stable tame isomorphism invariance of $\A(L)$}\label{sec:invofal}
We show that Legendrian isotopies of a Legendrian submanifold
$L\subset\C^n\times\R$ with a fixed spin structure does not change the
stable tame isomorphism class of its associated algebra $\A(L)$. More
precisely
\begin{thm}\label{thminv}
Let $L_t\subset\C^n\times\R$, $0\le t\le 1$ be a Legendrian isotopy
then the DGA $(\A(L_0),\pa_0)$ of $L_0$ and the DGA $(\A(L_1),\pa_1)$
of $L_1$ are stable tame isomorphic. (In particular the contact
homology of $L_0$ is isomorphic to that of $L_1$.)
\end{thm}
\begin{pf}
The theorem follows from Lemma \ref{lmaeasyinv} and Corollaries
\ref{corhsinv} and Lemma~\ref{ststi} which are proved below.
\end{pf}
Our approach is to study the bifurcations in moduli spaces under
generic isotopies. During such deformations three events effect the
moduli-spaces (and therefore possibly the differential): they
undergo Morse modifications, there appear disks of formal dimension
$-1$ (so called handle slide disks because of analogous phenomena in
Morse theory), and there appear self tangency instances. To prove
invariance it is sufficient to study these events separately. To
show invariance in the first case we use a parameterized moduli
space. To show invariance in the second we use an auxiliary
Legendrian embedding of $L\times [0,1]$ with certain boundary
conditions at $L\times\{0\}$ and $L\times\{1\}$, and its
differential graded algebra. We mention these proofs have an
advantage over the more straight forward proofs of invariance given
in \cite{ees2}. In particular, we avoid the need for ``degenerate
gluing'', which is technically much more difficult than the gluing
results needed to prove $\partial^2=0.$
\subsubsection{Invariance under trivial isotopies}
Let $L_t\subset\C^n\times\R$, be a generic $1$-parameter family of
Legendrian submanifolds such that $L_0$ and $L_1$ are generic and such
that there are neither handle slide moments nor self tangency moments
during the isotopy.
that there are neither handle slide moments nor self tangency moments
during the isotopy.
\begin{lma}\label{lmaeasyinv}
$$
(\A(L_0),\pa_0)=(\A(L_1),\pa_1).
$$
\end{lma}
\begin{pf}
Note that since the complex angles in the added directions are taken
to be smaller than any complex angle of a double point of $L_t$, the
capping orientations are continuous in $t$. Therefore, the
parameterized moduli spaces are naturally oriented. In particular,
the rigid disks on $L_0$ and $L_1$ equals the oriented boundary of
the one-dimensional parameterized moduli space. This oriented
cobordism immediately gives $\pa_1=\pa_0$.
\end{pf}

\subsubsection{An auxiliary Legendrian submanifold}
We associate to a $1$-parameter family of Legendrian embeddings
$\phi_t\colon L\to\C^n\times\R$, $0\le t\le 1$,
Legendrian embeddings
$\Phi_f^\delta\colon L\times\R\to\C^{n+1}\times\R$, depending on
$\delta>0$ and a positive Morse function $f\colon\R\to\R$.
Let $\phi_t\colon L\subset\C^n\times\R$, $t\in[-1,1]$ be a Legendrian
isotopy. For small $\delta>0$, fix smooth non-decreasing functions
\begin{equation}\label{eqalphadelta}
\alpha^\delta\colon[-1,1]\to[-\delta,\delta]
\end{equation}
such that $\alpha^\delta(\pm t)=\pm\delta$ for $\frac 34\le t\le 1$,
and such that $\alpha^\delta(t)=\delta t$ for $-\frac14\le
t\le\frac14$. Note that $\alpha^\delta\to 0$ as $\delta\to 0$.
Fix standard coordinates
$$
\bigl((x_1,y_1,\dots,x_n,y_n),z\bigr)=(x,y,z)
$$
on $\C^n\times\R$. Define $\phi_t^\delta$, $t\in\R$ as
$$
\phi_t^\delta=\begin{cases}
\phi_{-\delta}   &\text{for $t\in(-\infty,-1]$,}\\
\phi_{\alpha^\delta(t)} &\text{for $t\in[-1,1]$,}\\
\phi_{\delta}   &\text{for $t\in[1,\infty)$.}
\end{cases}
$$
Write
$$
\phi_t^\delta(q)=\bigl(x_t(q),y_t(q),z_t(q)\bigr),\quad q\in L.
$$
Fix a {\em positive} Morse function $f\colon\R\to\R$ and
$\delta>0$. Let $f'(t)=\frac{df}{dt}$ denote the derivative of $f$.
Define $\Phi_f^\delta\colon \R\times L\to\C\times\C^n\times\R$,
$$
\Phi_f^\delta(t,q)=
\bigl(x_0(t,q),y_0(t,q),x(t,q),y(t,q),z(t,q)\bigr),\quad
(t,q)\in\R\times L,
$$
where
\begin{align*}
x_0(q,t) &= t,\\
y_0(q,t) &= f(t)\left(\frac{\pa z_{t}}{\pa t}
-{y_j}_{t}(q){\frac{\pa {x_j}_{t}}{\pa t}}\right)
+f'(t)z_{t}(q),\\
x(q,t) &= x_{t}(q),\\
y(q,t) &= f(t)y_{t}(q),\\
z(q,t) &= f(t)z_{t}(q).
\end{align*}
It is straightforward to check that $\Phi_f^\delta$ is a Legendrian
embedding.
Assume that the Morse function $f\colon\R\to\R$ above has local minima
at $\pm 1$ and no critical points in the region
$(-\infty,-1)\cup(1,\infty)$. Then the $x_0$-coordinate $c_0$ of each Reeb
chord $c$ of $\Phi$ satisfies $|c_0|\le 1$.
Let $u$ be a holomorphic disk with boundary on $\Phi_f^\delta$ and one
positive puncture.
\begin{lma}\label{lmamin}
If the positive puncture of $u$ maps to a Reeb chord $c$ of $\Phi$
with $c_0=\pm 1$. Then the image of $u$ lies in $\{x_0=\pm1\}$.
\end{lma}
\begin{pf}
For definiteness assume the positive puncture of $u$ maps to $c$ with
$c_0=1$. Project $\Phi_f^{\delta}$ to the $(x_0,y_0)$-plane. The
image of this projection is contained in the region
$$
\bigl\{-\alpha |x_0-1|\le y_0\le\alpha|x_0-1|\bigr\}
$$
for some $\alpha$. If the projection $u_0$ of $u$ to the
$(x_0,y_0)$-plane is non-constant then it covers at least one of the
regions
$$
\bigl\{-\alpha |x_0-1|> y_0\bigr\}
\cap B_r\bigl((1,0)\bigr)\quad\text{ or }\quad
\bigl\{\alpha|x_0-1|<y_0\bigr\}
\cap B_r\bigl((1,0)\bigr),
$$
for some ball $B_r\bigl((1,0)\bigr)$. Since $u$ has boundary on
$\Phi_k^\delta$, $u_0$ takes no boundary point to the line
$\{x_0=1\}$. This and the above covering property contradicts $u_0$
being bounded in the $y_0$-direction. The lemma follows.
\end{pf}

\begin{lma}\label{lmaalg}
The image of every holomorphic disk with boundary on $\Phi_f^\delta$
is contained in the region $\{|x_0|\le 1\}$.
\end{lma}
\begin{pf}
Arguing as in the proof of Lemma \ref{lmamin} we find that the projection to
the $(x_0,y_0)$-plane of a holomorphic disk with boundary
on $\Phi_f^\delta$ can not intersect the lines $\{x_0=\pm 1\}$ in interior
points. It follows that the image of any disk lies entirely in one of
the regions $\{x_0\le -1\}$, $\{|x_0|\le 1\}$, or $\{x_0\ge
1\}$. However, a disk with image in $\{x_0\ge 1\}$ ($\{x_0\le -1\}$)
must have its positive puncture at a Reeb chord $c$ with $c_0=1$
($c_0=-1$). The lemma follows from Lemma~\ref{lmamin}.
\end{pf}
Assume now that $\phi_\delta(L)$ is generic for each $\delta\ne
0$. Then it is a consequence of Lemma~6.25 in \cite{ees2} that any rigid
disk with boundary on $\Phi_f^\delta$ and positive corner at some Reeb
chord $c$ with $c_0=\pm 1$ is transversely cut out. Moreover by Lemma~6.12
in \cite{ees2} transversality of the $\bar\pa$-equation can be achieved
by perturbation near the positive puncture of a disk and it follows
that there exists (arbitrarily small) perturbations of $\Phi_f^\delta$ which
are supported in the region $\{|x_0|<1\}$ and which makes every moduli space
(of formal dimension $\le 1$) transversely cut out. We fix such a
perturbation of $\Phi_f^\delta$ but keep the notation $\Phi_f^\delta$
for the perturbed Legendrian embedding.
Let $\A(\Phi_f^\delta)$ denote the algebra over
$\Z[H_1(\R\times L)]=\Z[H_1(L)]$ generated by the Reeb chords of
$\Phi_f^\delta$ as in Subsection~\ref{defofal} and define the map (differential)
$\pa$ of $\A(\Phi_f^\delta)$ as there.
\begin{lma}\label{lmaalg'}
The map $\pa\colon \A(\Phi_f^\delta)\to\A(\Phi_f^\delta)$ satisfies
$\pa\circ\pa=0$.
\end{lma}
\begin{pf}
In the light of Lemma \ref{lmaalg},
a word by word repetition of the proof of Theorem \ref{thmd^2=0}
establishes the lemma.
\end{pf}

\subsubsection{Invariance under handle slides}
Let $\phi_t\colon L\to\C^n\times\R$, $-1\le t\le1$ be a Legendrian
isotopy such that $L_0$ is a generic handle slide moment. That is,
there exists one handle slide disk in some $\M_A(a;{\mathbf b})$,
which is the only non-empty moduli space of formal negative dimension,
that all moduli spaces of holomorphic disk with boundary on $\phi_t(L)=L_t$,
$t\ne 0$ of negative formal dimension are empty, and that all moduli
spaces of rigid disks are transversally cut out. We choose notation so
that $\{b_1,\dots,b_r,a,c_1,\dots,c_s\}$ are the
Reeb chords of $L_0$ and so that
$$
\ZZ(b_1)\le\dots\le\ZZ(b_r)\le\ZZ(a)\le\ZZ(c_1)\le\dots\le\ZZ(c_s).
$$
Recall for a Reeb chord $c,$ $\ZZ(c)$ is the difference in $z$-coordinates of its endpoints
in $\R^{2n+1}.$
Let $f\colon\R\to\R$ be a positive Morse function with local minima at
$\pm 1$, no critical points in the region
$(-\infty,-1)\cup(1,\infty)$, and one local maximum at $0$.
\begin{lma}\label{lmaRchords}
For all sufficiently small $\delta>0$ the Reeb chords of
$\Phi_f^\delta$ are
$$
\Bigl\{b_j[-1], b_j[1], b_j[0]\Bigr\}_{j=1}^r\cup
\Bigl\{a[-1], a[1], a[0]\Bigr\}\cup
\Bigl\{c_j[-1], c_j[1], c_j[0]\Bigr\}_{j=1}^s,
$$
where for any Reeb chord $c$ of $L_0$, $|c[-1]|=|c[1]|=|c[0]|-1=|c|$.
\end{lma}
\begin{pf}
It is easy to see that for $\delta=0$ the Reeb chords are as described
above and that the corresponding double points in $\C\times\C^n$ are
transverse. This shows that the Reeb chords are as claimed for all
sufficiently small $\delta$.
The second statement in the lemma is a straightforward consequence of
the following grading formula from \cite{ees1}.
Let $c^\pm$ be the $z$-coordinates of the upper and lower points in the front projection corresponding to
the Reeb chord $c.$ Assume $c^+$ is above $c^-$. Near $c^+$ we can represent the front as the graph of a function
$h_+: U\to \R,$ $U\subset \R^n$ with $h_+(q)=c^+$ for some $q\in U.$
We can similarly find a function $h_-$ for the front near $c^-.$
Let $h=h_+-h_-.$ Since $c^\pm$ correspond to a double point in the Lagrangian projection,
$h$ has a critical point at $q.$ If $c$ is a transverse double point $q$ is a non-degenerate
critical point. From \cite{ees1} we have
\begin{equation}\label{eq:frontform}
|c|=D(\gamma)-U(\gamma)+ \text{Morse Index}(d^2 h) - 1,
\end{equation}
where $\gamma$ is a path in the front connecting $c^+$ to $c^-$ and $D(\gamma)$ and $U(\gamma)$ is
the number of down- and up-cusps of $\gamma.$
\end{pf}
We call $b_j[0]$, $a[0]$, and $c_j[0]$, $[0]$-Reeb chords, and
$b_j[\pm 1]$, $a[\pm 1]$, and $c_j[\pm 1]$, $[\pm 1]$-Reeb chords. As
above we perturb $\Phi_f^\delta$ slightly in the region ${|x_0|<1}$ to
make it generic with respect to holomorphic disks. Note that the
$x_0$-coordinate of $[\pm 1]$-Reeb chord equals $\pm 1$ and that the
$x_0$-coordinate of a $[0]$-Reeb chord is very close to $0$ for small
$\delta>0$.
Consider a sequence of functions $f_k$ as above with $f_k\to 1$ as
$k\to\infty$ (i.e. each $f_k$ has a non-degenerate maximum at $0$ and
non-degenerate local minima at $\pm 1$). Fix $k$ and pick $\delta>0$
sufficiently small so that $\Phi_{f_k}^{\delta}$ satisfies
Lemma \ref{lmaRchords}. Let $\Phi_k^\delta=\Phi_{f_k}^{\delta}$.

We next note that as $\delta\to 0$, $\Phi_k^\delta\to\Phi_k^0$
where
$$
\Phi_k^0(t,q)=\Bigl(t,f_k'(t) z(q),x(q),f_k(t)y(q), f_k(t)z(q)\Bigr),
$$
with $(x(q),y(q),z(q))=\phi_0(q)$.
\begin{lma}\label{lmaconnstrip}
There exists $k_0$ such that for all $k>k_0$ there exists a $\delta_k>0$
such that for all $\delta<\delta_k$ and any Reeb chord $c$ the
following holds. The moduli spaces $\M(c[0],c[1])$ and
$\M(c[0],c[-1])$ of holomorphic
disks with boundary on $\Phi_k^\delta$ consists of exactly one point
which is a transversely cut out rigid disk.
Moreover the sign of the rigid disk in $\M(c[0],c[1])$ and that of the
disk in $\M(c[0],c[-1])$ are opposite.
\end{lma}
\begin{pf}
First consider the case $\delta=0$. It is easy to find rigid disks in
the $(x_0,y_0)$-plane with positive puncture at $c[0]$ and negative
puncture at $c[\pm 1]$. Moreover, by \cite{ees2} Lemma 6.25 these disks are
transversely cut out.
To see that these are the only disks, let $U$ and
$V$ be neighborhoods of the endpoints of the Reeb chord $c$ in $L_0$
and consider the projections of $\Phi_k^0([-1,1]\times U)$ and
$\Phi_k^0([-1,1]\times V)$ to $\C^n$. For sufficiently large $k$,
these projections intersect only at $0$ and it follows that there
exists a positive $h>0$ such that the area of the projection of any
disk with boundary on $\Phi_k^0$, positive puncture at $c[0]$, and
negative at $c[\pm 1]$ is either equal to zero or larger than
$h$. Since $\ZZ(c[0])\to\ZZ(c[\pm 1])$ as $k\to\infty$ it follows that
for $k$ large enough the disks in the $(x_0,y_0)$-plane are the only
ones.
We also check the statement about
signs in the case $\delta=0$. To this end, note that the trivialized
boundary conditions of the two disks in the $(x_0,y_0)$-plane are
identical and that multiplication by $-1$ is a
holomorphic automorphism relating them. Since multiplication by $-1$
reverses the orientation of the kernel of the linearized problem, it
follows that their signs are opposite.
Finally, we note that the fact that the moduli space
$\M(c[0],c[\pm 1])$ corresponding to $\Phi_k^0$ is transversely cut
out implies that the statement of the lemma holds also for
$\Phi_k^\delta$ for all sufficiently small $\delta$ (where the
smallness depends on $k$).
\end{pf}
We next note that as $k\to\infty$, $\Phi_k^0$ approaches the
Legendrian submanifold
$$
\Phi(t,q)=(t,0,x(q),y(q),z(q)).
$$
The projection of this Legendrian submanifold to $\C$ is simply the
$x_0$-axis and its projection to $\C^n$ agrees with that of $L_0$.
\begin{lma}\label{lmacn=a}
There exists $k_0$ such that for all $k>k_0$ there exists a $\delta_k>0$
such that for all $\delta<\delta_k$ and any Reeb chord $c\ne a$ the
following holds. If the moduli space $\M_A(c[0];{\mathbf e})$, where
${\mathbf e}$ is a word constant in the $[0]$-generators and
${\mathbf e}\ne c[\pm 1]$, has formal dimension $0$ then it is empty.
\end{lma}
\begin{pf}
Again we start with the case $\delta=0$. Consider a disk $u$ as above
with boundary on $\Phi_k^0$. As $k\to\infty$, $\Phi_k^0\to\Phi$
and the projection of $u$ converges to a broken disk $\{v^j\}_{j=1}^m$
with boundary on
$L_0$. The components $v^j$ of such a broken disk either have formal
dimension at least $0$, or equals the handle slide disk. Also, any
Reeb chord $b$ appearing as a puncture of some $v^j$ has
$\ZZ(b)\le\ZZ(c)$ and exactly one component of the broken disk must
have its positive puncture at $c\ne a$. This component has formal
dimension at least $0$ (since it is not the handle slide disk). For
a disk $v^j$ let $|v^j_+|$ is the grading of its positive puncture and
$|v^j_-|$ the sum of gradings of its negative punctures and the
negative of the grading of the homology data. Then the
formal dimension of (the moduli space of) $v^j$ is
$|v^j_+|-|v^j_-|-1$. The above implies
that
$$
N=\sum_{j=1}^m(|v^j_+|-|v^j_-|)\ge 1.
$$
Since the positive puncture of $u$ is its only $[0]$-puncture it
follows that the formal dimension of $u$ equals $N$.
The statement of the lemma follows for $\delta=0$. Since emptyness of
a moduli space is an open condition the lemma follows in general.
\end{pf}

Let $\Omega$ be the map from $\A(\Phi_k^\delta)$ to $\A(L_0)$ which maps
$c[\pm 1]$ to $c$ and $c[0]$ to $0$ for any Reeb chord $c$ of $L$.
\begin{lma}\label{lmac=a}
There exists $k_0$ such that for all $k>k_0$ there exists a $\delta_k>0$
such that for all $\delta<\delta_k$ the following holds.
If $u$ is a holomorphic disk with boundary on $\Phi_k^\delta$
in $\M_C(a[0],{\mathbf e})$, where
${\mathbf e}$ is a word constant in the
$[0]$-generators and ${\mathbf e}\ne a[\pm 1]$, and if this moduli
space has formal dimension $0$
then $C=A$ and $\Omega{\mathbf e}={\mathbf b}$.
\end{lma}
\begin{pf}
Consider first the case $\delta=0$. Taking the limit as $k\to\infty$
and arguing as in the proof of Lemma \ref{lmacn=a} we see that the
projection of $u$ converges to a broken disk $\{v^j\}_{j=1}^m$, that
all Reeb chords $b$ appearing
as a puncture of some $v^j$ satisfies $\ZZ(a)\ge\ZZ(b)$, and that
there is a unique component with its positive puncture at $a$. If this
component is not the handle slide disk then the argument in the proof
of Lemma \ref{lmacn=a} shows that the formal dimension of $u$ is at least
$1$. If, on the other hand, this component is the handle slide disk then
the formal dimension of $u$ equals $0$ only if the broken disk has no
other components. This shows the lemma for $\delta=0$. Again since the
condition that a moduli space is empty is open the lemma follows in
general.
\end{pf}
Fix $k$ sufficiently large and $\delta>0$ sufficiently small so that Lemmas
\ref{lmaconnstrip}, \ref{lmacn=a}, and \ref{lmac=a} holds for
$\Phi_k^\delta$. We also assume that $\Phi_k^\delta$ is generic with
respect to holomorphic disks. Let $\Phi=\Phi_k^\delta$
Let $\hat\A=\A(\Phi)$. We denote the differential of $\hat\A$ by
$\Delta$, see Lemma \ref{lmaalg'}. There are natural inclusions
$\A_\pm=\A(L_{\pm\delta})\subset\hat\A$.
Lemma \ref{lmamin} implies that this is an inclusion of DGA's in other
words,
$$
\Delta c[\pm 1]= \Gamma_{\pm}(\pa_{\pm} c),
$$
where $\Gamma_{\pm}:\A_\pm\to\hat \A$ is the map defined on generators by
$\Gamma_\pm(c)=c[\pm 1]$, and where $\pa_\pm$ is the differential on
$\A_\pm$.
For generators $b_j[0]$ we have by Lemmas \ref{lmaconnstrip} and \ref{lmacn=a}
\begin{equation}\label{eqDeltab[0]}
\Delta b_j[0]= b_j[1]- b_j[-1] + \beta_1^j + \Ordo(2),
\end{equation}
where $\beta_1^j$ is linear in the $c[0]$-generators and $\Ordo(2)$
denotes a linear combination of monomials which are at least quadratic
in the $[0]$-generators.
For the generator $a[0]$ we have by Lemmas \ref{lmaconnstrip} and \ref{lmac=a}
\begin{equation}\label{eqDeltaa[0]}
\Delta a[0]= a[1]-a[-1] + \epsilon + \alpha_1 + \Ordo(2),
\end{equation}
where $\Omega(\epsilon)=mA{\mathbf b}$, where $m\in\Z$ and where
$\alpha_1$  is linear in the $[0]$-generators.
For generators $c_j[0]$ we have by Lemmas \ref{lmaconnstrip} and \ref{lmacn=a}
\begin{equation}\label{eqDeltac[0]}
\Delta c_j[0]= c_j[1]-c_j[-1] + \gamma_1^j + \delta_1^j(a[0])+ \Ordo(2),
\end{equation}
where $\gamma_1^j+\delta_1^j(a[0])$ is linear in the
$[0]$-generators, where $\delta_1^j(a[0])$ lies in the ideal
generated by $a[0]$, and where $\gamma_1^j$ is constant in the
$a[0]$ generator.
Below we will consider $\pa_+$ and $\pa_-$ as different differentials
on the algebra $\A$. Let $\epsilon$ be as in \eqref{eqDeltaa[0]} and write
$\theta=\Omega(\epsilon)$. Consider the tame isomorphism $\psi$
of $\A$ defined on generators as
$$
\psi(c)=
\begin{cases}
c &\text{ if $c\ne a$},\\
a+\theta &\text{ if $c=a$}.
\end{cases}
$$
If $v\in A$ and $c$ is a generator of $\A$ then let
$$
\left(\frac{v}{c}\right)\bullet\colon \A\to\A
$$
be the map defined on monomials by replacing each occurrence of $c$ by
$v$. Then with $\pa_+^\psi=\psi^{-1}\circ\pa_+\circ\psi$ denoting the
induced differential, a straightforward calculation gives
\begin{align}
\pa_+^\psi(c)=
\begin{cases}
\left(\frac{a-\theta}{a}\right)
\bullet(\pa_+ c) &\text{ if $c\ne a$},\\
\pa_+ a +\pa_+\theta &\text{ if $c=a$}.
\end{cases}
\end{align}
\begin{lma}\label{lmahsinv}
The algebra $(\A,\pa_-)$ is isomorphic to the algebra
$(\A,\pa_+^\psi)$.
\end{lma}
\begin{pf}
We prove that the two differentials agree on generators.
By Lemma \ref{lmaalg'}, $\Delta^2=0$. Thus, summing the terms constant
in the $[0]$-generators after acting by $\Delta$
in \eqref{eqDeltab[0]} we find
\begin{equation}\label{eqDelta0b[0]}
0=\pa_+b_j[1] -\pa_-b_j[-1]
+ \left(\Delta\beta_1\right)_0,
\end{equation}
where $\left(\Delta\gamma_1\right)_0$ denotes the part of
$\Delta\gamma_1$ which is constant in the $[0]$-generators. Since the
constant part of $\Delta b_k[0]$ equals $b_k[1]-b_k[-1]$ it follows
that
$$
\Omega\left(\Delta\beta_1\right)_0=0.
$$
Therefore, applying $\Omega$ in \eqref{eqDelta0b[0]}, we conclude
\begin{equation}\label{eqdeltab}
\pa_- b_j =\pa_+ b_j =\left(\frac{a-\theta}{a}\right)\bullet(\pa_+ b_j),
\end{equation}
since no monomial in $\pa_+ b_j$ contains $a$.
Applying $\Delta$ to \eqref{eqDeltaa[0]} we find similarly
\begin{equation}\label{eqDelta0a[0]}
\pa_- a[-1] =\pa_+ a[1] +\Delta\epsilon +\left(\Delta\alpha_1\right)_0.
\end{equation}
The first equality in \eqref{eqdeltab} implies that
$$
\Omega(\Delta\epsilon)=\pa_+\theta.
$$
Since every $[0]$-generator in
$\alpha_1$ is for the form $b_j[0]$ we find, as
with $\beta_1$ above, that $\Omega\left(\Delta\alpha_1\right)_0=0$.
We conclude
\begin{equation}\label{eqdeltaa}
\pa_- a =\pa_+ a +\pa_+\theta.
\end{equation}
Applying $\Delta$ to \eqref{eqDeltac[0]} gives
\begin{equation}\label{eqDelta0c[0]}
\pa_-c_j[-1]=\pa_+c_j[1] + \left(\Delta\gamma_1^j\right)_0 +
\left(\Delta\delta_1^j(a[0])\right)_0.
\end{equation}
Applying $\left(\frac{a[-1]-\epsilon}{a[1]}\right)\bullet$ to both
sides in \eqref{eqDelta0c[0]} and noting that no monomial in
$\pa_-c_j[-1]$ contains an $a[1]$ generator we get
\begin{equation}\label{eqhsc}
\begin{split}
\pa_- c_j[-1] =&\left(\frac{a[-1]-\epsilon}{a[1]}\right)\bullet(\pa_+c_j[1])+
\left(\frac{a[-1]-\epsilon}{a[1]}\right)\bullet
\left(\Delta\gamma_1^j\right)_0\\
&+\left(\frac{a[-1]-\epsilon}{a[1]}\right)\bullet
\left(\Delta\delta_1^j(a[0])\right)_0.
\end{split}
\end{equation}
Each term in $\left(\Delta\delta_1^j(a[0])\right)_0$
arises by replacing $a[0]$ in every monomial $\xi a[0]\eta$ of
$\delta_1^j(a[0])$ with $(a[1]-a[-1]+\epsilon)$ yielding
$\xi(a[1]-a[-1]+\epsilon)\eta$. When
$\left(\frac{a[-1]-\epsilon}{a[1]}\right)\bullet$ is a applied to
$\xi(a[1]-a[-1]+\epsilon)\eta$ the result is
$$
\xi(a[-1]-\epsilon-a[-1]+\epsilon)\eta=0.
$$
Thus, the last term in \eqref{eqhsc} vanishes.
Since the $[0]$-generator of any monomial in $\gamma_1^j$ equals
either $c_k[0]$ for some $k$, or $b_r[0]$ for some $r$ and since the
constant part of $\Delta c_k[0]$ equals $c_k[1]-c_k[-1]$, and the
constant part of $\Delta b_k[0]$ equals $b_k[1]-b_k[-1]$, we conclude
that
$$
\Omega\left(\frac{a[-1]-\epsilon}{a[1]}\right)\bullet
\left(\Delta\gamma_1^j\right)_0=0.
$$
Thus, applying $\Omega$ to \eqref{eqhsc} we arrive at
\begin{equation}\label{eqdeltac}
\pa_- c_j=\left(\frac{a-\theta}{a}\right)\bullet\pa_+c_j.
\end{equation}
The lemma follows from \eqref{eqdeltab}, \eqref{eqdeltaa},
\eqref{eqdeltac}.
\end{pf}
\begin{cor}\label{corhsinv}
If $L_t\subset\C^n\times\R$, $-1\le t\le 1$, is a Legendrian isotopy
with a generic handle slide at $t=0$ as above then
the stable tame isomorphism classes of $(\A(L_{-1},\pa_{-1}))$ and $(\A(L_1),\pa_1)$ are the
same. \qed
\end{cor}

\subsubsection{Invariance under self tangencies}
We now turn our attention to self tangency instances. As shown in
\cite{ees2} there is no loss of generality in assuming that our isotopy
near the self tangency instant has standard from. Consider a
$1$-parameter family of Legendrian submanifolds $L_t\subset\C^n\times\R$,
$t\in[-1,1]$. We assume that $L_t$ is constant in $t$
outside $B_{2r}(0)\times\R$ for some ball $B_{2r}(0)\subset\C^n$
or radius $2r$ around the origin. The intersection of $L_t$ with
$B_r(0)\times\R$ consists of two sheets $L_t^+$ and $L_t^-$. Let
$(x,y,z)=(x+iy,z)$,
$x,y\in\R^n$, and $z\in\R$ be coordinates on $\C^n\times\R$. Assume
that the sheet $L_t^-$ is constant in $t$ and that it satisfies
$$
L_t^-=B_{2r}(0)\times\R\cap\{y=0,z=0\}.
$$
The second sheet is moving with $t$. Let
$x=(x_1,x_2)\in\R\times\R^{n-1}$,
$y=(y_1,y_2)\in\R\times\R^{n-1}$, and let $\la\,,\ra$ denote the
standard inner product on $\R^k$. Then $L_t^+\cap B_r(0)\times\R$ is
given by the map
$[-1,1]\times\R\times\R^{n-1}\to\C^n\times\R$,
$(t,q_1,q_2)\mapsto(x_1+iy_1,x_2+iy_2,z)$ where
\begin{align*}
x_1(t,q_1,q_2)&=q_1,\\
y_1(t,q_1,q_2)&=3q_1^2+ t,\\
z(t,q_1,q_2)&=q_1(q_1^2+ t)+c+\la q_2, q_2\ra,\\
x_2(t,q_1,q_2)&= q_2,\\
y_2(t,q_1,q_2)&= 2q_2,
\end{align*}
where $c>0$ is a constant. In the region
$L_t^+\cap (B_{2r}(0)\setminus B_r(0))\times\R$ the isotopy
interpolates by the one above and the constant isotopy.
We denote this $1$-parameter family of
Legendrian embeddings $\phi_t\colon L\to\C^n\times\R$,
$t\in[-1,1]$.
Let $\delta>0$. Note that $\phi_{-\delta}$ has two Reeb chords
more than $\phi_{\delta}$ and that the extra Reeb chords of
$\phi_{-\delta}$ converges to the self tangency Reeb chord $o$ of
$\phi_0$. We choose notation so that the Reeb chords of
$\phi_{\delta}$ are
$$
\left\{b_1,\dots,b_s,a_1,\dots,a_l\right\},
$$
those of $\phi_0$ are
$$
\left\{b_1,\dots,b_s,o,a_1,\dots,a_l\right\},
$$
and those of $\phi_{-\delta}$ are
$$
\left\{b_1,\dots,b_s,b,a,a_1,\dots,a_l\right\},
$$
and so that
$$
\ZZ(b_1)\le\dots\le\ZZ(b_r)\le\ZZ(b)\le
\ZZ(o)\le\ZZ(a)\le\ZZ(a_1)\le\dots\le\ZZ(a_s).
$$
Fix a positive Morse function $f$ with local minima at $\pm 1$ and no
critical points in $(-\infty,-1)\cup(1,\infty)$ and with one local
maximum at $\beta\in(0,1)$, for some small $\beta$.
Assume that $\phi_0(L)$ is a generic self tangency moment
which means that all moduli spaces of negative formal dimension are
empty and that all rigid disks with boundary on $\phi_0(L)$ are
transversely cut out.
\begin{lma}\label{lmastReebchords}
Fix $f$ as above. Then there exists $\delta_0$ such that for all
$\delta<\delta_0$ the Reeb chords of $\Phi^\delta_f$ are
$$
\{b_j[-1],b_j[1],b_j[0]\}_{j=1}^r\cup
\{b[-1],a[-1]\}\cup
\{a_j[-1],a_j[1],a_j[0]\}_{j=1}^s,
$$
where the $x_0$-coordinates of the $[\pm 1]$-Reeb chords are $\pm 1$, and
where the $x_0$-coordinates of the $[0]$-Reeb chords are close to $0$.
\end{lma}
\begin{pf}
Consider first the case $\delta=0$. The Reeb chords of $\Phi^0_f$ are
easily seen to be
$$
\{b_j[-1],b_j[1],b_j[0]\}_{j=1}^r\cup
\{o[-1],o[1],o[0]\}\cup
\{a_j[-1],a_j[+1],a_j[0]\}_{j=1}^s,
$$
and all except $o[\pm1]$ and $o[0]$ have transverse tangent planes at
their ends. We conclude that for $\delta>0$ sufficiently small
$\Phi^\delta_f$ has the Reeb chords
$$
\{b_j[-1],b_j[1],b_j[0]\}_{j=1}^r\cup
\{a_j[-1],a_j[+1],a_j[0]\}_{j=1}^s,
$$
and possibly some Reeb chords in a neighborhood of $o[\pm 1]$ and of
$o[0]$. To find these we take a closer look at $\Phi^\delta_f.$
If $o^-$ denotes the lower endpoint of $o$ then
in a neighborhood of $\{o^-\}\times\R$, $\Phi^\delta_f$ is
simply the embedding (with notation as above)
$$
\Phi(t,q)=(t,0,q,0,0).
$$
In neighborhoods of $(\pm 1,\{o^+\})$, $\Phi^\delta_f$ is given by
\begin{align*}
x_0(t,q_1,q_2)&=t,\\
y_0(t,q_1,q_2)&=f'(t)\Bigl(q_1(q_1^2\pm\delta)+c
+\la q_2, q_2\ra\Bigr)\\
x_1(t,q_1,q_2)&=q_1,\\
y_1(t,q_1,q_2)&=f(t)(3q_1^2\pm\delta),\\
x_2(t,q_1,q_2)&= q_2,\\
y_2(t,q_1,q_2)&= f(t)2q_2\\
z(t,q_1,q_2)&=f(t)\Bigl(q_1(q^2_1\pm\delta)+c+\la q_2,q_2\ra\Bigr),
\end{align*}
and we find that there are two Reeb chords $a[-1]$ and $b[-1]$ in a
neighborhood of $o[-1]$ and no Reeb chords in a neighborhood of
$o[1]$.
In a neighborhood of $(0,\{o^+\})$ it is given by
\begin{align*}
x_0(t,q_1,q_2)&=t,\\
y_0(t,q_1,q_2)&=f'(t)\Bigl(q_1(q_1^2+\delta t)+c
+\la q_2, q_2\ra\Bigr)+f(t)\delta q_1\\
x_1(t,q_1,q_2)&=q_1,\\
y_1(t,q_1,q_2)&=f(t)(3q_1^2+\delta t),\\
x_2(t,q_1,q_2)&= q_2,\\
y_2(t,q_1,q_2)&= f(t)2q_2\\
z(t,q_1,q_2)&=f(t)\Bigl(q_1(q^2_1+\delta t)+c+\la q_2,q_2\ra\Bigr),
\end{align*}
To have a Reeb chord we note first that the equation $3q^2+\delta t=0$
must hold. Since $\delta>0$, the $t$-coordinate of any Reeb chord thus
satisfies $t<0$. Moreover, $|t|<1$ implies that at a Reeb chord
$q_1=\Ordo(\delta^{\frac12})$ and from the expression for $y_2$,
$q_2=0$. The final condition is $y_0=0$, which implies
$$
f'(t)(c+\Ordo(\delta^{\frac32})) + f(t)\Ordo(\delta^{\frac32})=0.
$$
By the choice of $f$, the local maximum of $f$ lies at $t=\beta>0$
thus $f'(t)>0$ for $-1<t\le 0$ and letting $\delta\to 0$ we see that
$\Phi_f^\delta$ does not have any Reeb chords near $o[0]$. This
finishes the proof.
\end{pf}
Consider a sequence of functions $f_k$ as above with $f_k\to 1$ as
$k\to\infty$ (i.e. each $f_k$ has a non-degenerate maximum at $\beta_k>0$ and
non-degenerate local minima at $\pm 1$). Fix $k$ and pick $\delta>0$
sufficiently small so that $\Phi_{f_k}^{\delta}$ satisfies
Lemma \ref{lmastReebchords}. Let $\Phi_k^\delta=\Phi_{f_k}^{\delta}$.
\begin{lma}\label{lmastconnstrip}
There exists $k_0$ such that for all $k>k_0$ there exists a $\delta_k>0$
such that for all $\delta<\delta_k$ and a Reeb chord
$c\in\{b_j\}_{j=1}^r\cup \{a_k\}_{k=1}^s$  the
following holds. The moduli spaces $\M(c[0],c[1])$ and
$\M(c[0],c[-1])$ of holomorphic
disks with boundary on $\Phi_k^\delta$ consists of exactly one point
which is a transversely cut out rigid disk.
Moreover the sign of the rigid disk in $\M(c[0],c[1])$ and that of the
disk in $\M(c[0],c[-1])$ are opposite.
\end{lma}
\begin{pf}
The proof is a word by word repetition of the proof of Lemma
\ref{lmaconnstrip}.
\end{pf}
We next note that as $k\to\infty$, $\Phi_k^0$ approaches the
Legendrian submanifold
$$
\Phi(t,q)=(t,0,x(q),y(q),z(q)).
$$
The projection of this Legendrian submanifold to $\C$ is simply the
$x_0$-axis and its projection to $\C^n$ agrees with that of $L_0$.
\begin{lma}\label{lmastc=b}
There exists $k_0$ such that for all $k>k_0$ there exists a $\delta_k>0$
such that for all $\delta<\delta_k$ and any Reeb chord $b_j$,
$j=1,\dots,r$, the following holds. If the moduli space
$\M_A(b_j[0];{\mathbf e})$, where ${\mathbf e}$ is a word constant in
the $[0]$-generators and ${\mathbf e}\ne b_j[\pm 1]$, has formal
dimension $0$ then it is empty.
\end{lma}
\begin{pf}
For $k$ sufficiently large,
$\ZZ(b_j[0])<\ZZ(a[-1])<\ZZ(b[-1])$. Therefore, for such $k$, a disk
with its positive puncture at $b_j[0]$ must have negative punctures
mapping to Reeb chords in the set $\{b_i[0],b_i[\pm
1]\}_{i=1}^s$. After this observation the lemma follows from the proof
of Lemma \ref{lmacn=a}.
\end{pf}
\begin{lma}\label{lmastc=a}
There exists $k_0$ such that for all $k>k_0$ there exists a $\delta_k>0$
such that for all $\delta<\delta_k$ and any Reeb chord $a_j$,
$j=1,\dots,s$ the following holds. If the moduli space
$\M_A(a_j[0];{\mathbf e})$, where ${\mathbf e}$ is a word constant in
the $[0]$-generators and ${\mathbf e}\ne a_j[\pm 1]$, has formal
dimension $0$ then $a[-1]$ appears at least once as a letter
in the word ${\mathbf e}$.
\end{lma}
\begin{pf}
Let $u$ be a disk in $\M_A(a_j[0];{\mathbf e})$. As in the proof of
Lemma \ref{lmac=a} its projection to $\C^n$ converges to a broken disk with
boundary on $L_0$ as $k\to\infty$. Let $\{v^j\}_{j=1}^m$ be the
components of this broken disk.
As above we let $o$ denote the degenerate Reeb chord of $L_0$. Let
$|o|=|b[-1]|=|a[-1]|-1$.  Then the formal
dimension of a disk $v_j$ satisfies the following.
\begin{itemize}
\item If the positive puncture of $v_j$ does not equal $o$ then the formal
dimension of $v_j$ equals $|v_j^+|-|v_j^-|-1$, where $|v_j^+|$ is the grading
of the Reeb chord at its positive puncture, and where $|v_j^-|$ is the
sum of the gradings at its negative corners and the negative of
the grading of its homology data.
\item If the positive puncture of $v_j$ equals $o$ then the formal
dimension of $v_j$ equals $|v_j^+|-|v_j^-|.$
\end{itemize}
Since each of the disks $v_j$ have non-negative formal dimension and
since at least one of them does not have its positive puncture at $o$,
we find
$$
\sum_{j=1}^m(|v_j^+|-|v_j^-|)\ge 1.
$$
Thus, assuming that none of the negative punctures of $u$ map to
$a[-1]$ we find that the formal dimension of $u$ is at least $1$ which
contradicts it being rigid. It follows that at least one of the negative
punctures of $u$ map to $a[-1]$.
\end{pf}

Fix $k$ sufficiently large and $\delta>0$ sufficiently small so that
Lemmas \ref{lmastc=b} and \ref{lmastc=a} hold for $\Phi^\delta_k$. Let
$\Phi=\Phi^\delta_k$.
Consider the algebra $\hat\A=\A(\Phi)$ and let its differential be
$\Delta$. Note that
$\A_+=\A(L_{\delta})$ and $\A_-=\A(L_{-\delta})$ can be considered as
subalgebras of $\hat\A$ via the map which takes $a_j$ and $b_j$ to
$a_j[\pm 1]$ and $b_j[\pm 1]$ and $a$ and $b$ to $a[-1]$ and $b[-1]$,
respectively. It is a consequence of Lemma \ref{lmamin} that $\A_\pm$ are
differential graded subalgebras of $\hat A$. In other words, denoting
their respective differentials $\pa_+$ and $\pa_-$ we have
\begin{align*}
\Delta(b_j[\pm 1]) &=\Gamma_\pm(\pa_\pm b_j),\quad\text{ for }j=1,\dots,r,\\
\Delta(a[-1]) &=\Gamma_-(\pa_-a), &\\
\Delta(b[-1]) &=\Gamma_-(\pa_-b),&\\
\Delta(a_j[\pm 1]) &=\Gamma_\pm(\pa_\pm a_j),\quad\text{ for }j=1,\dots,s,
\end{align*}
where $\Gamma_\pm a_j=a_j[\pm 1]$, $\Gamma_\pm b_j=b_j[\pm 1]$,
$\Gamma_-a=a[-1]$, and $\Gamma_-b=b[-1]$.
It follows from Lemma \ref{lmastc=b} that
\begin{equation}\label{eqstDb}
\Delta b_j[0]=b_j[1]-b_j[-1]+\beta_1^j+\Ordo(2),
\end{equation}
where $\beta_1^j$ denotes term which is linear in the $[0]$-generators
and $\Ordo(2)$ the term which is quadratic and higher.
It follows from Lemma \ref{lmastc=a} that
\begin{equation}\label{eqstDa}
\Delta a_j[0]=a_j[1]-a_j[-1]+\gamma(a[-1])+\alpha_1^j+\Ordo(2),
\end{equation}
where $\gamma(a[-1])$ lies in the ideal generated by $a[-1]$ and is
constant in the $[0]$-generators and where $\alpha_1^j$ is the linear
in the $[0]$-generators.
Consider the stabilized algebra $S(\A_+)$ with extra generators $e_0$
and $e_1$, $|e_0|=|a|$ and $|e_1|=|b|=|a|-1$ and define the algebra
homomorphism $\Phi_0\colon\A_-\to S(A_+)$
$$
\Phi_0(c)=
\begin{cases}
e_0 &\text{ if $c=a$},\\
e_1-v &\text{ if $c=b$},\\
c &\text{otherwise},
\end{cases}
$$
where $v$ is the unique element in $\A_+$ such that $\pa_- a= b + v$
(in $\A_-$). For the existence of such an element see \cite{ees2}, Lemma~1.16.
Let  $\tau\colon S(\A_+)\to\A_+$ be the natural
projection and let $\A_j\subset\A_-$ be the subalgebra generated by
$\{b_1,\dots,b_s,b,a,a_1,\dots,a_j\}$. Then
\begin{lma}\label{lmaalg1}
\begin{equation}\label{eqalg1}
\Phi_0\circ\pa_- w=\pa_+^s\circ\Phi_0 w,
\end{equation}
for all $w\in\A_0$ and
\begin{equation}\label{eqalg2}
\tau\circ\Phi_0\circ\pa_- =\tau\circ\pa_+^s\circ\Phi_0.
\end{equation}
\end{lma}
\begin{pf}
Let $\Omega$ be the map which takes $b_j[-1]$ and $b_j[1]$ to
$b_j$ and which takes $a_j[1]$ and $a_j[-1]$ to $a_j$.
Then, for generators $b_j$, \eqref{eqalg1} follows by
applying $\Omega\circ\Delta$ to \eqref{eqstDb}.
We also have
$$
\Phi_0\pa_- b=\Phi_0(-\pa_- v)=\pa_+^s\Phi_0 b,
$$
and
$$
\Phi_0\pa_- a=\Phi_0 (b + v)= e_1=\pa_+\Phi_0 a.
$$
Thus, \eqref{eqalg1} holds and it is sufficient to prove
\eqref{eqalg2} for $a_j$-generators to conclude it holds in general.
Applying $\Delta$ to \eqref{eqstDa} and considering the constant term
we find
$$
\pa_+a_j[1]=\pa_- a_j[-1] -\Delta(\gamma(a[-1]))
- (\Delta\alpha^j_1)_0.
$$
Letting $\tilde\Phi_0$ be the map which takes $a[-1]$ to $e_0$ and
$b[-1]$ to $e_1-\Gamma_-(v)$ we find that
$$
\pa_+ a_j[1]=\tilde\Phi_0\pa_-a_j[-1]
-\tilde\Phi_0(\Delta(\gamma (a[-1]))
-\tilde\Phi_0(\Delta\alpha^j_1)_0.
$$
Since the constant part of $\Delta b_j[0]$ equals $(b_j[1]-b_j[-1])$
we note that $\Omega$ annihilates each polynomial in
$(\Delta\alpha_1^j)_0$ which originates from a monomial in
$\alpha_1^j$ in the ideal generated by $b_j[0]$. Moreover, the
constant part of $\Delta a_j[0]$ equals
$a_j[1]-a_j[-1]+\gamma(a[-1])$. Since $\tau$ annihilates
$\Phi_0(\gamma(a[-1]))$ and since $\Omega$ annihilates
$a_j[1]-a_j[-1]$ and since $\tau$ and $\Omega$ commutes we find that
$\tau\circ\Omega$ annihilates the last term.
Finally, each term in $(\Delta(\gamma(a[-1]))$ contains either
$a[-1]$ or (when $a[-1]$ is differentiated) $b[-1]+\Gamma_-v$.
Hence, $\tau$ annihilates $\Phi_0(\Delta\gamma(a[-1]))$ and we
conclude
$$
\tau \pa_+ a_j =\tau\Phi_0\pa_- a_j,
$$
as claimed.
\end{pf}
\begin{lma}\label{ststi}
$\A(L_-)$ is stable tame isomorphic to $\A(L_+)$.
\end{lma}
Given Lemma~\ref{lmaalg1} the proof of this lemma is standard, see \cite{Chekanov, ees2, ENS}.
\subsection{Change of spin structure}\label{subsec:cos}
In our construction of Legendrian contact homology in Subsection~\ref{defofal}
we assume that our Legendrian manifolds $L\subset\C^n\times\R$
are spin and we fix a spin structure on $L$. In general the DGA
$(\A(L),\pa)$ of $L$ depends on the fixed spin structure. We explain
here the exact form of this dependence.
\subsubsection{Spin structures}
Let $M$ be an oriented $n$-manifold which is spin. In
Subsection~\ref{OMS} we viewed spin structures as trivializations of
${\tilde T} M$ over the $1$-skeleton that extends over the
$2$-skeleton. Given a spin structure $\mathfrak{s}_0$ on $M$ we can
identify the set of spin structures $\text{Spin}(M)$ with
$H^1(M;\Z_2).$ Specifically, let $\mathfrak{s}$ be another spin
structure. Denote the trivialization associated to $\mathfrak{s}$ by
$\tau_{\mathfrak{s}}.$ Isotop $\tau_{\mathfrak{s}_0}$ and
$\tau_{\mathfrak{s}}$ to be the same on the 0-skeleton
$\Delta^{(0)}$ of $M.$ Now along each edge in the 1-skeleton
$\Delta^{(1)}$ of $M$ we can compare $\tau_{\mathfrak{s}_0}$ to
$\tau_{\mathfrak{s}}.$ This gives us a loop in
$\pi_1(SO(n+2))=\Z_2.$ Thus we have a map
$d(\mathfrak{s}_0,\mathfrak{s}):\Delta^{(1)}\to \Z_2.$ One may
easily check that this gives a well defined cohomology class
$[d(\mathfrak{s}_0,\mathfrak{s})]\in H^1(M;\Z_2).$ It is a standard
fact that
\[d(\mathfrak{s}_0,\cdot):\text{Spin}(M)\to H^1(M;\Z_2)\]
is a one to one correspondence. We denote this map
$d_{\mathfrak{s}_0}.$
\subsubsection{The change of the differential}
Let $L\subset\C^n\times\R$ be a Legendrian submanifold with a fixed
spin structure $\mathfrak{s}_0$. Let $\mathfrak{s}$ be another spin structure on
$L$. Let $\pa_{\mathfrak{s}_0}$ and $\pa_{\mathfrak{s}}$ denote the
differentials on $\A(L)$
induced from the spin structures $\mathfrak{s}_0$ and $\mathfrak{s}$
respectively. Let
$A\in H_1(L)$.
\begin{thm}\label{changeos}
Let $a$ be any Reeb chord of $L$. Assume that
$$
\pa_{\mathfrak{s}_0} a = \sum_j m_j A_j{\mathbf b}_j,
$$
where $m_j\in \Z$, $A_j\in H_1(L)$, and where ${\mathbf b}_j$ is a Reeb
chord word. Then
$$
\pa_{\mathfrak{s}} a = \sum_j \sigma(\langle d_{\mathfrak{s}_0}({\mathfrak{s}}), A_j\rangle) m_j A_j{\mathbf b}_j,
$$
where $\sigma\colon\Z_2\to\{1,-1\}$ is the non-trivial homomorphism.
\end{thm}
\begin{pf}
Take a triangulation of $L$ containing all the double points in the
complex projection in the 0-skeleton
and containing the capping paths for the double points in the
1-skeleton. (We can assume the
capping paths are all disjoint arcs in $\Delta^{(1)}.$)
Note that we may find trivializations $\tau_{\mathfrak{s}_0}$ and
$\tau_{\mathfrak{s}}$
corresponding to $\mathfrak{s}_0$ and $\mathfrak{s}$, respectively, of
${\tilde T}L$ over $\Delta^{(2)}$ such that the trivializations agree over all
capping paths. Let $u$ be any rigid holomorphic disk. Consider the
trivializations of the Lagrangian boundary conditions on the closed
disk which arises when $u$ is capped off and which is induced from
$\tau_{\mathfrak{s}_0}$ and $\tau_{\mathfrak{s}}$ respectively. The
difference of these
trivializations arises in differences of the framings along the
boundary paths of $u$. Since $\tau_{\mathfrak{s}_0}$ and
$\tau_{\mathfrak{s}}$ agrees on the
capping paths it follows that this framing difference is exactly
$\sigma(\langle d_{\mathfrak{s}_0}({\mathfrak{s}}), A\rangle)$, where
$A$ is the homology class encoding
the boundary of $u$. By \cite{FOOO} the canonical orientation on the
determinant line of homotopically different
trivializations of a Lagrangian boundary condition on the closed disk
are opposite. The lemma follows.
\end{pf}
\subsubsection{Change in the DGA}
While it seems crucial that a Legendrian submanifold be spin to
define an oriented version of contact homology, there is no real
dependence on the spin structure,  when we define contact homology
over $\Z[H_1(L)].$
\begin{thm}
Let $L$ be a Legendrian submanifold of
$\R^{2n+1}.$ The DGA's of $L$ over $\Z[H_1(L)]$ associated to any two spin structures are tame
isomorphic.
\end{thm}
\begin{proof}
Given two spin structures $\mathfrak{s}$ and $\mathfrak{s}'.$ Let $\mathcal{A}$ and $\mathcal{A}'$
be the DGA's associated to $L$ using the two spin structures. Using the notation from
Theorem~\ref{changeos} define $\phi:\mathcal{A}\to\mathcal{A}'$ to be the identity of the generators
of the DGA but for $A\in H_1(L)$ let $\phi(A)=\sigma(\langle
d_{\mathfrak{s}}({\mathfrak{s}'}), A \rangle)A$ (thus the isomorphism $\phi$ arises from an
isomorphism of the base ring $\Z[H_1(L)]$). One may easily check that $\phi$ is a chain map and
hence a tame isomorphism from $\mathcal{A}$ to $\mathcal{A}'.$
\end{proof}
It is quite interesting to note that if one merely uses orientations to define the DGA for $L$ over
$\Z$ then there is a dependence on spin structures.
\begin{thm}
The DGA's of a Legendrian submanifold $L$ defined over $\Z$ using two different spin structures are
not necessarily stable tame isomorphic.
\end{thm}
\begin{proof}
If the DGA's are stable tame isomorphic then the contact homology associated to $L$ with the two
spin structures would be isomorphic. Let $L$ be the Legendrian unknot in $\R^3$ whose Lagrangian
projection has one double point. Thus $\mathcal{A}=\Z\langle a\rangle,$ where $a$ is the double point in the
projection. If we use the spin structure on $L$ defined in Subsection~\ref{3Drel} then
$\partial a=2$ so the contact homology is $\Z_2.$ By Theorem~\ref{changeos} we see that using the
other spin structure on $L$ will give a differential $\partial a=0.$ So the contact homology with
this spin structure is a copy of $\Z$ in gradings 0 and 1. Thus the DGA's associated to $L$ using
the two spin structures are not stable tame isomorphic. To get examples in higher dimensions one can
use the spinning construction, see \cite{ees1}.
\end{proof}
\subsection{The three dimensional case}\label{3Drel}
In this section we show that with the proper choice of spin
structure on $S^1$ the DGA we associate to a Legendrian knot in
$\R^3$ is the same as the combinatorially defined one given in
\cite{ENS}. We also deduce an alternative combinatorial description
and demonstrate that these two in a certain sense constitute a
complete list of possible combinatorial definitions.
Recall that our construction of orientations on the moduli spaces
relevant to contact homology depend on choices. Specifically (see
Section~\ref{COCD}) we chose an orientation on $\R^n,$ on the
capping operators, on spaces of conformal structures $\mathcal{C}_m$
and automorphisms $\A_m,$ and on $\C$. This last choice was largely
hidden in previous sections: we were simply using the natural
complex orientation $\C$. However, there is no real need to choose
this orientation and, as we shall see, the choice matters. We will
call the choices listed above {\em the choice of basic
orientations}.
\subsubsection{Combinatorial descriptions}
\begin{figure}[ht]
  \relabelbox \small {\epsfxsize=4.5in\centerline{\epsfbox{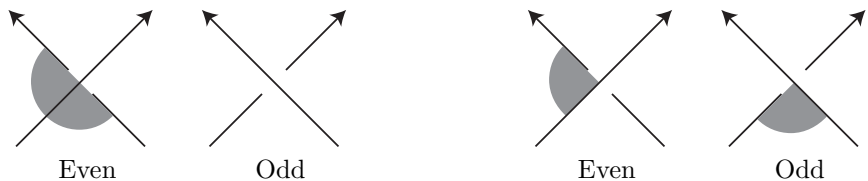}}}
  \relabel{1}{Even}
  \relabel{2}{Odd}
  \relabel{3}{Even}
  \relabel{4}{Odd}
  \endrelabelbox
        \caption{Two sign rules}
        \label{signrule}
\end{figure}

Let $L\subset\C\times\R$ be an oriented Legendrian knot and consider
a double point of its Lagrangian projection. Near this double point
the Lagrangian projection subdivides the plane into four quadrants.
We describe two shading rules.
\begin{itemize}
\item[{\rm (A)}]
A quadrant of a double point with even grading (see section \ref{subsec:orion1d})
is {\em A-shaded} if it is adjacent to the incoming edge of the
overcrossing, the other two quadrants are {\em A-unshaded}. All
quadrants of odd double points are {\em A-unshaded}. See the left
hand side of Figure~\ref{signrule}.
\item[{\rm (B)}]
A quadrant of a double point with even grading (see section \ref{subsec:orion1d})
is {\em B-shaded} if it is adjacent to the incoming edge of the
overcrossing and to the outgoing edge of the undercrossing, the
other three quadrants are {\em B-unshaded}. A quadrant of an odd
double points is B-shaded if it adjacent to the incoming edges of
both the over-  and the undercrossing, the other three are
B-unshaded. See the right hand side of Figure~\ref{signrule}.
\end{itemize}
The main goal of this section is to prove the following theorem.
\begin{thm}\label{thm:3Dsigns}
Let $L$ be an oriented  Legendrian knot in $\R^3$
equipped with the Lie group spin structure of $L=S^1$. Then there
exists a choice of basic orientations such that for any Reeb chord
$a$ of $L$,
\[\partial a = \sum \left( \sum_{u\in \M(a;{b_1\ldots b_n})} (-1)^{s_A(u)}
b_1\ldots b_n\right),\] where the first sum is over words $b_1
\ldots b_n$  in  the double points of the Lagrangian projection of
$L$ with $|a|-|b_1|-\dots-|b_n|=1$ and $s_A(u)$ is the number of
A-shaded corners in the image of $u.$
Moreover, there exists another choice of
basic orientations (where the orientation of $\C$ is opposite, for
more detail see Lemma \ref{lmacombor}) such that
\[\partial a = \sum \left( \sum_{u\in \M(a;{b_1\ldots b_n})} (-1)^{s_B(u)}
b_1\ldots b_n\right),\] where the first sum is over words $b_1
\ldots b_n$  in  the double points of the Lagrangian projection of
$L$ with $|a|-|b_1|-\dots-|b_n|=1$ and $s_B(u)$ is the number of
B-shaded corners in the image of $u.$
\end{thm}
\begin{rmk}
The sign rule presented in \cite{ENS} is the one corresponding to
the A-shading.
\end{rmk}
\begin{rmk}\label{rmk:3Dsigns}
Note that we have not explicitly identified which orientation on
$\C$ gives which orientation convention. We also note that other
changes in the basic orientations give differentials closely related
to those given here. (For reasonable definitions of basic
orientations they differ by an over-all sign.)
\end{rmk}
\begin{rmk}
To compute the differential with respect to the other spin structure
on $S^1$ (the null-cobordant one) one can appeal to
Theorem~\ref{changeos}. However, there is also a simple way to
compute the differential in this case. Start with the trivialization
of the stabilized tangent bundle to $S^1$ corresponding to the Lie
group spin structure and add a $\pi$-rotation in a small
neighborhood of a point $p \in S^1.$ We may assume that no capping
path used in computing the algebra and differential contains $p.$ If
$u\in\mathcal{M}(a;b_1,\ldots, b_n)$ then define $I(u)$ to be the
number of times $u(\partial D_{n+1})$ intersects $p.$ As in the
proof of Theorem~\ref{changeos} we see that the differential
$\partial'$ corresponding to the new spin structure is
\[\partial' a = \sum \left( \sum_{u\in \M(a;b_1,\ldots, b_n)}
(-1)^{s_\ast(u)+I(u)} b_1\ldots b_n\right),\] where $\ast=A$ or
$\ast=B$.
\end{rmk}
\subsubsection{Stabilization and trivialization}
Recall from Subsection~\ref{OMS} that when given a chord generic
Legendrian knot $L\subset\C\times\R$ we first consider the
Lagrangian projection $\Pi_\C$ of $L$ to $\C$ then the inclusion of
$\C$ into $\C^3$ as the first coordinate. In the $1$-dimensional
case the complex angle has only one component and therefore a
$1$-dimensional stabilization is sufficient. We thus consider a
$\C^2$ bundle over $L$ and a field of Lagrangian subspaces of the
bundle by assigning to each point $p$ in $L$ the Lagrangian subspace
\begin{equation}\label{dbforR3}
(t(p), e^{i\theta(p)})
\end{equation}
where $t(p)=\Pi_\C(\tau_p)$ for the unit tangent vector $\tau_p$ to
$L$ at $p$ and where $\theta:L\to [-\theta_0,\theta_0],$ for some
small $\theta_0,$ and $\theta(p)=\theta_0$ for $p$ in a neighborhood
of each upper end of a Reeb chord and $\theta(p)=-\theta_0$ in a
neighborhood of each lower end point of a Reeb chord.
The auxiliary linearized problem for a holomorphic disk
$u\colon D_{m+1}\to\C$ is then the $\bar\pa$-problem with boundary
condition given by the above plane field along $u(\pa D_{m+1})$. We note
that this auxiliary linearized problem is split:
$\bar\pa=\bar\pa_1\oplus\bar\pa_2$ where $\bar\pa_j$ acts on the
$j^{\rm th}$ coordinate of a section. Moreover,
$$
\ix(\bar\pa_2)=0,
$$
and $\bar\pa_2$ is an isomorphism.
To fit the above into the orientation scheme presented in previous
sections, we need a trivialization of the Lagrangian plane field
which meet the conditions presented in
Subsection~\ref{cappingrulegeq3} at the double points of $L$. (Note
that \eqref{dbforR3} gives a trivialization which does not
necessarily meet the conditions at Reeb chords.) To achieve this we
change the trivialization in \eqref{dbforR3} in a neighborhood of
the upper end of each {\em odd} Reeb chord as follows. Following the
orientation of the knot we add a $\pi$-rotation to the
trivialization right before we come to the upper Reeb chord end and
a $(-\pi)$-rotation right after (here we think of the Lagrangian
plane field as oriented by the trivialization presented above). With
this modification the trivialization does meet the necessary
conditions and can be capped off with $R_{ne}$, $R_{no}$, $R_{pe}$,
and $R_{po}$.
\subsubsection{Two bundles}
Consider the space $X_{m+1}$ of holomorphic immersions $D_{m+1}\to\C$ with $(m+1)$ convex corners on
the boundary. One called positive and the rest called negative.
\begin{lma}
The weak homotopy type of $X_{m+1}$ equals that of $SO(2)$. (That is $X_{m+1}$ is a
$K(\Z;1)$-space.)
\end{lma}
\begin{pf}
Each disk immersion with a marked point on its boundary can be contracted through immersed disks to
a small standard immersed disk near its marked point. In the present set up we use the positive
corner as the marked point and note that there is no problem keeping control of the rest of the
corners during such a deformation (parameterized by a compact space).
\end{pf}
We associate to each element $u\in X_{m+1}$ a Lagrangian boundary condition $\Lambda_{u}$ for the
$\bar\pa$-operator on $\C^2$-valued functions on $D_{m+1}$ by defining, with $\tau_\zeta$ denoting
the positive tangent vector of $\pa D_{m+1}$ at $\zeta,$
$$ \Lambda(\zeta)=\spa(du(\tau_\zeta),e^{i\theta(\zeta)})\subset\C^2,$$
where $-\theta_0\le\theta(\zeta)\le\theta_0$ and where $\theta=\theta_0$ ($\theta=-\theta_0$)
near the negative punctures along the incoming (outgoing) part of $\pa D_{m+1}$. Along the corresponding
parts of $\pa D_{m+1}$ near a positive puncture we let $\theta(\zeta)$ have the opposite signs.
Note that the $\bar\pa$-operator just mentioned has index $-(m-2)$ and that the dimension of its
kernel (cokernel if $m<2$) equals $0$ over each $u\in X_{m+1}$. Let $D\to X_{m+1}$ be the vector
bundle with fiber over $u$ equal to the cokernel (kernel) of this operator. Since the generator of
$\pi_1(X_{m+1})$ can be represented by a rigid rotation of a convex disk around its positive
puncture it is easy to see that $D$ is orientable.
We consider also the bundle $E\to X_{m+1}$, the fiber of which over
$u$ is the tangent space of conformal structures (automorphisms if
$m<2$) of the source $D_{m+1}$ of $u$. This is a fiber bundle of
fiber dimension $m-2$. Moreover, the natural linearization map
$i\circ \pa u\circ \gamma$, where $\gamma\in\End(T D_{m+1})$ is a
variation of the conformal structure gives a fiberwise isomorphism
$E\to D$ (this is a consequence of the general transversality
properties of the $\bar\pa$-equation in dimension $1$: rigid disk
are automatically transversely cut out).
The natural orientation of spaces of conformal structures induces an orientation on $E$. Assume that
$D$ is oriented. Then orientations of $D$ and $E$ either agrees or disagrees over every $u\in
X_{m+1}$.
\subsubsection{Diagram orientations and diagram trivializations}
We consider special types of trivializations of $\Lambda_u.$ These trivializations model
those coming from a knot diagram as mentioned above.
Let $\pa D_{m+1}=I_1\cup\dots\cup I_{m+1}$, be a subdivision into connected components where $I_1$
($I_{m+1}$) has its negative (positive) end at the positive puncture. For $u\in X_{m+1}$ we define a
{\em diagram orientation} of $u^\ast T\C$ over $\pa D_{m+1}$ to be a trivialization of this pull
back bundle which has the form
$$ t(\zeta)=\pm du(\tau_\zeta), \quad\text{ for }\zeta\in I_j, \text{ all }j.$$
Let $p$ be the positive corner of $u$ then we call $p$ {\em even (odd)} if a positive rotation of the incoming trivialization vector with
magnitude equal to the exterior angle of $u$ at $p$  gives the
negative (positive) outgoing trivialization vector. Let $q$ be a negative corner of $u$ then we call
$q$ {\em even (odd)} if a positive rotation of the incoming trivialization vector with magnitude equal to the exterior angle of $u$ at $q$
 gives the positive (negative) outgoing trivialization vector.
Note that a corner is even or odd according to whether or not the corresponding Reeb chord has even
or odd grading (see Lemma~\ref{eochords}), thus the number of odd corners of any diagram orientation
of a rigid disk in $X_{m+1}$ is odd.
We associate a {\em diagram trivialization} of the Lagrangian boundary condition $\Lambda_u$ to a
diagram orientation as follows. Fix a small neighborhood $U_j\subset I_j$ of the positive (negative)
end point of each $I_j$ which has its positive (negative) end point at a negative (positive) {\em
odd} puncture. The restriction of the diagram trivialization to the complement of the union of these
fixed neighborhoods is simply
$$ (t(\zeta),e^{i\theta(\zeta)}).$$
To complete the definition of the diagram trivialization we proceed as follows. On $U_j$
corresponding to a negative odd puncture we first make a positive $\pi$-rotation inside $\Lambda$
begining at
$$ (t(\zeta),e^{i\theta(\zeta)})$$
and ending at
$$ (-t(\zeta),-e^{i\theta(\zeta)}),$$
and then continue like that to the puncture. On $U_1$ corresponding
to a positive odd puncture we start out at the puncture with the
framing
$$ (-t(\zeta),-e^{i\theta(\zeta)})$$
make a negative $\pi$-rotation inside $\Lambda$ ending up at
$$ (t(\zeta),e^{i\theta(\zeta)}),$$
and then continue like that until we get into the region where the trivialization was already
defined.
\subsubsection{Inducing orientations on $D$}
Recall from Subsection~\ref{COCD} the capping disk $R_{no}$,
$R_{po}$, $R_{ne}$, $R_{pe}$. The determinant bundles over these
types of boundary conditions are in general {\em not} orientable.
Consider for example the split boundary condition $R_{ne}$ with
one-dimensional kernel and zero-dimensional cokernel. The kernel is
spanned by a function with second coordinate equal to $0$. Now,
applying a uniform $\pi$-rotation to the first or second line in the
split Lagrangian boundary condition brings us back to the original
boundary condition. Transporting the orientation along such a path
changes it when the first coordinate is rotated and does not change
it when the second coordinate is rotated. Thus the bundle is
non-orientable. However, the determinant bundles over subspaces of
the spaces of capping disks are orientable.
Let $Y_{no}$, $Y_{ne}$, $Y_{po}$, and $Y_{pe}$ be the spaces of
capping disks such that the second component of the trivialization
at the puncture equals $e^{\pm i\theta_0}$.
\begin{lma}
The determinant bundles of the $\bar\pa$-operator over $Y_{no}$,
$Y_{ne}$, $Y_{po}$, and $Y_{pe}$ are orientable.
\end{lma}
\begin{pf}
Note that all spaces above are homotopy equivalent to $SO(2)$ and
that the monodromy of a generating loop in $\pi_1(SO(2))$ preserves
orientation.
\end{pf}
Fix orientations on the bundles over $Y_{po}$ and $Y_{pe}$. This
determines orientations over $Y_{no}$ and $Y_{ne}$, respectively by
requiring that the orientations (of two glueable representatives)
glue to the canonical orientation of the determinant of the
resulting trivialized boundary condition over the closed disk.
\begin{lma}
Diagram orientaions of $\Lambda_u$ which vary continuously with
$u\in X_{m+1}$ induce an orientation on $D$ by gluing capping disks
(with orientations as above) to the corners of $u$.
\end{lma}
\begin{pf}
This construction has already been discussed, see
Subsection~\ref{cappingrulegeq3}. The result of gluing capping disks
$R_{no}$, $R_{ne}$, $R_{po}$, and $R_{pe}$ to the corresponding
punctures is a trivialized boundary condition on the closed disk. In
the gluing sequence the operators of the capping disk have oriented
determinant bundles which together with the canonical orientation on
the closed disk induce an orientation on $\bar\pa_{\Lambda_u}$. This
construction is clearly continuous in $u$.
\end{pf}
\subsubsection{Basic orientations and sign computation}
To compute signs in contact homology we need to compare the
orientation on $D$ induced from the orientation on $E$ and the
orientation on $D$ induced by a diagram trivialization. To this end
we look at the details of basic orientations to derive the shading
rules.
We first list the basic orientation choices.
\begin{itemize}
\item First we choose an orientation of $\R^2$. This orientation will be fixed
throughout the discussion.
\item Second we choose an orientation on $\C$.
\item Third choose orientations on $Y_{po}$ and $Y_{pe}$.
\end{itemize}
Recall that the first two choices determines canonical orientations
on the determinant bundles over trivialized boundary conditions over
the closed disk. Then the third choices induced orientations on
$Y_{no}$ and $Y_{ne}$ respectively via gluing. With all capping
operators oriented we orient all the bundles $D$ over $X_{m+1}.$
We next describe how the choices of orientation of $\C$ changes
canonical orientations. Recall that the canonical orientation on the
determinant bundle over trivialized boundary conditions on
the $0$-punctured disk was defined by, after deformation, expressing
the operator as an operator over $D^2$ with constant $\R^2$ boundary
conditions and an operator over $\C P^1$ with complex kernel or
cokernel attached at the origin of $D^2$. The kernel and cokernel of
the latter gets oriented by viewing them as complex vector spaces.
Now if the orientation of $\C$ is switched then the orientation of
each odd dimensional complex vector space changes and the
orientation of each even dimensional complex vector space remains
the same. Since the dimension of the kernel and cokernel of the
operator over $D^2$ with constant $\R^2$ boundary conditions is
2 and 0, respectively, the effect on the canonical
orientations is as follows:
changing the orientation of $\C$ changes the canonical
orientation for each operator of index $4j$ and keeps the
orientation of each operator of index $4j+2$.
Recall the index of the operators $R_{ne},R_{no},R_{pe},R_{po}$ are
$1,0,-1,0,$ respectively.
Fix an orientation of $\C$. Consider $u\in X_1$ with one positive
odd corner and with a diagram orientation which agrees with that of
the boundary of the disk. Pick the orientation of the determinant
bundle over $Y_{po}$ such that when gluing $R_{po}$ to the
trivialized boundary condition just described the orientation
induced on $D$ agrees with the one induced on $D$ by the bundle $E$
of linearized conformal structures. Note that the operator obtained
by gluing $R_{po}$ has index $2$ and hence the same choice of
orientation works also for the other orientation of $\C$. Denote
this orientation $o_{po}$
\begin{lma}\label{lma1punct}
Let $u\in X_1$ be any disk with any diagram orientation. The
orientation induced on the fiber of $D$ over $u$ obtained by gluing
$R_{po}$ with orientation $o_{po}$ agrees with the orientation
induced by $E.$
\end{lma}
\begin{pf}
There are two possible diagram orientations for such a disk. For
diagram orientations which can be obtained continuously from the
diagram orientation given above the lemma is clear. To see that the
lemma holds also for other diagram trivializations we need only
check it for one. To this end we rotate the first coordinate of the
glued boundary condition by $\pi$. This transports the canonical
orientation to the canonical orientation and the orientation of the
capping disk to the right orientation. Moreover, it preserves the
orientation on the $2$-dimensional kernel of the $\bar\pa$-operator
on the one punctured disk. Thus the induced orientations on $u$ with
the two different diagram trivializations must be the same.
\end{pf}
Our main theorem (Theorem~\ref{thm:3Dsigns}) clearly follows from
the following lemma.
\begin{lma}\label{lmacombor}
Let $o_D$ be the orientation induced on $D$ from the diagram
trivialization corresponding to a diagram orientation and let $o_E$
be the orientation induced by $E$. For a choice of basic
orientations $o^A(\C)$ on $\C$, $o_{po}$ on $Y_{po}$, and $o^A_{pe}$
on $Y_{pe}$ we have
$$ o_D=(-1)^{s_A} o_E,$$
where $s_A$ is the number of A-shaded punctures in the diagram
orientation.
For the choice of basic orientations $o^B(\C)=-o^A(\C)$, $o_{po}$,
and $o^B_{pe}=o^A_{pe}$ we have
$$ o_D=(-1)^{s_B} o_E,$$
where $s_B$ is the number of B-shaded punctures in the diagram
orientation.
\end{lma}
\begin{pf}
We will prove this lemma by induction on the number of punctures in
the disk. As a ``base case'' we must understand the orientations on
one- and two-punctured disks. The lemma for the disk with one
positive odd puncture and no other punctures is simply a restatement
of Lemma \ref{lma1punct}.
Fix an orientation $\hat o$ on $\C$. To determine the sign of an
arbitrary diagram orientation we will use orientations of
$1$-dimensional moduli spaces and in particular
Lemma~\ref{lmaMorient}. More precisely, we consider holomorphic
disks with several convex and one concave corner. Exactly as above
one may use capping disks to orient the one-parameter family of
disks in which it lives. For such model families Lemma
\ref{lmaMorient} relating various orientations hold.
To determine the sign of the disk with one even positive corner and
one odd negative corner, consider the disk with one concave (even)
corner shown in Figure~\ref{fig:signs1}. It is easy to arrange a
one-parameter family
\begin{figure}[ht]
  \relabelbox \small {\epsfxsize=3.5in\centerline{\epsfbox{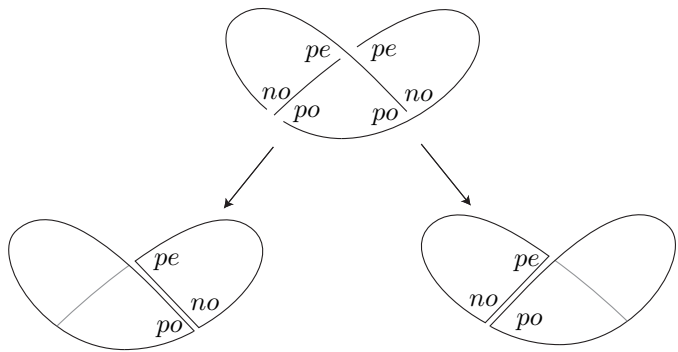}}}
  \relabel{1}{$pe$}
  \relabel{2}{$pe$}
  \relabel{3}{$pe$}
  \relabel{4}{$pe$}
  \relabel{5}{$po$}
  \relabel{6}{$po$}
  \relabel{7}{$po$}
  \relabel{8}{$po$}
  \relabel{9}{$no$}
  \relabel{10}{$no$}
  \relabel{a}{$no$}
  \relabel{12}{$no$}
  \endrelabelbox
        \caption{Splitting of a one dimensional space of disks.}
        \label{fig:signs1}
\end{figure}
which splits as follows. One splitting gives a one punctured disk
$D$ and a two punctured disk $B.$ The other splitting gives a one
punctured disk $D$ and a two punctured disk $B'.$ From the above we
know the signs on $D$ and $D'$ are both $+.$ Now
Lemma~\ref{lmaMorient} says we must have
\[\text{sign}(B)=-\text{sign}(B').\]
Note that $B$ and $B'$ are the two possible twice punctured disks
with a positive even puncture and a negative odd puncture. Thus we
know these disks have opposite signs. One choice of orientation
$o^A_{pe}$ yields the A-shading rule and the other one $o^B_{pe}$
the B-shading rule.
Note that the operator obtained by capping off $B$ has index $0$ and
hence its canonical orientation changes with $\hat o$. Since the
gluing of $R_{po}$ and $R_{no}$ also gives an operator of index $0$
also $o_{no}$ changes with the orientation of $\C$ and thus
$o^A_{pe}$ and $o^B_{pe}$ are independent of $\hat o$.

We have chosen orientations on the determinate line bundles over
$Y_{pe}$ and $Y_{po}$ so that the lemma is true for the cases
considered so far. To finish the ``base case'' of our induction we
are left to check that the lemma holds for disks $B''$ and $B'''$
with one positive odd puncture and one negative even puncture.
Arguing as in the above paragraph we see that the two configurations
of such disks must have opposite signs. Moreover, with respect to
one choice of signs on these disks the A-shading rule in the lemma
is correct with respect to the other the B-shading rule is. We
finish the base case by checking that a change of $\hat o$ changes
the signs of $B''$ and $B'''$. Note that the operator obtained by
capping off $B''$ has index index $2$. Thus its canonical
orientation is not affected by $\hat o$. However, the operator
obtained by gluing $R_{pe}$ and $R_{ne}$ has index $0$ so the
orientation $o^A_{ne}$ on $Y_{ne}$ induced from $o^A_{pe}$ changes
with $\hat o$. It follows that there is a choice $\hat o=o^A(\C)$ so
that the A-shading rule holds in all base cases. Moreover, with
$o^B(\C)=-o^A(\C)$ we find since $o^B_{pe}=-o^A_{pe}$ that the
induced orientation $o^B_{ne}$ on $Y_{ne}$ satisfies
$o^B_{ne}=o^A_{ne}$ and that the B-shading rule applies in all base
cases for $\hat o=o^B(\C)$.
To finish the proof we use the same argument in both cases. For
simplicity we give it only in the case of A-shading. Assume by
induction that the lemma is true for all disks with $n$ punctures.
Let $u$ be a disk with $n+1$ punctures.
If the puncture counterclockwise of the positive puncture is odd
then we can glue a once punctured disk $D$ to it. One of the four
possible cases if shown in Figure~\ref{fig:signs2}.
\begin{figure}[ht]
  \relabelbox \small {\epsfxsize=3.5in\centerline{\epsfbox{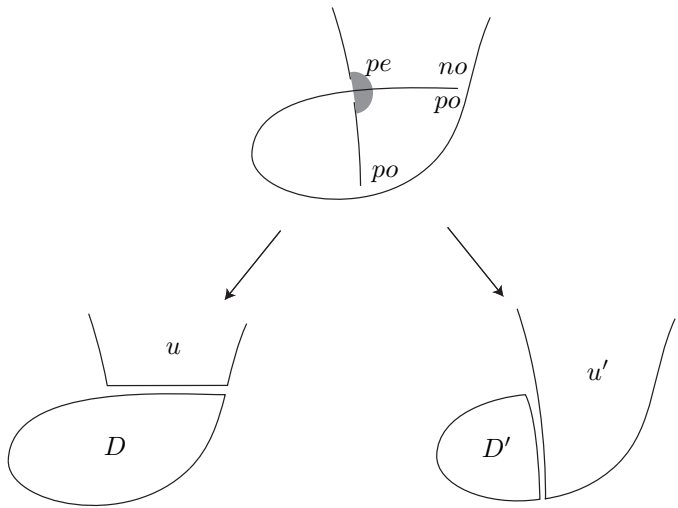}}}
  \relabel{1}{$po$}
  \relabel{2}{$po$}
  \relabel{no}{$no$}
  \relabel{pe}{$pe$}
  \relabel{D}{$D$}
  \relabel{D'}{$D'$}
  \relabel{u'}{$u'$}
  \relabel{u}{$u$}
  \endrelabelbox
        \caption{Inductive step with a negative corner odd. In the
        top drawing all corners are shaded.}
        \label{fig:signs2}
\end{figure}
We will finish the argument in this case. The remaining cases are
similar. The one dimensional moduli space formed by this gluing also
splits into an $n$ punctured disk $u'$ and a twice punctured disk
$D'.$ By induction we know the sign of $u'$ is $(-1)^{s'_A},$ where
$s'_A$ is the number of $A$-shaded regions in $u'.$ Moreover, the
sign of $D'$ is $+.$ Thus Lemma~\ref{lmaMorient} implies
\[\text{sign}(u)=\text{sign}(u)\text{sign}(D)=(-1)
\text{sign}(u')\text{sign}(D')=(-1)^{s'_A+1}.\]
Note the number of $A$-shaded regions of $u$ is $s'_A+1.$ Thus the
lemma holds for $u.$
An entirely analogous argument works if the puncture
counterclockwise of the positive puncture of $u$ is even. The only
difference is one must glue a twice punctured disk to $u.$ (Note in
this case one must actually have determined the signs on thrice
punctured disks, but this may be done as in the previous paragraph
by noting that at least one of the punctures must be odd.)\end{pf}

\section{Local Morse theory}\label{sec:lmt}
In this section we describe how to perturb highly degenerate
Legendrian submanifolds  into generic submanifolds. This will be our
main tool in constructing examples in the next section. In
Subsection~\ref{PDLS} we perturb a Legendrian submanifold, whose
complex projection has double points along the interior of a compact
codimension zero submanifold with boundary and a certain behavior at
this boundary using a Morse function on the submanifold, into a
generic Legendrian submanifold. For the perturbed Legendrian there
will be a double point in the complex projection for each critical
point of the Morse function. Adapting a construction of Floer
\cite{Floer89a} and Pozniak \cite{Pozniak},
we show that the contact homology boundary map for
these double points is related to the Morse ``gradient flow''
boundary map. In Subsection~\ref{sec:signsfordegen} we compare the
signs in these two boundary maps.

\subsection{Perturbing Degenerate Legendrian Submanifolds}\label{PDLS}
Consider a Legendrian submanifold $L$ in $\R^{2n+1}$. Let $U$ be
an $n$-manifold with boundary $\pa U$ and let $\phi\colon E\subset
J^1(U)\to\R^{2n+1}$ be a contact embedding (also respecting contact
forms) such that $\Pi_\C\colon U\to\C^n$ is an embedding. Let $\pa
U\times [0,2]$ be a collar neighborhood of $\pa U$ in $U$, with $\pa
U$ corresponding to $\pa U\times\{2\}$ . Let $V=U\setminus (\pa
U\times[1,2]),$ $W=U\setminus(\pa U\times[\frac32,2]),$ and
$W'=U\setminus\pa (U)\times[\frac12,2]$.
Assume that $L\cap \phi(E)$ consists of two sheets $L_1$ and $L_2$
represented in coordinates of $J^1(U)$ as the $1$-jet extensions
$j^1(f_1)$ and $j^1(f_2)$ of two functions $f_1$ and $f_2,$
respectively. Assume further that $f_2(x)=f_1(x)+c, c$ a positive
constant, for $x\in \overline{V}$, that $df_2\not= df_1$ for points
in $U\setminus \overline{V},$ and that $f_2(p,t)-f_1(p,t)$ is
monotone in $t$ for $(p,t)\in \pa U\times(1,2]$. We will frequently
think of $U$ and $V$ as subsets in $L_1\subset L.$ We will consider
functions $h:U\to \R$ satisfying
\begin{enumerate}
\item $h$ is a Morse-Smale function on $W$
\item the support of $h$ contains $V$ and is contained in $W$
\item all critical points of $f_2+h-f_1$ are critical points of $h$ and occur in $W'$
\item $h$ is real analytic near its critical points.
\end{enumerate}
Define the Legendrian submanifold $L_h$ to be the one obtained from
$L$ by replacing $j^1(f_2)$ in the front projection of $L$ with
$j^1{f_2+h}.$ The double points of $\Pi_{\C}(L_h)$  which lie inside
$\Pi_\C(E)$ correspond to the critical points of $h$ and double
which lie outside $\Pi_\C(E)$ corresponds to double points of
$\Pi_\C(L\setminus (j^1(f_1)(U)\cup j^1(f_2)(U)).$ Denote the part
of $\Pi_{\C}(L_h)$ corresponding to $j^1(f_1)(U)$ by $U_{f_1}$ (and
similarly for $U_{f_2+h}, V_{f_1},$ etc.). Having identified the
double points of $\Pi_\C(L_h)$ we know the generators of its algebra
$\A(L_h)$.
We next consider holomorphic disks. Let $x$ and $y$ be distinct
critical points of $h,$ thought of as double points of
$\Pi_\C(L_h).$ Let $(s,t)\in \Theta = \R\times[0,1]$ be conformal
coordinates for $D_2$ and $N$ a small tubular neighborhood of
$U_{f_1}$ in $\C^n$ that contains $U_{f_2+h}.$  Define
\[\mathcal{N}_h(x,y)=\{u:(\Theta,\partial \Theta)\to N, U_{f_1}\cup U_{f_2+h}) |
\text{ satisfying (1)--(3)} \}\]
where
\begin{enumerate}
\item $\dbar u=0$
\item $\lim_{s\to -\infty} u(s,t)=x$ and $\lim_{s\to \infty} u(s,t)=y$ and
\item $u(\R\times\{0\})\subset U_{f_1}$ and $u(\R\times\{1\})\subset U_{f_2+h}$.
\end{enumerate}
To define the operator $\dbar$ we use the standard complex structure on $\C^n.$
\begin{thm}\label{thm:lmt}
With the notation established above, there exists an $\epsilon>0$
such that if the $C^2$ norm of $h$ is less than $\epsilon$, and $N$
is contained in an $\epsilon$ neighborhood of $U_{f_1}$, for
positive $\lambda$ sufficiently small there is a one to one
correspondence between $\mathcal{N}_{\lambda h}(x,y)$ and gradient
flow lines of $\lambda h$ connecting $x$ to $y.$
\end{thm}
Our strategy for proving Theorem~\ref{thm:lmt} is to transplant our
problem into the cotangent bundle of some manifold and then use a
slightly modified version of a construction of Floer. To this end we
recall Floer's construction \cite{Floer89a}. Given a manifold $K$
consider its cotangent bundle $X=T^*K$ thought of as a symplectic
manifold. We will denote the zero section of $X$ by $K_0.$ Fix any
almost complex structure $J$ on $X$ that agrees with the canonical
one along the zero section. Let $h$ be a Morse function on $K$ and
define $K_h\subset X$ to be the graph of $dh.$ Let $\pi:X\to K$ be
the projection map and set $H(x)=h(\pi(x)).$ The Hamiltonian $H$
generates a flow $\phi_t$ on $X$ and $\phi_1(K_0)=K_h.$ Define the
time dependent almost complex structure $J_t=\phi_t^* J$ and
\[\mathcal{N}'_h(x,y)=\{u(\Theta,\partial \Theta)\to X |
\text{ satisfying (1)--(3)} \}\]
where
\begin{enumerate}
\item[(1)] $\frac{\partial u(s,t)}{\partial s}+ J_t(u(s,t)) \frac{\partial u(s,t)}{\partial s}=0$
\item[(2)] $\lim_{s\to -\infty} u(s,t)=x$ and $\lim_{s\to \infty} u(s,t)=y$ and
\item[(3)] $u(\R\times\{0\})\subset K_0$ and $u(\R\times\{1\})\subset K_h$.
\end{enumerate}
In \cite{Floer89a} it was shown that if the $C^2$ norm of $h$ is
sufficiently small then the map $u\mapsto u(s,0)$ from
$\mathcal{N}_h'(x,y)$ to $C^\infty(\R, K)$ is a bijection onto the
set of bounded trajectories of the gradient flow of $h$ connecting
$x$ and $y.$ Moreover if $h$ is a Morse-Smale function then the
moduli space $\mathcal{N}_h'(x,y)$ is transversely cut out by its
defining equation.
Define $\mathcal{N}_h''$ just as $\mathcal{N}_h'$ except instead of
the equation in condition (1) use the equation $\dbar_J u=0.$ The
modification of Floer's construction we need is given in the
following theorem.
\begin{lem}\label{lem:interp}
If the $C^2$ norm of $h$ is sufficiently small then there is a one
to one correspondence between $\mathcal{N}_{\lambda h}'(x,y)$ and
$\mathcal{N}_{\lambda h}''(x,y)$ for any $\lambda$ sufficiently
small.
\end{lem}
\begin{proof}
Let $\phi^\lambda_t$ be the flow of the Hamiltonian $\lambda
h(\pi(x))$ and $J^\lambda_t=(\phi^\lambda_t)^*J.$ In local Darboux
coordinates one may compute that
\[J^\lambda_t=J+\begin{pmatrix} -t\lambda d^2 h & 0\\ t^2\lambda^2 (d^2 h)^2 & t\lambda
d^2h\\\end{pmatrix}\]
where $d^2h$ is the Hessian of $h.$
Now consider the two parameter family of complex structures
\[J^{\lambda,s}_t=J+\begin{pmatrix} -ts\lambda d^2 h & 0\\ t^2s^2\lambda^2 (d^2 h)^2 & ts\lambda
d^2h\\\end{pmatrix}.\] Note $J^{\lambda, 0}_t=J$ (independent of
$\lambda$). So for a fixed $\lambda,$ $J^{\lambda,s}_t$ interpolates
between $J^{\lambda}_t$ and $J.$ Define $\mathcal{N}'_{\lambda h,
s}(x,y)$ just as $\mathcal{N}'_{\lambda h}(x,y)$ (including boundary
conditions) except use the complex structure $J^{\lambda,s}_t$ in
condition (1).
\begin{clm}
For $\lambda$ sufficiently small the spaces $\mathcal{N}'_{\lambda
h, s}(x,y)$ are transversely cut out for all $s\in[0,1].$
\end{clm}
This claim establishes the lemma since it will set up a one to one
correspondence between $\mathcal{N}_{\lambda h}'(x,y)=
\mathcal{N}'_{\lambda h, 1}(x,y)$ and $\mathcal{N}''_{\lambda
h}(x,y)=\mathcal{N}'_{\lambda h, 0}(x,y).$
To prove the claim note that since $h$ is Morse-Smale we know
$\mathcal{N}'_{h}(x,y)$ is transversely cut out by its defining
equation \cite{Floer89a}. Said another way, $\dbar_{J_t}$ is a
regular operator. If the claim is not true then we will find a
sequence of complex structures converging to $J_t$ that are not
regular. But this is a contradiction since by Gromov compactness
(and the upper semi-continuity of the dimension of the kernel) we
know there is an $\epsilon'$ neighborhood of $J_t$ in the space of
almost complex structures that contain only regular complex
structures.
Assume the claim is false. So for all $\lambda>0$ there is some $s$
so that $\mathcal{N}'_{\lambda h, s}(x,y)$ is not transversely cut
out. Thus there is always a solution $\widetilde{u}\in
\mathcal{N}'_{\lambda h, s}(x,y)$ for which $D_{\widetilde{u}}
\dbar_{J^{\lambda,s}_t}$ is not surjective. Since $\widetilde{u}$
satisfies
\[\dbar_{J^{\lambda,s}_t}\widetilde{u}= d \widetilde{u}+ J^{\lambda,s}_t\circ d\widetilde{u}\circ
j=0\] the map $u=\phi^{1-\lambda}_t (\widetilde{u})$ has boundary
conditions $K_0$ and $K_h$ and satisfies
\[du+ \widetilde{J}^{\lambda,s}_t \circ du \circ j =0.\]
In addition, since $\phi^{1-\lambda}_t$ intertwines the complex
structures $J^{\lambda,s}_t$ and $\widetilde{J}^{\lambda,s}_t$ we
see that $\widetilde{J}^{\lambda,s}_t$ is not regular. One may
compute that $\widetilde{J}^{\lambda,s}_t- J_t$ is of order
$\lambda$ so for $\lambda$ small enough
$\widetilde{J}^{\lambda,s}_t$ must be regular. Thus for $\lambda$
small enough $J^{\lambda, s}_t$ is regular for all $s\in[0,1].$
\end{proof}
To finish the proof of Theorem~\ref{thm:lmt} consider the manifold
$K=\overline{U}\cup_{\partial \overline{U}} -\overline{U}$ and its
cotangent space $X=T^*K.$ We can find a symplectomorphism $\psi$
from a neighborhood $N$ of $U_{f_1}$ in $\C^n$ to a neighborhood
$N'$ of $U \subset X,$ thought of as part of the zero section, that
takes $U_{f_1}$ to $U$ and sends the standard complex structure
along $U_{f_1}\subset \C^n$ to the canonical complex structure along
the zero section. (To ensure that $U_{f_2+h}$ still sits in $N$ one
merely needs to make sure the $C^1$ norm of $h$ is sufficiently
small.) Now let $J$ be the complex structure on $X$ that extends the
one pushed forward from $N$ and is standard along the entire zero
section. There is a function $g:K\to \R$ that is constant on
$V\subset K$ and satisfies $\Gamma_g\cap T^*U=\psi(U_{f_2}).$ We can
choose $U$ to be a small enough neighborhood of $\overline{V}$ so
that the function $g$ has $C^2$ norm small enough for Floer's
results to hold. We can further find a function $\widetilde{h}:K\to
\R$ such that $\Gamma_{g+\widetilde{h}}\cap T^*U=\psi(U_{f_2+h}).$
In addition there is some $\epsilon>0$ so that if
$||h||_{C^2}<\epsilon$ then $g+\widetilde{h}$ will also have $C^2$
norm small enough for Floer's results to hold. According to
Lemma~\ref{lem:interp} we can now find a $\lambda$ so that
$\mathcal{N}'_{\lambda(g+\widetilde{h})}(x,y)$ is in one to one correspondence
with $\mathcal{N}''_{\lambda(g+\widetilde{h})}(x,y).$ It is easy to see that all
the holomorphic curves in $\mathcal{N}''_{\lambda(g+\widetilde{h})}(x,y)$ lie in
$N'$ and correspond, via $\psi,$ to holomorphic curves in
$\mathcal{N}_{h}(x,y).$ This establishes Theorem~\ref{thm:lmt}.

\subsection{Comparison with the Morse-Witten complex}\label{sec:signsfordegen}
We briefly recall how to compute the homology of a manifold via a Morse function.
Given a Morse-Smale function $f:M\to \R$ we get a chain complex generated
by the critical points of $f$. The boundary map of this chain complex comes
from counting isolated gradient flow lines between critical points. We assign
a sign to each such flow line as follows. Let $x$ and $y$ be two critical points.
The Morse index of $x$ being one larger than the Morse index of $y.$ Let $U_x$ be
the unstable manifolds of $x$ (which is the set of points in $M$ that under the
gradient flow, flow to $x$ as time goes to $-\infty$). and $S_y$ be the stable
manifold of $y$ (which is the set of points that flow to $y$ as time goes to $\infty$).
The gradient flow lines connecting $x$ to $y$ are exactly $U_x\cap S_y.$
If we chose orientations on $U_x$ for all critical points $x$ and assume $M$ is oriented
then we get induced orientations on $S_x$ for all critical points $x.$ If such a choice has
been made then $U_x\cap S_y$ is an oriented manifolds. If it is a one dimensional manifold
then it also gets an orientation by $-\nabla f.$ Thus for isolated flow lines connecting
$x$ and $y$ we get a sign by comparing these orientations. The boundary map in
the chain complex comes from the signed count of these flow lines.
The homology of this complex agrees with the homology of the manifold.
For more details see \cite{SalamonNotes}.
Consider the set up in Section~\ref{sec:lmt}: a Legendrian with an open set $V$ of double points
in its Lagrangian projection and a Morse-Smale function $h$ on $V$ with which we perturb the Legendrian.
We now wish to compare the signed count of flow lines above with the signs of holomorphic
disks associated to the flow lines by Theorem~\ref{thm:lmt}. To this end choose a tree of gradient
flow lines of $h$ that connect all the critical points of $h$ (we are assuming $V$ is connected so this
can be done). Each of the moduli spaces of flow lines can be assigned an orientation using the
coherent orientations from Section~\ref{OMS}, by identifying it with the space of holomorphic curves
using Theorem~\ref{thm:lmt}. We may now choose orientations on the unstable manifolds of the critical
points so that the orientations induced by intersecting the stable and unstable manifolds agrees with
the coherent orientations on the flow lines in the chosen tree. This is possible since we chose a tree
of flow lines. Note that the orientations on all other moduli spaces of flow lines are determined by
the orientations on the spaces of flow lines in the chosen tree and gluing constructions.
Since the orientations on the moduli space of flow lines coming from
Section~\ref{OMS} and the procedure outlined in the previous paragraph are both ``coherent'' they are
both determined by the orientations on the flow lines in the tree and gluing constructions. We have
proved
\begin{thm}\label{thm:sfh}
Under the correspondence between holomorphic disks and gradient flow
lines set up in Theorem~\ref{thm:lmt} there is a choice on
orientation on the unstable manifolds of $h$ so that the Morse-Smale
orientations and orientations from Section~\ref{OMS} agree on all
moduli space of flow lines.
\end{thm}

\subsection{Examples}\label{sec:ex}
In this section we give examples showing that oriented
contact homology is a finer invariant than contact homology over
$\Z_2.$ To keep computations simple, we give examples showing that the
linearized contact
homology over $\Z$ will distinguish Legendrian links in the 1-jet
space of $T^n$ not distinguished by the full contact over $\Z_2.$
Similar constructions can be applied to construct Legendrian submanifolds in
$\R^{2n+1};$ however, the computation of the full contact homology
is more difficult. (It relies on a
Morse theoretic description of all holomorphic disks relevant to
contact homology, involving gradient flow trees with cusps.) We
therefore defer discussions of this to a forthcoming paper.

\begin{figure}[ht]
  \relabelbox \small {\epsfxsize=3.5in\centerline{\epsfbox{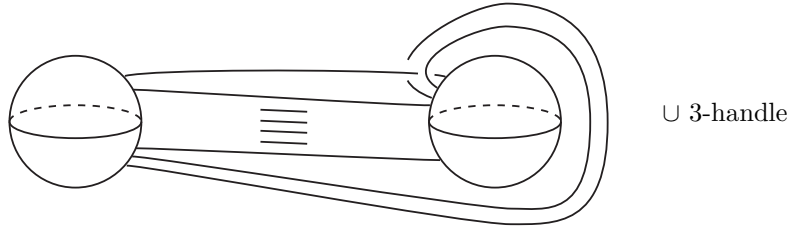}}}
  \relabel{A}{$\cup \text{ 3-handle}$}
  \endrelabelbox
        \caption{Handle decomposition of $D^4$. One of the 2-handles goes $p$ times over the
        1-handle the other goes once.}
        \label{fig:b4}
\end{figure}

Recall $J^1(T^n)=T^*T^n\oplus \R$ has a natural contact structure
$\alpha=\lambda +dt$ where $\lambda$ is the Liouville 1-from on
$T^*(T^n)$ pulled back to $J^1(T^n)$ and $t$ is the coordinate in
the $\R$ factor. Projecting out the $\R$ factor is called the
Lagrangian projection and is analogous to projecting out the
$z$-coordinate in $\R^{2n+1}.$
Note that since $T^*(T^n)$ admits a complex structure coming from
$\C^n$ by quotienting by a holomorphic action the analytic set up of
contact homology in \cite{ees1} for $\R^{2n+1}=J^1(\R^n)$ also works
for $J^1(T^n).$ Moreover, the proofs in \cite{ees1} carry over
word-for-work to this setting. (Note, in general, setting up contact
homology in jet spaces requires more work, see \cite{eesjet}.) With
the above understood we can take the contact homology of Legendrian
links in $J^1(T^n)$ to be well-defined.
Our main theorem of this subsection is
\begin{thm}
For any $n\geq 4$ there are infinitely many Legendrian links,
topologically isotopic to two copies of the zero section, in
$J^1(T^n)$ that are distinguished by their oriented linearized
contact homology, but have the same classical invariants and tame
isomorphic contact homology DGA's over $\Z_2.$
\end{thm}
Let $L_0$ be the zero section of $J^1(T^n)$ and $f$ a non-negative
self-indexing Morse function on $T^n$ with minimal number of
critical points. We will now alter $f$ near the unique index 0
critical point. For $p$ odd Figure~\ref{fig:b4}
gives a handle decomposition of $D^4.$
There are clearly analogous decompositions for all $n\geq 4$.
>From this we get a Morse function $h$ on
$D^n$ with critical points $c_0, c_1, c_2, c_2', c_3.$
The critical points are labeled by their Morse
index. We can assume that $h(c_0)<h(c_1)<h(c_2)<q<h(c_2')<h(c_3)$ and
$h^{-1}(\frac12)=\partial D^4.$ Moreover, the boundary map for the
Morse-Witten complex is
\begin{eqnarray*}
                \partial c_0=0& \quad \quad &\partial c_2'=c_1\\
                \partial c_1=0& \quad \quad &\partial c_3=pc_2'-c_2\\
                \partial c_2=pc_1& \quad \quad &
\end{eqnarray*}
Let $f_p$ be the Morse function on $T^n$ equal to $f$ on
$f^{-1}((-\infty, \frac12]),$ that is, outside a neighborhood of the
index 0 critical point of $f,$ and in the neighborhood equal to $h.$
We can assume $f_p$ is smooth. Let $L_p$ be the 1-jet of
$\epsilon(f_p-q)$ in $J^1(T^n):$
\[L_p=\{(x,\epsilon df_p(x), \epsilon(f_p(x)-q))\}.\]
We now fix the index 0 critical point and measure all gradings
relative to this. According to Theorem~\ref{thm:lmt} and the
discussion following it we know, for $\epsilon$ small enough, all
the holomorphic disks coming in to the computation of the boundary
map come form gradient flow lines. Thus one may easily compute the
boundary map of $L_0\cup L_p$ on the generators associated to $h$
are
\begin{eqnarray*}
              \partial c_0=0& \quad \quad &\partial c_2'=0\\
                \partial c_1=pc_2& \quad \quad &\partial c_3=pc_2'\\
                \partial c_2=0& \quad \quad &
\end{eqnarray*}
Moreover, the boundary map on the other generators comes from the
Morse boundary map for $f.$ Thus when everything is reduced mod 2
the boundary map is independent of $p$ but the contact homology over
$\Z$ (with relative grading chosen so that $c_0$ has grading $n$)
of the link $L_0\cup L_p,$ when $n\geq 8,$ is
\[ LCH_i(L(p))=
    \begin{cases}
      \oplus_{n+1} \Z, & i = 1 \\
      \oplus_{(n+1)n} \Z \oplus \Z_p & i=2\\
      \oplus_{(n+1)n(n-1)} \Z & i=3\\
      \oplus_{(n+1)n(n-1)(n-2)+1} \Z & i=4\\
      \oplus_{\binom{n+1}{i+1}} \Z & 4<i<n-3\\
      \oplus_{(n+1)n(n-1)+1} \Z & i=n-3\\
      \oplus_{(n+1)n} \Z \oplus \Z_p& i=n-2\\
      \oplus_{n+1}\Z & i=n-1\\
      \oplus_2 \Z & i=n.
    \end{cases}\]
One obtains a similar answer when $4\leq n<8.$ Note these gradings
are only relative gradings, but this is irrelevant as the $p$
torsion distinguishes the Legendrian links independent of grading.
The gradings were computed using the Morse index, see \cite{ees1}.
\section{Double points of exact Lagrangian immersions}\label{sec:dbpt}
\subsection{Double point estimates for Legendrian submanifolds with good algebras}
We describe how to use contact homology to derive a lower
bound on the number of double points of an exact Lagrangian
immersion. Contact homology has also been used to study
intersections of a pair of immersions \cite{Akaho}.
Let $f\colon L\to\C^n$ be an exact Lagrangian immersion. Then after
small perturbation we may assume that the Legendrian lift of $L$ is
an embedding which is chord generic. Let $(\A,\pa)$ denote its DGA.
Recall a DGA is augmented if the differential of no generator
contains a constant.
An {\em augmentation} of an algebra is a graded map $\epsilon:\A\to\Z$ such
that $\epsilon(1)=1$ and $\epsilon\circ\partial=0.$ Given an augmentation
$\epsilon$ the graded algebra tame isomorphism $\phi_\epsilon(a)=a+\epsilon(a)$
will conjugate $(\A,\pa)$ to an augmented algebra.
A DGA is called {\em good} if it admits an augmentation, and is hence tame
isomorphic to an augmented DGA. We show the following theorem where
we use the algebra $\A(L)$ with coefficients other than $\Z$, see
Remark \ref{rmkcoeff}.
\begin{thm}\label{dpegood}
Let $f\colon L\to \C^n$ be an exact Lagrangian immersions and let $(\A,\pa)$
be the DGA associated to an embedded chord generic Legendrian lift of $f$.
If $(\A,\pa)$ is good then $f$ has at least
$$
\frac12\dim(H_\ast(L;\Lambda))
$$
double points, where $\Lambda=\Q$ or $\Lambda=\Z_p$ for any prime $p$
if $L$ is spin and $\Lambda=\Z_2$ if $L$ is not spin.
\end{thm}
\begin{pf}
To simplify notation we identify $L$ with its image under the
embedding which is the lift of $f$ and write
$L\subset\C^n\times\R$. Let $L'$ be a copy of $L$ shifted a large
distance in the $z$-direction, where as usual $z$ is a coordinate in
the $\R$-factor. Then $L\cup L'$ is a Legendrian link. Moreover,
assuming that the shifting distance in the $z$-direction is
sufficiently large, shifting $L'$ $s$ units in the $x_1$-direction
gives a Legendrian isotopy of $L\cup L_s'$. After a large such shift
$L\cup L_s'$ projects to two distant copies of $\Pi_\C(L)$ and it is
evident that an augmentation for $L$ gives an augmentation for
$L\cup L'_s$. Moreover, the linearized contact homology of $L\cup
L'$ equals the set of sums of two vector spaces from the linearized
contact homology of $L$.
We will next compute the linearized contact homology of $L\cup L'$ in
a different manner. Let $g\colon L\to\R$ be a Morse function on
$L$ and use $g$ to perturb $L$ in $U\subset J^1(L)$, where $U$ is a
small neighborhood of the $0$-section. After identification of $U$
with a neighborhood of $L'$ in $\R\times\C^n$ we use this isotopy to
move $L'$ to $L''$. The projection of $L''$ into $\C^n$ then agrees (locally)
with an exact deformation of $L$ in its cotangent bundle and there is
a symplectic map from that cotangent bundle to a neighborhood of
$\Pi_\C(L)$ in $\C^n$. Pulling back the complex structure from $\C^n$
we get an almost complex structure on $T^\ast L$.
The intersection points of $L$ and $L''$ are of three types.
\begin{itemize}
\item[(1)] Critical points of $g$.
\item[(2)] Pairs of intersection points between $L$ and $L''$ near the
self-intersections of $L$.
\item[(3)] Self intersection points of $L$ and of $L''$ near self
intersections of $L$.
\end{itemize}
Fix augmentations of $\A(L)$ and of $\A(L'')$. If $\pa$ is the
differential of $\A(L\cup L'')$ it is easy to see that any monomial in
$\pa c$, where $c$ is a Reeb chord of type (1) or (2) must contain an
odd number of Reeb chords of type (1) and (2). Therefore the
augmentations of $\A(L)$ and $\A(L'')$ give an augmentation for
$\A(L\cup L'')$ that is trivial on double points of type (1) and (2).
Denote by $d$ the linearized differential induced by the
augmentations chosen and by $E_i$ the span of the double points of
type ($i$), $i=1,2,3.$ Suppose $a$ is a type (3) double point, then
$\partial a$ has no constant part and its linear part has no double
points of type (1) or (2), since each holomorphic disk with a
positive puncture at $a$ must have an even number of negative
corners of type (1) or (2). Thus $d(E_3)\subset E_3.$ If $b$ is of
type (1) or (2) then the linear part of $\partial b$ involves only
double points of type (1) and (2). Denote by $\pi_i$ the projection
onto $E_i, i=1,2,$ and $d_i=\pi_i\circ d.$ Then $d=d_1+ d_2$ on
$E_1\oplus E_2.$
Consider $d_1\colon E_1\to E_1$. We claim that for sufficiently small
perturbation $g$, $d_1\circ d_1|E_1=0$. To show this we consider gluing of two
(two-punctured) disks contributing to $d_1$. This gives a $1$-parameter
family two-punctured disks. Now, for
sufficiently small perturbation, no Reeb chord of
type (2) has length lying between the lengths of two Reeb chords of
type (1). Moreover, every Reeb-chord of type (3) has length bigger
than the difference of the lengths of two Reeb chords of type
(1). This shows that the $1$-parameter family must end at another
pair of broken disks with corners of type (1). It follows that $d_1^2=0$.
It follows that $d_1|E_1$ agrees with the Floer differential of $\hat L\cup
{\hat L}_g$, where $\hat L\subset T^\ast L$ is the $0$-section and
where ${\hat L}_g\subset T^\ast L$ is the graph of $dg$. Hence,
$$
\krn(d_1|E_1)/\img(d_1|E_1)\approx H_\ast(L;\Lambda).
$$
Write $E_1=W\oplus V$, where $W=\krn d_1|E_1$ and let $W'$ be a direct
complement of $d_1(V)\subset W$. Then $\dim W'=\dim H_\ast(L;\Lambda)$.
Fix the augmentations for $L$ and $L'$ which gives the element of
the linearized contact homology of $L$ which has the largest
dimension. By the above discussion we find that
$\krn(d_3)/\img(d_3)$ equals a direct sum of two copies of this
maximal dimension vector space. It follows that the contribution to
the linearized contact homology involving double points between
$L''$ and $L$ must vanish. We check how double points of type (2)
kill off the double points of type (1) that exist in the homology of
$(E_1, d_1).$
We compute
\begin{align*}
0&=d(d(W'))=\pi_1(d(d(W')))\\
&=\pi_1(d(d_2(W'))) = d_1(d_2(W')),
\end{align*}
where the third equality is due to the fact that $W'\subset E_1$ is in $\krn(d_1|E_1).$
It follows that
$\img(d_2|W')\subset\krn(d_1|E_2).$
Moreover, notice an element $e$ in $W'$ is
a non zero element in the linearized contact homology if and only if
$d_2 e=0$ and $e\notin\img(d_1|E_2).$
Thus if $d_2 e=0$ then $e$ is in $\img(d_1|E_2),$ showing that
$\krn d_2|W'\subset \img d_1|E_2.$
We find
\begin{align*}
\dim(E_2)&=\dim(\krn d_1|E_2)+\dim(\img d_1|E_2)\ge\\
&\dim(\img d_2|W')+\dim(\krn d_2|W')=\dim(W'),
\end{align*}
and conclude that
$$
2\cdot\sharp\{\text{double points}\}=\dim(E_2)\ge
\dim(W')=\dim(H_\ast(L;\Lambda)).
$$
\end{pf}
\subsection{Improving double point estimates}
In this section we show how to remove a constant from a double
point estimate.
\begin{thm}\label{improve}
Suppose there is a constant $K$ such that for any exact Lagrangian immersion
$f\colon L\to \C^n,$ $f$ has at least
$$
\frac12\dim(H_\ast(L;\Lambda))-K
$$
double points, where $\Lambda=\Q$ or $\Lambda=\Z_p$ for any prime $p$
if $L$ is spin and $\Lambda=\Z_2$ if $L$ is not spin. Then $f$ has
at least
$$
\frac12\dim(H_\ast(L;\Lambda))
$$
double points.
\end{thm}
We first prove a simple lemma that is a generalization of a ``spinning''
 operation in \cite{ees1}.
\begin{lma}\label{lmamultsph}
Let $f\colon L\to \C^n\times\R$ be a chord generic Legendrian
embedding with $R(f)$ Reeb chords and Maslov number $m_f$. Then, for
any $k\ge 1$ there exists a Legendrian embedding $F_k\colon L\times
S^k\to\C^{n+k}\times\R$ with $2R(f)$ Reeb chords and $m_{F_k}=m_f$.
\end{lma}
We call the Legendrian embedding of $L\times S^k$ the $k$-spin of $L.$
\begin{pf}
For $q\in L$, let $f(q)=(x(q),y(q),z(q))$. Note that translations in
the $x_j$-direction, $j=1,\dots,n$ and that the scalings $x\mapsto
kx$, $z\mapsto kz$, $k\ge 0$ are Legendrian isotopies which
preserves the number of Reeb chords. We may thus assume that $f(L)$
is contained in $\{(x,y,z)\colon |x|\le\epsilon\}$, where $\epsilon$
is very small. For convenience, we write $\R^n=\R\times\R^{n-1}$
with coordinates $x=(x_0,x_1)$ and corresponding coordinates
$(x_0,y_0,x_1,y_1)$ in $T^\ast\R^n=\C^n$.
Consider the embedding $S^k\subset\R^{k+1}\subset\R^{k+n}$, where
$S^k$ is the unit sphere in $\R^{k+1}$. Let
$(\sigma,x_0,x_1)\in S^k\times\R_+\times\R^{n-1}$
$$
(\sigma,x_0,x_1)\mapsto x_0\cdot\sigma+x_1,
$$
be polar coordinates on $\R^{k+n}$. Fix a Morse function $\phi$, with one
maximum and one minimum on $S^k$ which is an approximation of the
constant function with value $1$. Define $F\colon S^k\times
L\to\C^{n+k}\times\R$, $F=(F_x,F_y,F_z)$ as follows
\begin{align*}
F_x(\sigma,q)&=(1+x_0(q))\cdot\sigma + x_1(q),\\
F_y(\sigma,q)&=(1+x_0(q))^{-1}\nabla_{S^k}\phi(\sigma)
+\phi(\sigma)(y_0(q)\cdot\sigma+y_1(q)),\\
F_z(\sigma,q)&=\phi(\sigma)z(q),
\end{align*}
where we think of the gradient $\nabla_{S^k}\phi(\sigma)$ as a
vector in $\R^{k+1}$ tangent to $S^k$ at $\sigma$. It is then easily
verified that $F$ is a Legendrian embedding. Moreover, the Reeb
chords of $F$ occur between points $(q,\sigma)$ and $(q',\sigma')$
such that $\sigma=\sigma'$, $x_j(q)=x_j(q')$, $y_j(q)=y_j(q')$,
$j=0,1$, and either $z(q)=z(q')$ or $\nabla_{S^k}\phi(\sigma)=0$.
However, these conditions are incompatible with $f$ being an
embedding unless $\nabla_{S^k}\phi(\sigma)=0$ and we conclude that
the number of double points of $F$ are as claimed. The statement
about the Maslov number is straightforward.
\end{pf}
\begin{proof}[Proof of Theorem~\ref{improve}]
Given $K$ in the statement of the theorem, choose $l$ so that
$2^l>K.$ For any immersed exact Lagrangian $f\colon L\to \C^n$ lift
$f$ to an embedded Legendrian in $\C^n\times \R$ and $k$-spin this
Legendrian $l$ times. The Lagrangian projection of the resulting
Legendrian gives a new exact Lagrangian immersion $F_k:L\times
S^k\times\ldots \times S^k\to \C^{n+lk}.$ Since Reeb chords
correspond to double points $F_k$ has $2^l$ time as many double
points as $f.$ We have
\begin{align*}
2^l R(f)&=R(F_k)\geq \frac12 \dim(H_\ast(L\times S^k\times\ldots\times S^k;\Lambda))-K\\
&=\frac 12 (2^l \dim(H_\ast(L;\Lambda))) - K.
\end{align*}
Thus
\[
R(f)\geq \frac12 \dim(H_\ast(L;\Lambda)).
\]
\end{proof}


\end{document}